\documentclass[10pt]{article} %-*-LaTeX-*-
\usepackage[dvips]{graphicx}
\usepackage{psfrag}
\usepackage{calc}  %%required to make lists work correctly
\usepackage{amsmath}  
\usepackage{amsthm}
\usepackage{amsfonts}
\usepackage{amsopn}
\usepackage{amssymb}
\usepackage{amstext}
\usepackage[matrix,tips,graph,curve,web]{xy}

\makeatletter

\def\@maketitle{%
  \newpage
  \null
  \let \footnote \thanks
    {\normalfont\sffamily\bfseries\Large\noindent\@title \par}%
    \vskip 1em%
    {\normalfont\sffamily\large
        \noindent
        \@author
        \par}
  \par
  \vskip 4em}
\def\@seccntformat#1{\csname the#1\endcsname{.\ }}
\renewcommand\section{\@startsection {section}{1}{\z@}%
                                   {-3.0ex \@plus -1ex \@minus -.2ex}%
                                   {1.5ex \@plus.2ex}%
                                   {\normalfont\sffamily\large\bfseries}}
\renewcommand\subsection{\@startsection{subsection}{2}{\z@}%
                                     {-2.75ex\@plus -1ex \@minus -.2ex}%
                                     {1.5ex \@plus .2ex}%
                                   {\normalfont\sffamily\large}}
\def\fnum@figure{\normalfont\footnotesize\figurename~\thefigure}
\renewcommand\tableofcontents{%
    \section*{\contentsname
        \@mkboth{%
           \MakeUppercase\contentsname}{\MakeUppercase\contentsname}}%
    \@starttoc{toc}%
    }
\renewcommand*\l@part[2]{%
  \ifnum \c@tocdepth >-2\relax
    \addpenalty\@secpenalty
    \addvspace{2.25em \@plus\p@}%
    \begingroup
      \setlength\@tempdima{3em}%
      \parindent \z@ \rightskip \@pnumwidth
      \parfillskip -\@pnumwidth
      {\leavevmode
       \large \bfseries #1\hfil \hb@xt@\@pnumwidth{\hss #2}}\par
       \nobreak
       \if@compatibility
         \global\@nobreaktrue
         \everypar{\global\@nobreakfalse\everypar{}}%
      \fi
    \endgroup
  \fi}
\renewcommand*\l@section[2]{%
  \ifnum \c@tocdepth >\z@
    \addpenalty\@secpenalty
    \addvspace{1.0em \@plus\p@}%
    \setlength\@tempdima{1.5em}%
    \begingroup
      \parindent \z@ \rightskip \@pnumwidth
      \parfillskip -\@pnumwidth
      \leavevmode \sffamily\bfseries
      \advance\leftskip\@tempdima
      \hskip -\leftskip
      #1\nobreak\hfil \nobreak\hb@xt@\@pnumwidth{\hss #2}\par
    \endgroup
  \fi}
\renewcommand*\l@subsection{\sffamily\@dottedtocline{2}{1.5em}{2.3em}}
\renewcommand*\l@subsubsection{\@dottedtocline{3}{3.8em}{3.2em}}
\renewcommand*\l@paragraph{\@dottedtocline{4}{7.0em}{4.1em}}
\renewcommand*\l@subparagraph{\@dottedtocline{5}{10em}{5em}}
\makeatother

\makeatletter 
\@addtoreset{equation}{section}
\makeatother

\renewcommand{\theequation}{\thesection.\arabic{equation}}

\theoremstyle{plain}

\newtheorem{lemma}[equation]{Lemma}
\newtheorem{proposition}[equation]{Proposition}
\newtheorem{conjecture}[equation]{Conjecture}
\newtheorem{fact}[equation]{Fact}
\newtheorem{facts}[equation]{Facts}
\newtheorem*{theoremA}{Theorem A}
\newtheorem*{theoremB}{Theorem B}
\newtheorem*{theoremC}{Theorem C}
\newtheorem*{theoremD}{Theorem D}

\theoremstyle{definition}
\newtheorem{definition}[equation]{Definition}

\newtheorem{remark}[equation]{Remark}
\newtheorem{remarks}[equation]{Remarks}

\newtheorem{example}[equation]{Example}

\newtheorem{question}[equation]{Question}

\newenvironment{slist}[1]%
  {\begin{list}{}%
    {%
    \settowidth{\labelwidth}{#1}%
    \setlength{\itemindent}{0pt}%
    \setlength{\labelsep}{1em}%
    \setlength{\leftmargin}{\labelwidth+\parindent+\labelsep}%
    \setlength{\itemsep}{0pt}%
    \setlength{\parsep}{.6ex}}}%
  {\end{list}}

  {\begin{list}{}%
    {%
    \settowidth{\labelwidth}{#1}%
    \setlength{\itemindent}{0pt}%
    \setlength{\labelsep}{1em}%
    \setlength{\leftmargin}{\labelwidth+\labelsep}%
    \setlength{\itemsep}{0pt}%
    \setlength{\parsep}{.6ex}}}%
  {\end{list}}

\makeatletter
\newenvironment{subequations*}{% same thing but without incrementing
  \begingroup % conservative approach
  \let\protect\@nx
  \edef\@tempa{\def\@nx\theparentequation{\theequation}}%
  \@xp\endgroup\@tempa
  \setcounter{parentequation}{\value{equation}}%
  \setcounter{equation}{0}%
  \def\theequation{\theparentequation\alph{equation}}%
  \ignorespaces
}{%
  \setcounter{equation}{\value{parentequation}}%
  \global\@ignoretrue
}

\newcommand{\eval}[2][\right]{\relax
  \ifx#1\right\relax \left.\fi#2#1\rvert}
\makeatother

\newcommand{\hide}[1]{}
\newcommand{\sdef}[1]{{\bf #1}}

\newcommand{\xra}{\xrightarrow}
\newcommand{\aut}[1]{#1}

\def\Q{\mathbb{Q}}
\def\Z{\mathbb{Z}}
\def\R{\mathbb{R}}
\def\C{\mathbb{C}}

\def\P{\mathbb{P}}

\def\cp{\mathbb{CP}}

\newcommand{\roots}{\Delta}

\def\la#1{{\mathfrak{#1}}}

\def\t{\la t}
\def\v{\la v}

\def\td{\la t^*}

\def\inv{^{-1}{}}

\def\hra{\hookrightarrow}

\def\onto{\twoheadrightarrow}
\def\ddt0{{\left.{\frac{d}{dt}}\right|_{t=0}}}
\def\dds0{{\left.{\frac{d}{ds}}\right|_{s=0}}}

\def\CA{\mathcal{A}}

\def\CL{\mathcal{L}}

\def\CR{\mathcal{R}}
\def\CS{\mathcal{S}}

\def\om{\omega}

\def\al{\alpha}
\def\vp{\varphi}

\def\bd{\partial}

\def\sg{\sigma}
\def\bt{\beta}
\def\Gm{\Gamma}

\def\ep{\epsilon}

\def\lb{\lambda}

\def\d/{/\mspace{-6.0mu}/}

\def\<{\left\langle}
\def\>{\right\rangle}
\def\({\left(}
\def\){\right)}
\def\iso{\cong}

\def\c.{\cdot}
\def\h..{\ldots}

\def\o{\circ}
\def\emptyset{\varnothing}

\def\half{\frac{1}{2}}
\def\ot{\otimes}
\def\ti{\times}
\def\ul{\underline}

\DeclareMathOperator{\rank}{rk}
\DeclareMathOperator{\HH}{H}

\def\tr{{\operatorname{tr}}}

\def\vol{\operatorname{vol}}
\def\deg{\operatorname{deg}}

\def\Lie{\operatorname{Lie}}
\def\pt{\operatorname{pt}}

\def\hcf{\operatorname{hcf}}

\newcommand{\im}{\operatorname{im}}

\def\choose#1#2{{\binom{#1}{#2}}}

\newdir{ (}{{}*!/-1ex/@^{(}}
\newdir_{ (}{{}*!/-1ex/@_{(}}

\def\d/{/\mspace{-6.0mu}/{}}
\newcommand{\syq}[3]{#1\d/#3(#2)}
\newcommand{\ig}[2]{\int_{#2}\negthickspace\negthickspace{#1}}
\newcommand{\loc}{\lambda}

\newcommand{\singularities}{singularities}
\newcommand{\flag}{\Theta}
\newcommand{\U}[1]{U_{#1}}
\newcommand{\Uk}{\U{k}}
\newcommand{\CU}{\mathcal{U}}

\newcommand{\tU}{\tilde U}
\newcommand{\otno}{\twoheadleftarrow}
\newcommand{\hla}{\hookleftarrow}
\DeclareMathOperator{\opht}{ht}

\DeclareMathOperator{\data}{data}
\DeclareMathOperator{\diag}{diag}

\newcommand{\outn}{\nu_{\text{out}}}
\newcommand{\wtl}[1]{\widetilde{#1}}
\newcommand{\ol}[1]{\overline{#1}}
\newcommand{\f}{\footnotesize}

\newcommand{\spre}[1]{{\sffamily\footnotesize #1\par}\paragraph{}} %for `section preamble'
\newcommand{\union}{\cup}
\newcommand{\flagt}{\Theta}
\newcommand{\sbt}{\tau} %% subtorus of T 
\newcommand{\sbs}{\CS} %% subspace of \td, isomorphic to Lie(T/\sbt)^*
\theoremstyle{plain}
\newtheorem*{theoremBp}{Theorem B${}'$}

\newcommand{\cob}{\mathcal{U}}

\begin{document}

\title{Transversality theory, cobordisms,\\ 
and invariants of symplectic quotients}
\author{Shaun Martin%
\footnote{Institute for Advanced Study, Princeton, NJ;
smartin@ias.edu; February, 1999.}}
\maketitle

\section*{Introduction}

\subsection*{Symplectic quotients and their invariants}

This paper gives methods for understanding invariants of symplectic
quotients. The symplectic quotients that we consider are compact
symplectic manifolds (or more generally orbifolds), which arise as the
symplectic quotients of a symplectic manifold by a compact torus. A
companion paper~\cite{skm:g-t} examines symplectic quotients by a
nonabelian group, showing how to reduce to the maximal torus.

Throughout this paper we assume $X$ is a symplectic manifold, and that
a compact torus $T\iso S^1\ti\h..\ti S^1$ acts on $X$, preserving the
symplectic form, and having moment map $\mu:X\to\td$, where $\td$
denotes the dual of the Lie algebra of $T$. We assume that $\mu$ is a
proper map. (For definitions and our sign conventions see the notation
section at the end of this introduction). 

For every regular value $p\in\td$ of the moment map, the inverse image
$\mu\inv(p)$ is a compact submanifold of $X$ which is stable under
$T$, and on which the $T$-action is locally free (that is, every point
in $\mu\inv(p)$ has finite stabilizer subgroup). The \emph{symplectic
quotient}, which we denote $\syq Xp{T}$, is defined by taking the
topological quotient by $T$
\begin{equation*}
\syq Xp{T} := \frac{\mu\inv(p)}{T},
\end{equation*}
and is a compact orbifold (it is a manifold if the stabilizer
subgroup is the same for every point in $\mu\inv(p)$). Moreover the
symplectic form on $X$ defines in a natural way a symplectic form on
$\syq Xp{T}$. 

Many celebrated theorems in this field relate invariants of the triple
$(X,T,\mu)$ to invariants of the quotients $\syq Xp{T}$.  For example,
the \textbf{Duistermaat-Heckman theorem} \cite{dui-hec:var-coh}  relates a certain
oscillatory integral over $X$ to the volumes of the symplectic
quotients $\syq{X}{p}{T}$. Another example is the \textbf{Guillemin-Sternberg
quantization theorem}~\cite{gui-ste:geo-qua}, which relates the `geometric
quantization' of $X$ to that of its symplectic quotients%
\footnote{the geometric quantization is the index of a certain
naturally-defined Dirac operator; in the case of a K\"ahler manifold
this equals the space of holomorphic sections of a certain holomorphic
line bundle}. A third example is the
\textbf{Atiyah-Guillemin-Sternberg convexity theorem}, which relates a
very simple invariant of $(X,T,\mu)$, namely the convex hull of the
finite set of points $\mu(X^T)$, to an even simpler invariant of $\syq
Xp{T}$, namely whether it is empty.  One common feature of
these results is that the relevant invariants of $(X,T,\mu)$ can be calculated
in terms of data localized at the $T$-fixed points $X^T\subset X$.

\subsection*{The scope of this paper}

This paper provides results concerning a larger class of invariants,
including the integrals of arbitrary cohomology classes (thus
generalizing the volume) and the indexes of arbitrary elliptic
differential operators (generalizing the geometric quantization). In
order to describe this class of invariants, we first note that any
invariant of $\syq Xp{T}$ is also an invariant of the pair
$(\mu\inv(p),T)$ (the converse is of course not true). The easiest way to
describe the results of this paper is in terms of the submanifolds
$\mu\inv(p)$, for $p$ any regular value of $\mu$.

The submanifold $\mu\inv(p)$ defines an equivalence class
$[\mu\inv(p)]$, defined in terms of certain equivariant
cobordisms, and the invariants accessible by the methods of this paper
are those invariants that only depend on the
class $[\mu\inv(p)]$.

Explicitly, let $X'\subset X$ denote the subset consisting of those
points whose stabilizer subgroup is finite. Then the submanifold
$\mu\inv(p)\subset X'$ defines the cobordism class
\begin{equation*}
[\mu\inv(p)] \in \cob^*_T(X'),
\end{equation*}
where representatives of $\cob^*_T(X')$ are given by $T$-equivariant
maps of oriented manifolds to $X'$, and equivalences are given by the
boundaries of $T$-equivariant maps of oriented manifolds-with-boundary.
Explicitly, if $W\hra X'$ is any oriented manifold-with-boundary
mapped $T$-equivariantly to $X'$, then $[\bd W\hra X']=0 \in
\cob^*_T(X')$. Note that since $X'$ has a locally free $T$-action, every
manifold and cobordism must also have a locally free $T$-action.  

An example of an invariant that only depends on the class
$[\mu\inv(p)]$ is described in terms of the natural ring
homomorphism
\begin{equation*}
\kappa:\HH_T^*(X;\Q) \onto \HH_T^*(\syq{X}{p}{T};\Q)
\end{equation*}
defined by restriction, followed by the natural identification of the
equivariant cohomology of $\mu\inv(p)$ with the regular cohomology of
its quotient. This map is often referred to as the `Kirwan map', and
is known to be surjective~\cite{kir:coh-quo}. Given classes
$a,b\in\HH_T^*(X)$, then Stokes's theorem implies that 
the `cohomology pairing'
\begin{equation*}
\begin{split}
\HH_T^*(X)\ot\HH_T^*(X) &\to \Q \\
a,b &\mapsto \ig{\kappa(a)\smile\kappa(b)}{\syq XpT}
\end{split}
\end{equation*}
is an invariant of the equivalence class $[\mu\inv(p)]$ (Stokes's
theorem is also valid for orbifolds, as we explain in
appendix~\ref{app:orb}). A similar map exists in $K$-theory, and again
only depends on the class $[\mu\inv(p)]$

\subsection*{The main result of this paper}

We now describe the main topological result of this paper: theorem~C
(which appears in section~\ref{sec:f-p-c}).  Theorem~C describes a
cobordism between $\mu\inv(p)$ and a collection of submanifolds of $X$
that lie near the $T$-fixed points:

\begin{theoremC}[Approximate version]
Suppose the fixed point set $X^T$ is finite.  Then for every regular
value $p$ of the moment map,
\begin{equation*}
[\mu\inv(p)] = \sum_{i\in I} [\CS(F_i)];
\end{equation*}
where each $F_i\in X^T$ is a fixed point, and $\CS(F_i)$ is a $d$-fold
product of odd-dimensional spheres, lying in a small neighbourhood of
$F_i$, with $d=\dim T$. 
\end{theoremC}
In general, $\CS(F_i)$ is a $d$-fold fibre product of sphere bundles over
a connected component $F_i\subset X^T$ of the fixed point set.
Recall that, by definition, the equivalence class $\mu\inv(p)$ is
defined in terms of submanifolds on which $T$ has a locally free
action. The quotient $\CS(F_i)/T$ is an orbifold, and can be described
as a $d$-fold `tower' of weighted projective bundles over $F_i$.

By describing the submanifolds $\CS(F_i)$  explicitly, 
we can calculate the cohomology pairings
described above in terms of data localized at the fixed
points. Theorem~D carries this out, giving cohomological formulae in
terms of characteristic classes.

It is also possible, by applying techniques in $K$-theory, to derive
formulae for the indices of elliptic operators: these formulae will
appear in another paper.

\subsection*{Overview of the paper}

This paper has four main results, theorems A, B, C, and D. Their
logical relationship is as follows (the numbers indicate sections)

\begin{equation*}
\xymatrix{ & *\txt{\textit{Topology}} & *\txt{\textit{Cohomology}} \\
*\txt{\textit{Walls}} 
& *\txt{Theorem~A (\ref{sec:con-w-c-c}--\ref{sec:thma})\\wall-crossing-cobordism}
\ar@{=>}[r]\ar@{=>}[d]
& *\txt{Theorem~B (\ref{sec:w-c-f}--\ref{sec:w-c-f2})\\ wall-crossing formula} \ar@{=>}[d] \\
*\txt{\textit{Fixed points}} 
& *\txt{Theorem~C  (\ref{sec:rels}--\ref{sec:f-p-c})\\
fixed point cobordism}\ar@{=>}[r] 
& *\txt{Theorem~D (\ref{sec:f-p-f})\\fixed point formula.}}
\end{equation*}

\vspace{2.5ex}

Theorem~A is the main topological construction in this
paper. Theorems~A and~C each give a cobordism between
$\mu\inv(p)$ and a collection of `simpler' spaces: in theorem~C each
such space is a $d$-fold fibre product of sphere bundles over a
component of $X^T$, where $d=\dim T$; in theorem~A each such space is
a sphere bundle over a submanifold of a manifold $X^H$, where
$H\subset T$ is a $1$-dimensional subtorus. In fact $X^H$ is a
symplectic manifold, with an action of the $(d-1)$-torus $T/H$, and
having a moment map $\mu'$. The submanifold of $X^H$ which appears in
theorem~A is $\mu'{}\inv(q)$, for $q$ some regular value of $\mu'$. 
Theorem~A forms the inductive step in the proof of theorem~C, and
the induction is carried out in sections \ref{sec:rels} and \ref{sec:f-p-c}. 
The main techniques used in the proofs of theorems~A and C are
transversality theory, and general results in the theory of Lie group
actions on manifolds. The symplectic geometry which is used boils down
to a single fact, fact~\ref{fact:inf-geo-mom}, which is illustrated in
figure~\ref{fig:geo-mom}.

Theorems~B and D result from applying cohomological techniques to the
cobordisms constructed in theorems~A and C. Whereas a naive
application of Stokes's theorem would result in formulae which were
computable in principle, but unwieldy in practise, the real content of
theorems~B and D is to show how such formulae can be reduced to
computable formulae, eventually in terms of only the fixed points of
$X$. This is explained in more detail at the beginning of
section~\ref{sec:w-c-f}. In the proofs of theorems~B and D, fairly
extensive use is made of techniques in equivariant cohomology. We also
use various facts about orbifolds, which are explained in
appendix~\ref{app:orb}, as well as formulae which calculate integrals over the
fibres of weighted projective bundles. These formulae are proved in
appendix~\ref{app:wtd-proj-b}, and generalize classical formulae
involving Chern classes and Segre classes.

Finally, sections~\ref{sec:s2n} and \ref{sec:cp2n} calculate some
explicit examples. In section~\ref{sec:s2n} we study the $n$-fold product
of $2$-spheres $(S^2)^n$. This is a symplectic manifold, with a
Hamiltonian action of $SO(3)$, and the symplectic quotient
$\syq{(S^2)^n}{0}{SO(3)}$ is a manifold when $n$ is odd. 
These symplectic quotients have
been studied extensively, beginning with Kirwan's determination of the
Betti numbers~\cite{kir:coh-quo,kap-mil:sym-geo,hau-knu:coh-rin}. 
We use theorem~B, together with an
integration formula which allows us to reduce from a symplectic
quotient by $SO(3)$ to a symplectic quotient by the maximal torus
$S^1$ (proved in a companion paper~\cite{skm:g-t}) to give the
following formula for integrals of arbitrary cohomology classes on the
symplectic quotient $\syq{(S^2)^n}{0}{SO(3)}$, for $n$ odd:
\begin{equation*}
\ig{v_1^{l_1}\smile v_2^{l_2} \smile\h.. \smile v_n^{l_n}}{\syq{(S^2)^n}{0}{SO(3)}} = 
 -\half (-1)^{\frac{n-1}{2}}
 \sum_{\substack{K\subset\{1\h..n-1\} \\ |K| = \frac{n-1}{2}}}
  (-1)^{|K\cap\{1\h..m\}|}
\end{equation*}
where $\sum_i l_i = n-3$ and $m$ is equal to the number of odd $l_i$,
and $v_i$ is the natural degree $2$ cohomology class arising from the
$i$-th sphere in the product.

In section~\ref{sec:cp2n} we consider the space $(\cp^2)^n$. This has
a Hamiltonian action of $SU(3)$, and we calculate the volume of the
symplectic quotient $\syq{(\cp^2)^n}{0}{SU(3)}$ (the formula is not very
enlightening, but the methods are an application of theorem~D).

\subsection*{Relationship to other results}

There are a number of relationships between the cohomological formulae
proved in this paper (theorems~B and~D) and results of other authors. 

The mathematics in this paper was worked out in 1994, in Oxford and at
the Newton Institute in Cambridge. The intervening years have been
partly spent trying (possibly unsuccessfully) to understand how to
turn raw mathematics into a comprehensible manuscript. However, this
is a first attempt at writing mathematics, and so I beg the readers
indulgence in judging it.

The nonabelian localization formula of Witten~\cite{wit:two-dim-rev}
and Jeffrey-Kirwan~\cite{jef-kir:loc-non} gives an alternative way of
calculating cohomology pairings on symplectic quotients, involving
residues when $T=S^1$, and a multidimensional generalization of the
residue when $\dim T>1$.  An alternative approach to the
Witten-Jeffrey-Kirwan cohomology formula was taken by Guillemin and
Kalkman~\cite{gui-kal:jef-kir}, following from earlier independent
work of Kalkman~\cite{kal:coh-rin}. Guillemin and Kalkman use
`symplectic cutting' and `reduction in stages', but the geometric
arguments bear a strong resemblance to some of the arguments of this paper.

Jeffrey and Kirwan used the wall-crossing formula (theorem~B in this
paper, also the main result in
Guillemin-Kalkman~\cite{gui-kal:jef-kir}), together with results in
the companion paper~\cite{skm:g-t} to give a mathematically rigorous
proof of Witten's formulae for cohomology pairings on the moduli
spaces of stable holomorphic bundles over a Riemann surface described
above.

Some independent results on cobordisms of symplectic manifolds have
also been announced by Ginzburg, Guillemin and
Karshon~\cite{gin-gui-kar:cob-the}.

\subsection*{Acknowledgements}

I have benefited from very many elightening conversations with Simon
Donaldson, Mario Micallef, Frances Kirwan, Mike Alder, Michael
Callahan, Stuart Jarvis, Allen Knutson, Rebecca Goldin, Haynes Miller,
and Victor Guillemin.  But above all, I owe a great debt of
gratitude to Dietmar Salamon who has been both a friend and
source of inspiration to me.

\newpage
\tableofcontents

\subsection*{Notation and conventions}
Fixed through the entire paper, are the following:
\begin{slist}{$\syq{X}{p}{T}$}
\item[$X$] is a fixed smooth symplectic manifold (with symplectic form $\om$);
\item[$T$] $\iso S^1\ti\h..\ti S^1$ is a compact torus acting smoothly
on $X$, preserving $\om$;
\item[$\t,\td$] are the Lie algebra of $T$ and its dual, respectively;
\item[$\mu:X\to\td$] is a moment map for the $T$ action on $X$ (we
will assume throughout that $\mu$ is proper).
\end{slist}
We will use the following notational conventions:
\begin{slist}{$\syq{X}{p}{T}$}
\item[$\syq{X}{p}T$] $={\mu\inv(p)}/{T}$ denotes the `symplectic
quotient of $X$ by $T$ at $p$';
\item[$X^H$] denotes the subset of points fixed by the subgroup
$H\subset T$;
\item[$\HH^*(-)$] will always denote cohomology with rational coefficients;
\item[$\HH_G^*(-)$] denotes $G$-equivariant cohomology (rational
coefficients) for $G$ a group;
\item[$\kappa : \HH_T^*(X)\to\HH^*(\syq{X}{p}{T})$] for $p$ a regular
value of the moment map, denotes the natural map given by first
restricting to $\mu\inv(p)$, and then applying the natural isomorphism
$\HH_T^*(\mu\inv(p))\iso\HH^*(\syq{X}{p}{T})$ (the point $p$ will
always be clear from the context). $\kappa$ is often
referred to as the Kirwan map.
\end{slist}

\subsection*{Sign conventions for the moment map}

Different authors use varying sign conventions for the moment
map. Ours will be as follows. Given a symplectic manifold $(X,\om)$
with an action of a torus $T\iso S^1\ti\h..\ti S^1$ by
symplectomorphisms, let $V:\t\to\Gm(TX)$ be the infinitesimal action
map, taking an element $\xi$ of the Lie algebra of $T$ to the
corresponding vector field $V(\xi)$ on $X$. Then $\mu:X\to\td$ is a
\sdef{moment map} if it intertwines the $T$-action on $X$ and the
coadjoint action of $T$ on $\td$ (which is trivial in our case, since
$T$ is abelian), and which satisfies
\begin{equation}\label{eq:moment-map}
\<d\mu_x(v),\xi\> = \om_x(V(\xi),v),\qquad \forall x\in X, v\in T_xX,
 \xi\in\t.
\end{equation}
An {almost complex structure} $J:TX\to TX$ is \sdef{compatible}
with $\om$ if
\begin{equation}
g(\c.,\c.):=\om(\c.,J\c.)
\end{equation}
defines a Riemannian metric on $X$ (i.e.\ if $g$ is symmetric and
positive-definite). 

In the case of $S^1\subset\C^*$ acting on $\C$ by multiplication, 
our conventions boil down to the
following. Letting $z=x+iy$, and choosing the symplectic form
\begin{equation*}
\om=dx\wedge dy,
\end{equation*}
then the standard complex structure on $\C$ is compatible with $\om$,
and a moment map for the $S^1$-action is given by 
\begin{equation*}
\mu(z)=-\half|z|^2.
\end{equation*}

Finally, we recall the standard orientation of a complex vector
space, as defined in algebraic geometry:
if $\{e_1,\h..,e_n\}$ is a complex basis, then
\begin{equation}\label{eq:complex-ort}
\{e_1,ie_1,e_2,ie_2,\h.. e_n,ie_n\}
\end{equation}
is a real oriented basis. Thus, if $X$ is a symplectic manifold and
$J$ is a compatible almost complex structure, the orientation induced
by $J$ agrees with the orientation given by the top power of the
symplectic form.

\newpage
\section{Constructing the wall-crossing-cobordism}\label{sec:con-w-c-c}

This section contains the main construction of the paper: the
construction of the `wall-crossing-cobordism'.  The tools needed for
this construction comprise one fact from symplectic geometry, and some
transversality theory. We begin by stating the fact from symplectic
geometry, and illustrating it with a simple example. We then go on to
the main construction.

\subsection{Prelude: the geometry of the moment map}

We begin by explaining the key fact from symplectic geometry that we
use in this paper: this fact relates submanifolds defined by the group
action to submanifolds defined by critical points of the moment map.

\begin{figure}[!hbtp]\label{fig:geo-mom}
\psfrag{i}{\f$\mu\inv(Z)$}
\psfrag{j}{\f$\mu\inv(p_0)$}
\psfrag{k}{\f$\mu\inv(p_1)$}
\psfrag{a}{\f$\syq{X}{p_0}{T}$}
\psfrag{b}{\f$\syq{X}{p_1}{T}$}
\psfrag{p}{\f$p_0$}
\psfrag{s}{\f$p_1$}
\psfrag{q}{\f$q_0$}
\psfrag{r}{\f$q_1$}
\psfrag{m}{\f$\mu$}
\psfrag{x}{\f$X$}
\psfrag{z}{\f$Z$}
\psfrag{w}{\f$W$}
\psfrag{d}{\f$/T$}
\psfrag{c}{\f$W/T$}
\psfrag{0}{\f$X^{H_0}$}
\psfrag{1}{\f$X^{H_1}$}
\psfrag{2}{\f$X^{H_2}$}
\psfrag{3}{\f$X^{T}$}
\psfrag{4}{\f$\mu(X^{H_0})$}
\psfrag{5}{\f$\mu(X^{H_1})$}
\psfrag{6}{\f$\mu(X^{H_2})$}
\psfrag{7}{\f$\mu(X^{T})$}
\psfrag{8}{\f$\Lie(T/H_0)^*$}
\psfrag{9}{\f$\Lie(T/H_1)^*$}
\psfrag{l}{\f$\Lie(T/H_2)^*$}
\psfrag{t}{$\td$}
\begin{center}
\includegraphics[width=5.5in]{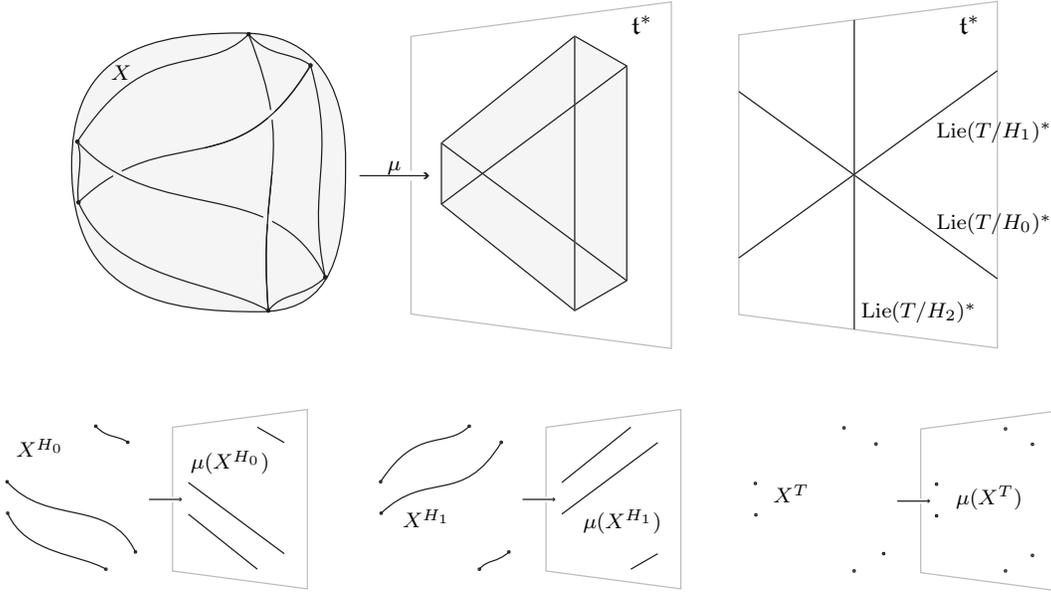}
\end{center}
\caption{A moment map, and its restriction to various submanifolds:
illustrating fact~\ref{fact:glo-geo-mom}. Here $T$ is $2$-dimensional,
and the subgroups $H_i$ are $1$-dimensional subtori.  The manifold $X$
and its submanifolds are only represented schematically: in the
concrete example from which this illustration is derived, $X$ is
$6$-dimensional, and each component of $X^{H_i}$ (represented by a
curved line in $X$) is a $2$-sphere (explained in
example~\ref{ex:geo-mom}).}
\end{figure}

Recall that $X$ denotes a symplectic manifold, acted on by
a torus $T\iso S^1\ti\h..\ti S^1$, with
associated moment map $\mu:X\to\td$, where $\td$
denotes the dual of the Lie algebra of $T$.  Let $\sbt\subset T$
be a subtorus. Then the short exact sequence of groups $\sbt\hra
T\onto T/\sbt$ induces the following exact sequences of Lie algebras
and their duals
\begin{equation*}
\begin{split}
\Lie(\sbt) &\hra \t \onto \Lie(T/\sbt) \\
\Lie(\sbt)^* &\otno \td \hla \Lie(T/\sbt)^*. \\
\end{split}
\end{equation*}
Hence for any subtorus $\sbt$ we will consider $\Lie(T/\sbt)^*$ to be
a subspace of $\td$ (of codimension $\dim \sbt$).

The key fact concerning the geometry of $\mu$ describes the way that
the derivative of $\mu$ encodes information about the $T$ action.  For
any point $x\in X$, letting $d\mu:T_xX\to\td$ denote the derivative,
we have
\begin{fact}[Infinitesimal version]\label{fact:inf-geo-mom} 
A subtorus $\sbt\subset T$ fixes $x$ if
and only if
\begin{equation*}
d\mu(T_xX)\subset\Lie(T/\sbt)^*.
\end{equation*}
\end{fact}

For example, if the action is locally free at $x$, then $d\mu_x$ must
be onto, and hence if $p\in\td$ is a regular value of $\mu$, then the
action on $T$ on $\mu\inv(p)$ is locally free.

The above fact has a global consequence. If $\sbt\subset T$ is a
subtorus, then we denote by $X^\sbt$ the set of points fixed by
$\sbt$: a local-coordinate argument shows that $X^\sbt$ is a closed
submanifold of $X$ (and an averaging argument shows that $X^\sbt$ is a
symplectic submanifold of $X$).

\begin{fact}[Global version]\label{fact:glo-geo-mom} The moment map 
$\mu$ maps each component of $X^\sbt$ to an
affine translate of $\Lie(T/\sbt)^*$ in  $\td$.
\end{fact}

For example, fixing a $1$-dimensional subtorus $H\iso S^1$ of $T$,
then $\mu$ maps each connected component of $X^H$ to an affine
hyperplane in $\td$, parallel to $\Lie(T/H)^*$. The images of such
submanifolds $X^H$, as $H$ varies through all $1$-dimensional subtori
of $T$, form `walls' which separate regions of regular values in
$\mu(X)$. At the other extreme, $\mu$ maps each connected component of
$X^T$ to a point in $\td$.

\begin{example}\label{ex:geo-mom}
Let $X$ be the set of $3\ti 3$ Hermitian matrices with eigenvalues
$0,1$ and $4$, and let $T\subset SU(3)$ be the maximal torus. Then $T$
acts on $X$ by conjugation, and a moment map for this action is given by
sending a matrix to its diagonal entries. Figure~\ref{fig:geo-mom}
illustrates some of the features of the moment map in this case (the
image of the moment map is accurate, but the illustration of $X$ is
schematic: $X$ is $6$-dimensional). The details in this illustration
are explained below.

We describe $X$ and $T$ explicitly as follows. 
Let $T$ be the diagonal matrices in $SU(3)$, that is,
$T=\{\diag(e^{i\theta_0},e^{i\theta_1},e^{i\theta_2})\mid
\theta_1+\theta_2+\theta_3=0\}$, and let $t\in T$ act on a matrix
$A\in X$ by $A\mapsto tAt\inv$. The map which takes $A\in X$ to its
diagonal entries $(a_{11},a_{22},a_{33})$ takes values in a
$2$-dimensional hyperplane in $\R^3$ (since $a_{11}+a_{22}+a_{33}=\tr
A = 5$), and this hyperplane can then be identified with $\td$ to give
a moment map for the $T$-action (the symplectic
form on $X$ is defined by identifying $X$ with a certain coadjoint
orbit%
\footnote{The map $A\mapsto iA$ identifies $X$ with an adjoint orbit
of $U(3)$; using an invariant inner product to identify
$\Lie(U(3))\iso\Lie(U(3))^*$ then identifies $X$ with a coadjoint orbit,
on which there is a natural symplectic form. A moment map for the
$T$-action is then given by the composition
$X\hra\Lie(U(3))^*\onto\Lie(T)^*$.}).

The set of $T$-fixed points in $X$ are the diagonal matrices: the
diagonal entries must be $0,1,4$ in some order, and so there are $6$ such
matrices. That is, $X^T$ consists of $6$ isolated points. These
points and their images under $\mu$ are depicted in the lower
right part of figure~\ref{fig:geo-mom}. The Atiyah-Guillemin-Sternberg
convexity theorem states that the image $\mu(X)$ equals the convex
hull of the image $\mu(X^T)$ of these points. Note that  the example
we are considering is atypical, because each point
of $\mu(X^T)$ defines a vertex of the polyhedron $\mu(X)$. In general,
not every point in $\mu(X^T)$ defines a vertex: some may map to the 
interior of $\mu(X)$.

Now consider the $1$-dimensional subtorus
$H_0:=\{\diag(e^{i\theta_0},e^{-i\theta_0/2},e^{-i\theta_0/2}\}\subset
T$. Then $H_0$ fixes the `block-diagonal' matrices of the form
\begin{equation*}
\begin{pmatrix}b&0&0\\0&*&*\\0&*&*\end{pmatrix}.
\end{equation*}
The entry $b$ must be one of the eigenvalues $0,1$ or $4$, and the
remaining $2\ti 2$ block has eigenvalues given by the other two. Thus
$X^{H_0}$ is made up of three components (each such component turns out
to be a $2$-sphere). A similar analysis holds for the subtori
$H_1:=\{\diag(e^{-i\theta_1/2},e^{i\theta_1},e^{-i\theta_1/2}\}$ and 
$H_2:=\{\diag(e^{-i\theta_2/2},e^{-i\theta_2/2},e^{i\theta_2}\}$.
There are infinitely many $1$-dimensional subtori of $T$: all the
others have as their fixed points only the points $X^T$.

In figure~\ref{fig:geo-mom} the subspaces $\Lie(T/H_i)\subset\td$ are
shown (here $0\le i\le 2$). Since each $H_i$ has dimension $1$, these
subspaces have \textit{codimension} $1$: they are hyperplanes. Each
submanifold $X^{H_i}$ has three components, each of which maps to an
affine translate of $\Lie(T/H_i)$ (shown for $i=0,1$, the picture for
$i=2$ is similar).
\end{example}

\subsection{The main lemma, and the resulting construction}

\begin{definition}\label{def:trans_path}
Let $p_0$ and $p_1$ be regular values of the moment map $\mu:X\to\td$.
A \sdef{transverse path} is a one-dimensional submanifold $Z\subset\td$,
with boundary $\{p_0,p_1\}$, such that $Z$ is transverse to $\mu$.
\end{definition}

It follows from transversality theory that $\mu\inv(Z)$ is a
submanifold of $X$, with boundary $\mu\inv(p_0)\sqcup\mu\inv(p_1)$ 
(the boundary of $Z$ is a submanifold of
$\td$ which is also transverse to $\mu$). The wall-crossing-cobordism,
which we define in \ref{def:w-c-c}, is constructed from the
submanifold $\mu\inv(Z)$. This construction 
is made possible by the following result.

\begin{proposition}\label{pr:muinvZ}
For any $x\in\mu\inv(Z)$,
the stabilizer subgroup of $x$ is either finite or $1$-dimensional.
If $H\subset T$ is any subgroup isomorphic to $S^1$ then the
submanifold $X^H$ of points fixed by $H$ is transverse to $\mu\inv(Z)$.
\end{proposition}

We prove proposition \ref{pr:muinvZ} below.  First, we use this result
to define the wall-crossing-cobordism:

\begin{figure}[!htbp]\label{fig:w-c-c}
\psfrag{i}{\f$\mu\inv(Z)$}
\psfrag{j}{\f$\mu\inv(p_0)$}
\psfrag{k}{\f$\mu\inv(p_1)$}
\psfrag{a}{\f$\syq{X}{p_0}{T}$}
\psfrag{b}{\f$\syq{X}{p_1}{T}$}
\psfrag{p}{\f$p_0$}
\psfrag{s}{\f$p_1$}
\psfrag{q}{\f$q_0$}
\psfrag{r}{\f$q_1$}
\psfrag{m}{\f$\mu$}
\psfrag{x}{\f$X$}
\psfrag{z}{\f$Z$}
\psfrag{w}{\f$W$}
\psfrag{d}{\f$/T$}
\psfrag{c}{\f$W/T$}
\psfrag{0}{\f$X^{H_0}$}
\psfrag{1}{\f$X^{H_1}$}
\psfrag{2}{\f$X^{H_2}$}
\psfrag{3}{\f$X^{T}$}
\psfrag{4}{\f$\mu(X^{H_0})$}
\psfrag{5}{\f$\mu(X^{H_1})$}
\psfrag{6}{\f$\mu(X^{H_2})$}
\psfrag{7}{\f$\mu(X^{T})$}
\psfrag{8}{\f$\Lie(T/H_0)^*$}
\psfrag{9}{\f$\Lie(T/H_1)^*$}
\psfrag{l}{\f$\Lie(T/H_2)^*$}
\psfrag{t}{$\td$}
\begin{center}
\includegraphics[width=5.5in]{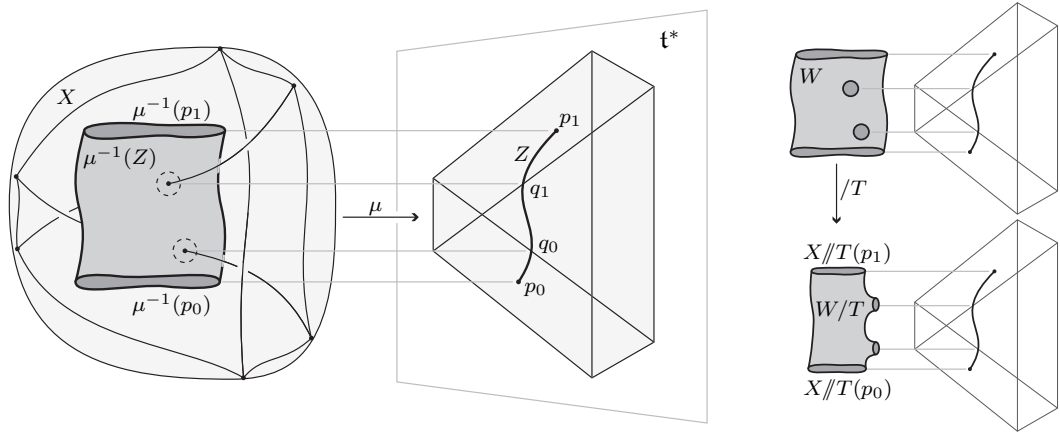}
\end{center}
\caption{A transverse path $Z$ and the resulting
wall-crossing-cobordism $W/T$. In the diagram on the left, the
submanifolds $X^{H_i}$ intersect $\mu\inv(Z)$, and the dashed circles
indicate open tubular neighbourhoods of these intersections: removing
these open neighbourhoods from $\mu\inv(Z)$ results in $W$.}
\end{figure}

\begin{definition}\label{def:w-c-c}
Let $W\subset X$ be the manifold-with-boundary given by removing open
subsets of $\mu\inv(Z)$ as follows. Fix a $T$-invariant metric on
$\mu\inv(Z)$, and set
\begin{equation*}
W := \mu\inv(Z) \setminus \bigsqcup_{H\iso S^1} N_\ep( X^{H})\cap\mu\inv(Z)
\end{equation*}
where $H$ runs through all $S^1$-subgroups of $T$, and $N_\ep(X^H)$ 
is the open $\ep$-tubular neighbourhood of $X^H$. We
choose an $\ep$ small enough to ensure that these subsets of
$\mu\inv(Z)$ have disjoint closures (it follows from proposition
\ref{pr:muinvZ} that this is possible). We define the
\sdef{wall-crossing-cobordism} to be the quotient
orbifold-with-boundary $W/T$. This procedure is illustrated in
figure~\ref{fig:w-c-c}.
\end{definition}

\begin{remarks}
\begin{enumerate}
\item
Only finitely many subgroups $H\iso S^1$ actually contribute in the
above definition. This is because $\mu\inv(Z)$ is a compact
$T$-manifold (since the moment map is assumed to be proper), and thus
only finitely many subgroups of $T$ can occur as stabilizer subgroups
\cite{bre:int-com,kaw:the-tra}.
\item
We may choose $p_0$ or $p_1$ outside the image of $\mu$, in which
case the corresponding boundary component will be empty. For
example, moving $p_1$ to lie outside the image of $\mu$ removes the
boundary component $\mu\inv(p_1)$, but introduces an extra
wall-crossing, like so:\\
{\psfrag{i}{\f$\mu\inv(Z)$}
\psfrag{j}{\f$\mu\inv(p_0)$}
\psfrag{k}{\f$\mu\inv(p_1)$}
\psfrag{a}{\f$\syq{X}{p_0}{T}$}
\psfrag{b}{\f$\syq{X}{p_1}{T}$}
\psfrag{p}{\f$p_0$}
\psfrag{s}{\f$p_1$}
\psfrag{q}{\f$q_0$}
\psfrag{r}{\f$q_1$}
\psfrag{m}{\f$\mu$}
\psfrag{x}{\f$X$}
\psfrag{z}{\f$Z$}
\psfrag{w}{\f$W$}
\psfrag{d}{\f$/T$}
\psfrag{c}{\f$W/T$}
\psfrag{0}{\f$X^{H_0}$}
\psfrag{1}{\f$X^{H_1}$}
\psfrag{2}{\f$X^{H_2}$}
\psfrag{3}{\f$X^{T}$}
\psfrag{4}{\f$\mu(X^{H_0})$}
\psfrag{5}{\f$\mu(X^{H_1})$}
\psfrag{6}{\f$\mu(X^{H_2})$}
\psfrag{7}{\f$\mu(X^{T})$}
\psfrag{8}{\f$\Lie(T/H_0)^*$}
\psfrag{9}{\f$\Lie(T/H_1)^*$}
\psfrag{l}{\f$\Lie(T/H_2)^*$}
\psfrag{t}{$\td$}
\includegraphics{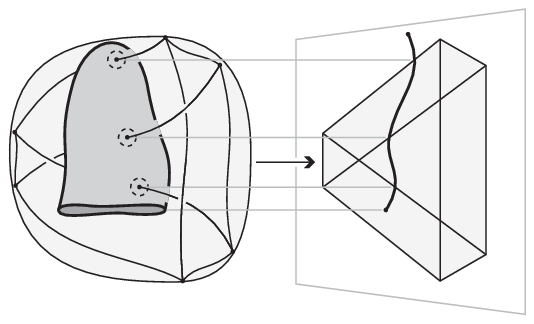}}
\end{enumerate}
\end{remarks}

\subsection{A combinatorial characterization of transverse paths, and
the proof of Proposition~\ref{pr:muinvZ}}

\begin{definition}\label{def:int}
We define a \sdef{wall} in $\td$ to be a connected component of the
image of $\mu(X^H)$, for some $H\iso S^1$. We define the \sdef{interior} of a
wall to be the set of points $q$ in the wall such that every point in
$\mu\inv(q)$ has stabilizer subgroup which is either $0$- or
$1$-dimensional.
\end{definition}

For example, in figure~\ref{fig:geo-mom} there are $9$ walls in total.
The arrangement of walls in $\td$ completely characterizes the set of
transverse paths:

\begin{lemma}[Geometry of $Z$ in $\td$]\label{lem:trans_path}
A path $Z$ is transverse to $\mu$ if and only if it intersects each wall
transversely in its interior.
\end{lemma}

\begin{proof}
We must show that, for every $x\in\mu\inv(Z)$, the tangent space
$T_{\mu(x)}\td$ is spanned by $d\mu(T_xX)$ and $T_{\mu(x)}Z$:
\begin{equation}\label{eq:span}
T_{\mu(x)}\td = d\mu(T_xX) + T_{\mu(x)}Z.
\end{equation}
We will
use the natural identification $T_{\mu(x)}\td\iso\td$. Let
$\sbt\subset T$ denote the subtorus given by the identity component of
the stabilizer subgroup of $x$: that is, $\sbt$ is the maximal
subtorus which fixes $x$. Then fact~\ref{fact:inf-geo-mom} implies
that $d\mu(T_xX)=\Lie(T/\sbt)^*$. Since $Z$ is $1$-dimensional, 
in order
for \eqref{eq:span} to hold $\sbt$
must be either $0$- or $1$-dimensional. This immediately implies that
every point of $Z$ must be either a regular value of $\mu$ or lie in
the interior of any wall which it is in.
If $\sbt$ is $0$-dimensional then
$d\mu(T_xX)$ already spans $\td$. If $\sbt$ is $1$-dimensional, then in
order for \eqref{eq:span} to hold, $T_{\mu(x)}Z$ must be complementary
to $\Lie(T/\sbt)^*$. Applying fact \ref{fact:inf-geo-mom}, 
this is the assertion that $Z$ is transverse to the
wall $\mu(X^\sbt)$ at $\mu(x)$.
\end{proof}

\begin{proof}[Proof of Proposition \ref{pr:muinvZ}]
In the course of proving lemma \ref{lem:trans_path}, we have already
seen that, for every point $x\in\mu\inv(Z)$, the stabilizer subgroup
of $x$ must be either $0$- or $1$-dimensional.

The statement that $Z$ is transverse to $\mu(X^H)$ (lemma
\ref{lem:trans_path}) 
is equivalent to the statement that the composition
\begin{equation}\label{eq:trans1}
T_xX^H \xra{d\mu} T_q\td \to \nu_q Z
\end{equation}
is surjective, for every $q$ in $Z\cap\mu(X^H)$, and for every
$x\in\mu\inv(q)\cap X^H$. 
Using the natural identification, via the pullback, of the normal
bundles:
\begin{equation*}
\mu^*:\nu Z \xra{\iso} \nu \mu\inv(Z),
\end{equation*}
then the composition \eqref{eq:trans1} can be factored
\begin{equation*}
T_xX^H \hra T_xX \to \nu_x \mu\inv(Z) \xra{\iso} \nu_q Z.
\end{equation*}
Since this map is surjective, it follows that the composition
$T_xX^H\to \nu_x\mu\inv(Z)$ is surjective, for every $q \in
Z\cap\mu(X^H)$, and for every $x\in\mu\inv(q)\cap X^H$, which gives
the result.
\end{proof}

\section{The data of a path, and how it describes the boundary of
the wall-crossing-cobordism}\label{sec:w-c-d}

\subsection{The data associated to a transverse path}

\begin{definition}\label{def:w-c-d}
Associated to each transverse path $Z\subset\td$ is a finite set
$\data(Z)$, which we refer to as the \sdef{wall-crossing data for
$Z$}.  We define $\data(Z)$ to be the set of pairs $(H,q)$, such that
$H\iso S^1$ is an oriented subgroup of $T$, and $q\in Z\cap\mu(X^H)$.
The orientation of $H$ is defined by the direction of the
wall-crossing: we orient $Z$ so that the positive direction goes from
$p_0$ to $p_1$; then a positive tangent vector in $T_qZ$, thought of
as an element of $\td$, defines a linear functional on $\t$, and this
restricts to a nonzero functional on $\la{h}$; and we orient $H$ to be
positive with respect to this functional.
\end{definition}

\begin{remarks}
\begin{enumerate}
\item We may also apply the above definition to a closed $1$-manifold
$Z\subset\td$, as long as $Z$ is oriented and transverse to $\mu$.
The wall-crossing data has a nontrivial interpretation in this case,
too.
\item We give an example to illustrate the orientation of $H$. 
Suppose our torus $T$ is the standard circle
$T=S^1=\R/\Z$, with Lie algebra and its dual identified with $\R$ in
the standard manner. In this case $p_0$ and $p_1$ are real numbers. If
$p_0<p_1$, then $Z$ must be the interval $[p_0,p_1]$, and each
wall-crossing induces the positive (i.e.\ standard) orientation on
$S^1$. If $p_1<p_0$, then $Z$ must be the interval $[p_1,p_0]$, and
each wall-crossing induces the negative orientation on $S^1$.)
\item
It is not possible for the same pair $(H,q)$ to appear twice in the
wall-crossing data, however we may have pairs $(H_0,q_0)$ and $(H_1,q_1)$
with $H_0=H_1$ while $q_0\ne q_1$: since $Z$ may cross the same wall more
than once; or $Z$ may cross different walls which are parallel and
thus correspond to the same subgroup.
And it is also possible for $q_0$ to equal $q_1$ (with
$H_0\ne H_1$). This is because a point $q$ may lie in the interior of two
different walls simultaneously. This happens when components of the
submanifolds $X^{H_0}$ and $X^{H_1}$ are disjoint in $X$, while their
images under $\mu$ both contain $q_0=q_1$. There are three points in
figure~\ref{fig:geo-mom} with this property.
\end{enumerate}
\end{remarks}

\subsection{The boundary of the wall-crossing-cobordism}

The wall-crossing data indexes the boundary components of $W$:

\begin{proposition}
\label{pr:bd_cpts}
The submanifold $W\subset X$ has boundary
\begin{equation*}
\quad \mu\inv(p_0) \quad\sqcup\quad \mu\inv(p_1) %
\quad \sqcup \quad \bigsqcup_{(H,q)\in\data(Z)} S_{(H,q)}
\end{equation*}
where
\begin{equation*}
S_{(H,q)} := \eval{S(\nu X^{H})}_{X^{H}\cap \mu\inv(q)}.
\end{equation*}
Here $S(\nu X^{H})$ denotes the unit sphere bundle in the normal
bundle of $X^{H}$ in $X$. Note that $S_{(H,q)}$ need not be connected: its
components correspond to the connected components of
$X^{H}\cap\mu\inv(q)$.
\end{proposition}

\begin{proof}
By proposition~\ref{pr:muinvZ}, each $X^{H}$ is transverse to
$\mu\inv(Z)$. Hence the
intersection $N_\ep(X^{H})\cap\mu\inv(Z)$ gives a tubular
neighbourhood of $X^{H}\cap\mu\inv(Z)$ in $\mu\inv(Z)$.
Similarly, the normal bundle to $X^{H}\cap\mu\inv(Z)$ in
$\mu\inv(Z)$ is the restriction of the normal bundle to $X^{H}$ in
$X$. By scaling, the unit sphere bundle is equivalent to the
$\ep$-sphere bundle.
\end{proof}

\begin{figure}[!hbtp]\label{fig:w-c-d}
\psfrag{j}{\f$\mu\inv(p_0)$}
\psfrag{k}{\f$\mu\inv(p_1)$}
\psfrag{a}{\f$\syq{X}{p_0}{T}$}
\psfrag{b}{\f$\syq{X}{p_1}{T}$}
\psfrag{p}{\f$p_0$}
\psfrag{s}{\f$p_1$}
\psfrag{q}{\f$q_0$}
\psfrag{r}{\f$q_1$}
\psfrag{m}{\f$\mu$}
\psfrag{x}{\f$X$}
\psfrag{z}{\f$Z$}
\psfrag{w}{\f$W$}
\psfrag{d}{\f$/T$}
\psfrag{c}{\f$W/T$}
\psfrag{0}{\f$S_{(H_0,q_0)}$}
\psfrag{1}{\f$S_{(H_1,q_1)}$}
\psfrag{2}{\f$P_{(H_0,q_0)}$}
\psfrag{3}{\f$P_{(H_1,q_1)}$}
\psfrag{8}{\f$\Lie(T/H_0)^*$}
\psfrag{9}{\f$\Lie(T/H_1)^*$}
\psfrag{l}{\f$\Lie(T/H_2)^*$}
\psfrag{t}{$\td$}
\begin{center}
\includegraphics[width=5.5in]{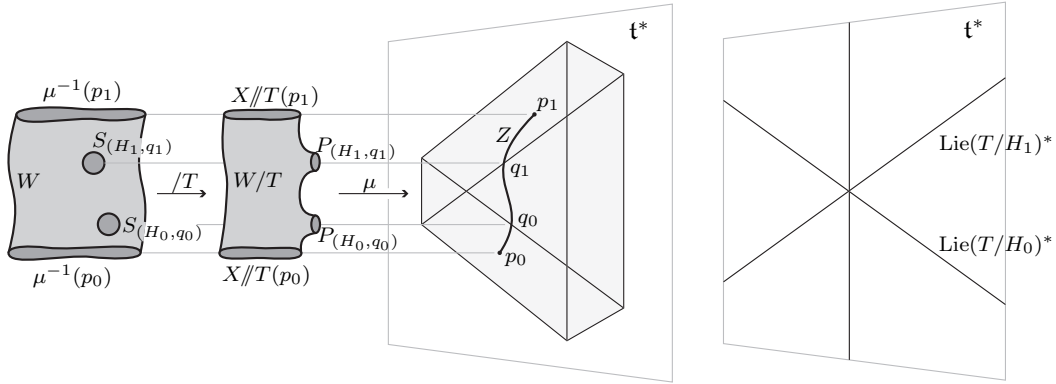}
\end{center}
\caption{A transverse path $Z$ with 
$\data(Z)=\{(H_0,q_0),(H_1,q_1)\}$, and the corresponding boundary
components of $W$ and of the wall-crossing-cobordism $W/T$.}
\end{figure}

Taking the quotient  by $T$ (which has a locally free action on $W$),
we thus have a description of the boundary of the
wall-crossing-cobordism:
\begin{equation}
\bd(W/T) \quad\iso\quad \syq{X}{p_0}T\quad\sqcup\quad \syq{X}{p_1}T%
 \quad\sqcup\quad %
 \bigsqcup_{(H,q)\in\data(Z)} P_{(H,q)},
\end{equation}
where
\begin{equation}\label{eq:PHq}
P_{(H,q)}:= S_{(H,q)}/T.
\end{equation}

\subsection{The boundary components as weighted projective bundles}

The rest of this section is devoted to giving a more explicit
description of the boundary components $P_{(H,q)}$.

The projection of the fibre bundle
\begin{equation}\label{eq:s-xh}
\pi:S_{(H,q)} \to X^{H}\cap\mu\inv(q)
\end{equation}
is a $T$-equivariant map, and by construction, $T$ acts with at most
finite stabilizers on the total space $S_{(H,q)}$. The subgroup $H$
acts trivially on the base, so that the $T$-action descends to an
action of $T/H$, and it follows from proposition \ref{pr:muinvZ} that
this $T/H$-action on the base is locally free.

We will first consider the quotient of the base of the fibre
bundle~\eqref{eq:s-xh}, and then we will state a proposition which
describes the quotient of the total space.

The submanifold $X^H\subset
X$ is a closed symplectic submanifold, stable under $T$, and the
restriction of the moment map $\mu$ to $X^H$ gives a moment map for
the $T$-action on $X^H$. Hence the quotient of the base can be
described as a symplectic quotient
\begin{equation*}
\(X^{H}\cap\mu\inv(q)\)/T = \syq{X^H}{q}{T}.
\end{equation*}
This looks like a singular kind of quotient: $q$ is not a regular
value of $\mu$, for instance. But the appearance of singularity is an
illusion: we know a priori that $H$ acts trivially on the manifold
$X^H$, and so by Fact~\ref{fact:glo-geo-mom} we know that the image
under $\mu$ of each component of $X^H$ must lie in some affine
hyperplane $\mathcal{S}\subset\td$ (parallel to $\Lie(T/H)^*$). Now
$q$ \emph{is} a regular value in $\mathcal{S}$ for the restriction of
$\mu$, thought of as a map to $\mathcal{S}$ (the fact that $q$ is 
a regular value in this sense is equivalent to the condition that $Z$ cross
each wall in its interior). Hence $\mu\inv(q)\cap X^H$ is a compact
closed submanifold of $X^H$, and its quotient $\syq{X^H}{q}{T}$ is a
compact symplectic orbifold. This kind of symplectic quotient is
explained in more detail in section~\ref{sec:rels}.

\begin{proposition}\label{pr:PHq}
There exists a complex vector orbibundle 
\begin{equation*}
\nu \to \syq{X^H}{q}{T}
\end{equation*}
together with an action of $H$ on $\nu$, covering the trivial action
on $\syq{X}{q}{T}$, and such that the set of fixed points equals the
zero section, such that 
\begin{equation*}
P_{(H,q)} \iso S(\nu)/H \to \syq{X^H}{q}{T}.
\end{equation*}
Here $S(\nu)$ denotes the unit sphere bundle in $\nu$ (relative to a
choice of invariant metric).  In the case that the symplectic quotient
$\syq{X}{q}{T}$ is a free quotient, $\nu$ is a vector bundle, induced
by the normal bundle $\nu X^H$.
\end{proposition}

The vector bundle $\nu\to \syq{X}{q}{T}$ is not uniquely defined: to
defined it we choose a complementary subgroup $T'\subset T$ so that
$T=T' \ti H$. Then $T'$ defines a lift of the action of $T/H$ on $X^H$
to its normal bundle $\nu X^H$, and we let this action define the
induced vector orbibundle (as defined in appendix \ref{app:orb}) over
$\syq{X^H}{q}{T}$.  Then the $H$-action on $\nu X^H$ naturally
descends to an action on $\nu$, and we will show that $S_{(H,q)}/T =
S(\nu)/H$.

\begin{proof}[Proof of Proposition \ref{pr:PHq}]
The space $S_{(H,q)}$ is
formed from $\nu X^H$ by the operations of taking the sphere bundle,
restriction, and taking the quotient by $T$.  By decomposing $T$ as
$T'\ti H$ we can take the quotient by $T$ in two stages. The proof
then amounts to permuting the order of these operations (and seeing
that the result is indepent of the order of operations).
Explicitly, 
$S_{(H,q)}/T = 
\(S(\eval{\nu X^{H}}_{X^{H}\cap\mu\inv(q)})/T'\)/H
= S(\nu)/H.$

The complex structure on $\nu$ is induced by fixing an invariant
almost complex structure on $X$, compatible with the symplectic form
$\om$ (see e.g. \cite[Proposition 2.48]{mcd-sal:int-sym}). It follows by
$T$-invariance that the normal bundle $\nu X^H$ is an invariant complex vector
bundle, so that the complex structure descends to the quotient $\nu$.
\end{proof}

\begin{remarks}\begin{enumerate} 
\item In general $H$ will act on the fibres of $\nu$ with both
positive and negative weights (recall that $H$ is oriented, and so has
a natural identification with $S^1$) and we can thus decompose $\nu$
into the positive and negative weight subbundles
$\nu=\nu^+\oplus\nu^-$. Letting $\bar{\nu}^-$ denote the same
underlying real vector bundle as $\nu^-$, but with the conjugate
complex structure (that is, with multiplication by $i$ replaced by
multiplication by $-i$), then $S(\nu)/H$ can be identified with a
weighted projectivization of $\nu^+\oplus\bar{\nu}^-$.  Although
this describes the diffeomorphism type of $S(\nu)/H$, the natural
orientation of $S(\nu)/H$ (definition \ref{def:ort-PHq}) is not given by
this description.
\item The vector bundle $\nu$ depends on the choice of $T'$. However
the quotient $S(\nu)/H$ is independent of this choice, as we can see
from its description as $S_{(H,q)}/T$.  Changing the choice of $T'$
has the effect of tensoring $\nu$ with a certain line bundle, but this
change doesn't affect the quotient $S(\nu)/H$. This can be seen as a
generalization of the fact that the projectivization of a complex
vector bundle bundle is invariant under tensoring the vector bundle
with a line bundle.
\end{enumerate}\end{remarks}

\section{The orbifold singularities and orientation of the wall-crossing-cobordism}

\subsection{Orbifold singularities in the wall-crossing-cobordism}

We now address the question of the orbifold {\singularities} in the
wall-crossing-cobordism $W/T$.  These arise from points in $W$ whose
stabilizer subgroup is nontrivial. To be more precise, since we allow
for the possibility that there is some finite subgroup of $T$ which
stabilizes every point in $X$, the orbifold {\singularities} arise
from points in $W$ whose stabilizer subgroup is larger than the
generic one.

\begin{lemma}\label{lem:orb_sing}
Let $F\subset T$ be a finite subgroup, and $X^F\subset X$ the subset
of points fixed by $F$. Then $X^F$ is a closed symplectic submanifold of $X$,
transverse to $\bd W$, and also to the interior of $W$. It follows
that the wall-crossing-cobordism $W/T$ is an orbifold-with-boundary.
\end{lemma}

This lemma gives both coarse information, and very fine information. 
The coarse information provided by this lemma is that the
wall-crossing-cobordism is an orbifold-with-boundary, which we will see
is oriented, and hence satisfies Stokes's theorem.

However, this lemma actually makes it possible to determine the
structure of the orbifold singularities quite accurately. This is
because each $X^F$ is a closed symplectic submanifold of $X$, and it
follows that the restriction of the moment map $\mu$ to $X^F$ gives a
moment map for the action of $T$ on $X^F$, where the
walls and chambers of the image of $\mu$ for $X^F$ being a
subset of the corresponding walls and chambers for $X$. Thus we can treat all
the arguments in this paper as applying simultaneously to $X$ and to
$X^F$: each symplectic quotient $\syq{X}{p}{T}$ contains the symplectic
quotient $\syq{X^F}{p}{T}$, as does each wall-crossing-cobordism, and
so on. (We won't have cause to carry out such a detailed analysis in
this paper.)

\begin{proof}[Proof of lemma \ref{lem:orb_sing}]
The fact that $X^F$ is a closed manifold is a standard result of the
theory of compact group actions on manifolds (proved using an
equivariant exponential map, see e.g. \cite{bre:int-com}), and an easy
averaging argument shows that the restriction of the symplectic form
$\om$ to $X^F$ is nondegenerate. 

Now, let $X^*\subset X$ denote the set of points with finite stabilizer
subgroup. Then $W\subset X^*$, by construction. Any $p\in\td$ is a
regular value for $\eval{\mu}_{X^*}$, and using the same argument as in
the proof of Proposition~\ref{pr:muinvZ} we see that $\mu\inv(p)\cap X^*$
is transverse to $X^F$. It follows $X^F$ is transverse to $W$, and to the
boundary components $\mu\inv(p_0)$ and $\mu\inv(p_1)$.

It remains to show transversality to the boundary components $S_{(H,q)}$.
The description of $S_{(H,q)}$ as the sphere bundle in a vector bundle makes
it clear that the stabilizer doesn't depend on the radius of the
sphere.

By varying the radius of the sphere, we can foliate $W$ locally by a
one-parameter family of submanifolds. Since $X^F$ is transverse to
$W$, to show transversality to one of the leaves of this foliation, we
must simply show that for every point in the intersection with $X^F$,
there is a tangent vector to $X^F$ which is transverse to the leaves
of the foliation.  But this follows from the fact that the
stabilizer subgroup is independent of the radius of the sphere.
\end{proof}

\subsection{Orienting the wall-crossing-cobordism}

In this subsection we define an orientation on the
wall-crossing-cobordism $W/T$. We then calculate the induced
orientations on its various boundary components. 

\begin{definition}\label{def:ort-w-c-c}
The orientation is extremely easy from a conceptual point of view:
$W/T$ is foliated by symplectic orbifolds $\syq{X}{p}{T}\cap W/T$, for
$p\in Z$, and the normal bundle to this foliation is identified with
$TZ$ by the moment map. Thus the symplectic orientation of the leaves,
combined with the orientation of $Z$ in which the positive direction
goves from $p_0$ to $p_1$, gives an orientation of the
wall-crossing-cobordism $W/T$.

To carry this out explicitly, we begin by fixing a metric on $W/T$. 
Let $x$ be a point in $W$, denote by $[x]$ the
corresponding point in $W/T$, and set $p=\mu(x)$. 
We assume for simplicity that $[x]$ is a smooth point of $W/T$ (but by
using orbifold metrics and orbifold differential forms, as described in
appendix \ref{app:orb}, this construction also works at the orbifold points).
By construction, $d\mu$ is surjective at $x$ (fact
\ref{fact:inf-geo-mom}). Moreover, since $\mu$ is $T$-invariant, it
descends to a map from $W/T$ to $Z$. 
We can thus decompose the tangent space  $T_{[x]}(W/T)$ 
into the kernel and the cokernel of $d\mu$. 
Identifying these spaces explicitly gives us
\begin{equation*}
T_{[x]}(W/T) \iso T_{[x]}\syq{X}{p}{T} \oplus T_pZ.
\end{equation*}
Now $\syq{X}{p}{T}$ is a symplectic orbifold, and we
denote its symplectic form by $\om_p$. Using the above decomposition,
we can extend $\om_p$ to $T_{[x]}(W/T)$. 
denoting the extension by
$\wtl{\om}_p$ (this $2$-form will not necessarily be closed, but it
will be nondegenerate on the tangent spaces to the leaves). 
Let $Z$ be parametrized by the variable $t$, with $t=0$ at $p_0$ and
$t=1$ at $p_1$.
Then
\begin{equation*}
\wtl{\om}^k_p\wedge \mu^*dt
\end{equation*}
(where $\dim \syq{X}{p}{T} = 2k$) defines a top-degree form, and hence
an orientation of $W/T$ at $[x]$. But the above construction can
be simultaneously applied to every smooth point of $W/T$, with
the resulting form varying smoothly, hence orienting $W/T$. 
\end{definition}

\begin{remark} 
In fact, the above definition can be enhanced in a straighforward
manner to define a `complex orientation' of $W/T$. We won't need it in
this paper, however.
\end{remark}

The rest of this section is taken up with describing in a precise way
an orientation on the wall-crossing boundary components of $W/T$, and
then stating the result that this orientation equals the induced
boundary orientation. We give two definitions, and then state this
result (the proof of which is given in appendix~\ref{app:ort-bd}).

\begin{definition}\label{def:ort-V}
Let $V$ be an oriented real vector space, and suppose the oriented group
$H\iso S^1$ acts on $V$, fixing only the origin. We define the \sdef{induced
orientation} of $S(V)/H$ (where $S(V)$ denotes the unit sphere in $V$
relative to an invariant metric).
Given a point $v\in S(V)$, denote by $H\c.v\in S(V)/H$ the associated
$H$-orbit. There is a natural isomorphism
\begin{equation*}
T_{H\c.v}(S(V)/H) \oplus \R^+\c.v \oplus \la{h} \iso T_vV \iso V,
\end{equation*}
where $\R^+\c.v$ denotes the ray from the origin through $v$.  We
define the orientation of $S(V)/H$ to be that orientation which is
compatible with the above isomorphism together with the given
orientations of $\R^+$, $\la{h}$, and $V$.
\end{definition}

For example, let $V=\C^n$, and let $H\iso S^1$ act with weight
$1$. Then $S(V)/H$ is naturally identified with complex projective
space, and the orientation we have defined agrees with the orientation
induced by the complex structure.  Similarly, if $H$ acts with
positive weights, then $S(V)/H$ is a weighted projective space, and
the above-defined orientation again agrees with the orientation
induced by the complex structure (see appendix \ref{app:wtd-proj-b}
for more details).

We now define an orientation on the boundary components of the
wall-crossing-cobordism corresponding to wall-crossings. We then prove
that this agrees with the induced boundary orientation.

\begin{definition}\label{def:ort-PHq}
Recall that proposition
\ref{pr:PHq} identifies the boundary component of the wall-crossing-cobordism
corresponding to the pair $(H,q)$ as the total space of the bundle
\begin{equation*}
S_{(H,q)}/T = S(\nu)/H \to \syq{X^H}{q}{T}.
\end{equation*}
where $\nu \to \syq{X^H}{q}{T}$ is a vector bundle induced by the
normal bundle $\nu X^H$ and a decompostion of $T$ as $T'\ti H$. 
We orient this space as follows. Since $X^H$ is a symplectic
submanifold, the symplectic orientations of $X$ and of $X^H$ 
induce a natural orientation on the normal
bundle $\nu X^H$, which descends (by invariance) to the induced bundle
$V$. Combining this with definition \ref{def:ort-V} and  the
orientation of $H$ given in definition \ref{def:w-c-d} gives an
orientation of the fibres of the bundle $S(\nu)/H \to
\syq{X^H}{q}{T}$. The base is a symplectic quotient, and we orient
it by its symplectic form. We then orient the total space $S(\nu)/H$ by
the product orientation.
(the order is irrelevant, since both the base
and fibre are even-dimensional).
\end{definition}

\begin{lemma}\label{lem:ort-bd}
\hide{\footnote{The short proof: 1. The
orientations of the two symplectic quotients $\syq{X}{p_i}{T}$ follow
immediately from the definition (\ref{def:ort-w-c-c}). 2. Given a
wall-crossing component, flip the symplectic form on the subbundle of
$\nu X^H$ carrying negative-weight action. (To flip the symplectic
form on the normal bundle to a symplectic submanifold, let $\pi:\nu
Y\to Y$ be the projection, and consider $2\pi^*\om_Y - \om$. This is
symplectic in a neighbourhood of the zero-section. Now consider the
submanifold to be the positive normal bundle over $X^H$.) With this
flipped symplectic form, the boundary component is now the same as if
we had placed the point $p_1$ just before the wall-crossing, so in
this case we know the orientation of the wall-crossing component. In
general, we have flipped the orientation of $X$ by $(-1)^d$, where
$2d$ is the dimension of the negative-weight normal bundle to
$X^H$. Thus this has flipped the orientation of the
wall-crossing-cobordism by the same factor, and thus we must flip the
orientation of the boundary component by the same factor to
compensate. But this is precisely taken care of by our fibre
orientation.
In an infinite-dimensional manifold this argument using the
orientation of $X$ doesn't work, but the implied
calculation, which can be carried out purely on the (finite
dimensional) wall-crossing-cobordism will give the same answer.}
}%%%edih
Let the wall-crossing-cobordism $W/T$ be oriented as in definition
\ref{def:ort-w-c-c}. Then the induced boundary orientation of 
$\syq{X}{p_0}{T}$ is $ -(\om_{p_0}^k)$, and of
$\syq{X}{p_1}{T}$ is $\om_{p_1}^k$
(where $\om_{p_i}$ denote the respective induced symplectic forms), 
and the induced boundary orientation of each $P_{(H,q)}$ is equal to
the product orientation defined in \ref{def:ort-PHq} above.
\end{lemma}

The proof is conceptually rather simple, but keeping track of the
various vector spaces involved in a comprehensible way makes it quite
long, and it has been relegated to appendix~\ref{app:ort-bd}. 

\section{Theorem A: a summary of the existence and properties of the
wall-crossing-cobordism.}\label{sec:thma}

\newcommand{\statetheoremA}{
\begin{theoremA}
Suppose $p_0,p_1\in\td$ are regular values of the
moment map $\mu$, and let
$Z\subset\td$ be path joining $p_0$ and $p_1$  which is transverse to
$\mu$. There there are two objects naturally associated to
$Z$. The first is a finite set
$\data(Z)$, consisting of
pairs $(H,q)$, where $H\iso S^1$ is a subgroup of $T$, and $q$ is a
point in $\td$. And the second object naturally associated to $Z$ is an
oriented cobordism, whose boundary equals
\begin{equation*}
-\syq{X}{p_0}{T}\quad\sqcup\quad \syq{X}{p_1}T%
 \quad\sqcup\quad %
 \bigsqcup_{(H,q)\in\data(Z)} P_{(H,q)}.
\end{equation*}
For each pair $(H,q)\in\data(Z)$ the space $P_{(H,q)}$ is the total
space of a bundle
over the compact symplectic orbifold $\syq{X^H}{q}{T}$, whose fibres
are weighted projective spaces.
\end{theoremA}
}
\statetheoremA

Moreover
\begin{enumerate}
\item The cobordism arises as the quotient, by $T$, of a
submanifold-with-boundary $W\subset X$, such that the $T$-action on
$W$ is locally free.
\item The points $p_0$ and $p_1$ need not lie in the image of
$\mu$. If either lies outside the image of $\mu$, then the associated
boundary component is empty.
\item The boundary component $-\syq{X}{p_0}T$ denotes $\syq{X}{p_0}T$
with the negative of its symplectic orientation.
\item Each space $P_{(H,q)}$ can be described as follows.
There exists a complex vector orbibundle $\nu \to
\syq{X^H}{q}{T}$, with an action of $H$ on the fibres, such that
\begin{equation*}
P_{(H,q)}=S(\nu)/H\to\syq{X^H}{q}{T}. 
\end{equation*}
The bundle $\nu\to\syq{X^H}{q}{T}$ is induced by
the normal bundle $\nu X^H$, and depends on the choice of a complement
to $H$ in $T$; however the bundle $P_{(H,q)}\to\syq{X^H}{q}{T}$ is
independent of this choice.
\item The induced orientations on the spaces $P_{(H,q)}$ are given by
the product of the symplectic  orientation of $\syq{X^H}{q}{T}$ and a
natural orientation on the fibres, defined in terms of the oriented
group $H$, and the oriented fibres of $\nu$.
\item The wall-crossing-data $\data(Z)$ is determined by the
arrangement of walls in $\td$ (which can be deduced from the fixed
point data of $(X,T,\mu)$), together with the path $Z$.
\end{enumerate}

\section{The localization map and the wall-crossing formula}
\label{sec:w-c-f}

In this section we fix our attention on a single wall-crossing. Fixing
notation, we suppose $p_0, p_1\in\td$ are regular values of $\mu$, joined by a
transverse path $Z$ having a single wall-crossing at $q$, and we let
$H\iso S^1$ be the oriented subgroup associated to the wall. 

Theorem~A says, roughly, that the symplectic quotients
$\syq{X}{p_0}{T}$ and $\syq{X}{p_1}{T}$ are in some way related by the
symplectic quotient $\syq{X^H}{q}{T}$. Theorem~B gives a
cohomologically precise version of this.

\begin{theoremB}
There is a map
\begin{equation*}
\loc_H:\HH_T^*(X) \to \HH_{T/H}^*(X^H)
\end{equation*}
such that, for any $a\in\HH_T^*(X)$, 
\begin{equation*}
\ig{\kappa(a)}{\syq{X}{p_0}T} - \ig{\kappa(a)}{\syq{X}{p_1}T} =
  \ig{ \kappa(\loc_H(\eval{a}_{X^H}))}{\syq{X^{H}}{q}{T}}.
\end{equation*}
(The maps $\kappa$ on the left hand side are the natural maps
\mbox{$\HH_T^*(X)\to\HH^*(\syq{X}{p_i}T)$}
and on the right hand side is the natural map
\mbox{$\HH_{T/H}^*(X^H) \to \HH^*(\syq{X^H}{q}{T})$.})

Moreover, for any component $X^H_i\subset X^H$, the restriction of
$\loc_H(a)$ to $X^H_i$ only depends on the restriction of $a$ to
$X^H_i$.
\end{theoremB}

Recall that $\syq{X^H}{q}{T}$ can be considered to be a symplectic
quotient of $X^H$ by the quotient group $T/H$ (expained in
section~\ref{sec:w-c-d}); and the various maps denoted by $\kappa$ are
defined by restriction to the relevant submanifold, followed by the
natural identification of the equivariant cohomology of this manifold
with the rational cohomology of its quotient.

We call $\loc_H$ the \sdef{localization map}: we first define
$\loc_H$, and then we prove theorem~B. In the next section we give an
explicit formula for $\loc_H$ in terms of characteristic classes.  The
localization map is the key to an inductive process, which will allow
us to localize calculations to the fixed points $X^T$. We will carry
out the induction in section~\ref{sec:f-p-c}.

\begin{definition}\label{def:loc}
The localization map $\loc$ depends on the triple $(X,T,H)$, where $X$
is a symplectic manifold, $T$ is a compact torus which acts on $X$
(preserving the symplectic form), and $H\iso S^1$ is an oriented
subgroup of $T$. In this section, $X$ and $T$ will be fixed, and we
will write $\loc_H$ to denote the dependence on the oriented subgroup
$H$ (in later sections will decorate the symbol $\loc$
with any data that is not obvious from the context.)

Given $X$ and $T$, then $\loc_H$ is the (degree-lowering) map
\begin{equation*}
\loc_H:\HH_T^*(X) \to \HH_{T/H}^*(X^H)
\end{equation*}
defined as follows.  Let $S(\nu X^H)$ denote the sphere bundle in the
normal bundle $\nu X^H$ to $X^H$ in $X$. We then denote by $p$ and
$\pi$ the projections
\begin{equation*}
\xymatrix{S(\nu X^H) \ar[rr]^{/H} \ar[rd]_p && 
  S(\nu X^H)/H \ar[ld]^{\pi} \\
  & X^H}
\end{equation*}
Let $\pi_*$ denote integration over the fibres of $\pi$ (where the
fibres are oriented according the definition \ref{def:ort-V}, using
the symplectic orientation of the normal bundle to $X^H$).
Then we let $\loc_H$ equal the
composition
\begin{equation*}
\xymatrix{ & \HH^*_{T}(S(\nu X^H)) \ar[r]^{/H}_{\iso} % 
  & \HH^*_{T/H}(S(\nu X^H)/H) \ar[d]^{\pi_*} \\
\HH^*_{T}(X) \ar[r]^{i^*} & \HH_T^*(X^H) \ar[u]_{p^*}  %
  & \HH^*_{T/H}(X^H)}
\end{equation*}
where $i:X^H\hra X$ denotes the inclusion, and the map $\HH^*_{T}(S(\nu
X^H)) \xra[\iso]{/H} \HH^*_{T/H}(S(\nu X^H)/H)$ is the natural map on
equivariant cohomology induced by the locally free quotient (see for example
\cite{ati-bot:mom-map}).
\end{definition}

\begin{proof}[Proof of theorem~B] 
The proof is a straightforward exercise involving
identifying the various maps involved, and repeatedly using the fact
that integration over the fibre commutes with restriction (together
with some general facts about equivariant cohomology.)

Let $j:W\hra X$ denote the inclusion. Then, for any $a\in\HH^*_T(X)$,
we have $j^*(a)\in\HH^*_T(W)$, and we write
\begin{equation*}
j^*(a)/T \in \HH^*(W/T),
\end{equation*}
for the corresponding naturally induced class (recall that the
$T$-action is locally free on $W$, and we are taking cohomology with
rational coefficients).

Since the wall-crossing-cobordism $W/T$ is an oriented
orbifold-with-boundary, it follows that
the boundary is homologous to zero (fact \ref{fact:orb-bd}), and hence
\begin{equation*}
\ig{j^*(a)/T}{\bd(W/T)} = 0.
\end{equation*}
Using the identification of the boundary of $W/T$ (theorem~A), we thus get
\begin{equation*}
- \ig{j^*(a)/T}{\syq{X}{p_0}{T}} + \ig{j^*(a)/T}{\syq{X}{p_1}{T}} 
+ \ig{j^*(a)/T}{P_{(H,q)}} = 0.
\end{equation*}
We rewrite this, letting $i:S_{(H,q)}\hra X$ denote the inclusion, and
identifying the maps $\kappa$:
\begin{equation*}
- \ig{\kappa(a)}{\syq{X}{p_0}{T}} + \ig{\kappa(a)}{\syq{X}{p_1}{T}} 
+ \ig{i^*(a)/T}{P_{(H,q)}} = 0.
\end{equation*}

Letting $\pi$ denote the projection 
\begin{equation*}
\pi:P_{(H,q)} \to X_{(H,q)} = \syq{X^{H}}{q}{T}
\end{equation*}
and $\pi_*$ denote integration over the fibres of $\pi$, then we have
\begin{equation*}
\ig{i^*(a)/T}{P_{(H,q)}} = \ig{\pi_*(i^*(a)/T)}{X_{(H,q)}}.
\end{equation*}

Thus we have been reduced to proving
\begin{equation}\label{eq:reduced-to}
\pi_*(i^*(a)/T) = \kappa( \loc_H(a)).
\end{equation}

We will now use two naturality properties of integration over the
fibre, for maps in the
commutative diagram
\begin{equation}\label{cd:pi}
\xymatrix@=2ex@C1.5ex{
S_{(H,q)} \ar@{^{(}->}[rr] \ar[dd]^{/H} && S(\nu X^H) 
  \ar[dd]^{/H} \\ %1
\\ %2
S_{(H,q)}/H \ar@{^{(}->}[rr] \ar[ddd]^{/(T/H)} \ar[ddr]^{\wtl{\pi}}
  && S(\nu X^H)/H \ar[rrdd]^{\wtl{\wtl{\pi}}} \\ %3
\\ %4
 & \mu\inv(q)\cap X^H \ar@{^{(}->}[rrr] \ar[ddd]^{/(T/H)} &&& X^H \\ %5
P_{(H,q)} \ar[ddr]^{\pi} \\ %6
\\ %7
& \syq{X^H}{q}{T} \\ %8
}
\end{equation}
Letting $i:S_{(H,q)}\hra X$ and $\wtl{\wtl{i}}:S(\nu
X^H)\hra X$ denote the inclusions, we have
\begin{equation*}
\pi_*(i^*(a)/T) = \wtl{\pi}_*(i^*(a)/H)/(T/H).
\end{equation*}
This is because of the first naturality property of integration over
the fibre: it commutes with simultaneous quotient of the base and the
total space. 

The second naturality property of integration over the fibre is that it
`commutes with restriction'. Concretely, in our case, this gives
\begin{equation*}
\wtl{\pi}_*(i^*(a)/H) = 
  \eval{\wtl{\wtl{\pi}}_*%
  (\wtl{\wtl{i}}{}^*(a)/H)}_{\mu\inv(q)\cap X^H}.
\end{equation*}

Now, in the diagram 
\begin{equation*}
\xymatrix{S(\nu X^H) \ar@{^{(}->}[rr]^{\wtl{\wtl{i}}} 
  \ar[rd]_p &&X\\
  & X^H \ar@{^{(}->}[ru]_k \\
}
\end{equation*}
an easy scaling argument shows that $\wtl{\wtl{i}}$ is
equivariantly homotopic to $k\o p$. Hence 
\begin{equation*}
\begin{array}{rl}
{\wtl{\wtl{\pi}}_*(\wtl{\wtl{i}}{}^*(a)/H)} %
  &= {\wtl{\wtl{\pi}}}_*(p^*(k^*(a))/H) \\
  &=  \loc_H(a) \\
\end{array}
\end{equation*}
since this turns out to be precisely the definition of $ \loc_H$, with
the data $(X,T,H)$. 

Putting all this together, we thus have
\begin{equation*}
\begin{array}{rl}
\pi_*(i^*(a)/T) &= %
  \(\eval{ \loc_H(a)}_{\mu\inv(q)\cap X^H}\)/(T/H) \\
  &= \kappa( \loc_H(a)) \\
\end{array}
\end{equation*}
by definition of $\kappa$. But this proves equation
\eqref{eq:reduced-to}, and hence, by the arguments preceding equation
\eqref{eq:reduced-to}, we have completed the proof.
\end{proof}

\section{The wall-crossing formula in terms of characteristic classes}
\label{sec:w-c-f2}

By giving an explicit formula for the localization map in terms
of characteristic classes, we can restate a more explicit version
of the wall-crossing formula (which we call theorem~B${}'$.)

We can give an explicit formula for the localization map $ \loc_H$ using
the definitions and results of appendix \ref{app:wtd-proj-b}. Using
this explicit formula, we can then recast theorem~B in a more explicit
form. Before carrying this out, we must give a definition, which will
help us account for the possibility that $X^H$ has a number of
components, and describe the way a decomposition of $T$ induces a
decomposition of a cohomology class.

\begin{definition}
Let $Y$ be a connected manifold with an action of $T$. We define $o_T(Y)$
to be the order of the maximal subgroup of $T$ which stabilizes every
point in $Y$ (so $o_T(Y)=1$ if and only if $T$ acts effectively on $Y$). (In
every case which we consider, this number will be finite).
We extend this definition to the case in which $Y$ may have a number of
components by defining $o_T(Y)$ to be the degree-$0$ cohomology class
which restricts to give this number on each component.
\end{definition}

Now suppose $T'\subset T$ is a complement to
$H$, so that $T=T' \ti H$. Then the restriction of any class
$a\in\HH^*_T(X)$ to $X^H$ decomposes
\begin{equation}\label{eq:decomp-a}
\eval{a}_{X^H} = \sum_{i\ge 0} a_i\ot u^i
\end{equation}
according to the natural isomorphism
\begin{equation*}
\HH_T^*(X^H)\iso\HH_{T'}^*(X^H)\ot\HH^*(BH),
\end{equation*}
where $u\in\HH^2(BH)$
is the positive generator (with respect to the orientation of $H$ defined
in \ref{def:w-c-d}).

We then have

\begin{proposition}\label{pr:expl}
Let $a\in\HH^*_T(X)$, and suppose $T'\subset T$ is a complement to
$H$, so that $T=T' \ti H$. Then
\begin{equation*}
 \loc_H(a) = \frac{o_T(X)}{o_{T/H}(X^H)}
\sum_{i\ge 0} a_i\smile  s^w_{i-r+1},
\end{equation*}
where the classes $a_i\in\HH^*_{T'}(X^H)$ are defined by the natural
decomposition of $a$ given in equation
\eqref{eq:decomp-a} above, and $s^w_i$ denotes the $i$-th $T'$-equivariant
weighted Segre class of $(\nu X^H,H)$ (definitions \ref{def:wtd-s-c}
and \ref{def:ewtd}), and $r$ is the function, constant on connected
components of $X^H$, such that $2r=\rank(\nu X^H)$. 
\end{proposition}

\begin{proof}
We will show how this proposition follows from the integration formula
proved in appendix \ref{app:wtd-proj-b}, namely proposition
\ref{pr:segre-int} (together with its `equivariant enhancement',
equation \eqref{eq-int}). 

Explicitly, we are using the vector bundle $\nu X^H\to X^H$ and the
groups $H$ and $T'$ in place of the vector bundle $V\to Y$ and the
groups $S^1$ and $G$ in appendix~\ref{app:wtd-proj-b}.

We need first to give the normal bundle $\nu X^H$ a complex structure
compatible with its symplectic form, so that the definition of the
orientation of $S(\nu X^H)/H$ used in  appendix~\ref{app:wtd-proj-b}
agrees with its natural orientation (definition~\ref{def:ort-PHq}). 
And second, we must show that our factor
${o_T(X)}/{o_{T/H}(X^H)}$ is equal, on each component of $X^H$, to
the factor $k$ in the appendix.

Firstly, general principles in symplectic topology imply that there
exists a $T$-invariant almost complex structure $J:TX\to TX$,
compatible with the symplectic form $\om$, and such that $TX^H$ is
stable under $J$ (see, for example, McDuff and Salamon
\cite[Proposition 2.48]{mcd-sal:int-sym}). Such an almost complex
structure gives the normal bundle $\nu X^H$ a complex structure, and
the orientations induced by the complex structure and the symplectic
form agree (equation \eqref{eq:complex-ort}), and thus we can apply
Proposition \ref{pr:segre-int} with this complex structure. 

Secondly, we need to show that, for each component $X^H_i$ of
$X^H$, we have 
\begin{equation*}
k={o_T(X)}/{o_{T/H}(X^H_i)},
\end{equation*}
where $k$ is the greatest common
divisor of the weights of the $H$-action on the fibres of $\nu X^H_i\to
X^H_i$. But using the decomposition $T=T'\ti H$, together with lemma
\ref{lem:orb_sing}, it is clear that $k=o_H(\nu X^H_i)$, and that
$o_T(X) = o_{T'\ti H}(\nu X^H_i) = o_{T'}(X^H_i)\c.o_H(\nu X^H_i)$.
\end{proof}

We can now rewrite theorem~B using this explicit identification.

\newcommand{\statetheoremBp}{
\begin{theoremBp}
Suppose $p_0, p_1\in\td$ are regular values of $\mu$, joined by a
transverse path $Z$ which has a single wall-crossing, at $q$. Let
$H\iso S^1$ be the subgroup associated to the wall, and choose
$T'\subset T$ so that $T=T'\ti H$.

Then there are characteristic classes $s^w_i\in\HH^{2i}_{T'}(X^H)$
(the equivariant weighted Segre classes of $\nu X^H$,
as defined in \ref{def:wtd-s-c}) such that, for any $a\in\HH_T^*(X)$,
\begin{equation*}
\ig{\kappa(a)}{\syq{X}{p_0}T} - \ig{\kappa(a)}{\syq{X}{p_1}T} =
  \ig{ \kappa\( \frac{o_T(X)}{o_{T/H}(X^H)} 
  \textstyle\sum_{i\ge 0} a_i\smile  s^w_{i-r+1}\)}{\syq{X^{H}}{q}{T'}}.
\end{equation*}
where $r$ is the function, constant on connected components of $X^H$,
such that $2r=\rank(\nu X^H)$; and the classes $a_i\in\HH^*_{T'}(X^H)$
are defined by restricting $a$ to $X^H$ and decomposing, as in
equation \eqref{eq:decomp-a} above.  (The map $\kappa$ on the left
hand side of the main equation is the natural map
\mbox{$\HH_T^*(X)\to\HH^*(\syq{X}{p_i}T)$} and on the right hand side
is the natural map\\ \mbox{$\HH_{T'}^*(X^H) \to
\HH^*(\syq{X^H}{q}{T'})$.})
\end{theoremBp}
}

\statetheoremBp

\section{A generalization of a transverse path and its data}
\label{sec:rels}

This is the first of three sections in which we apply the preceding
results inductively, ending up with results concerning the $T$-fixed
points of $X$. In this section we generalize the notion of a
transverse path, and the associated data. In section \ref{sec:f-p-c}
we show how this generalized data corresponds to cobordisms involving
the fixed points $X^T$, and in section \ref{sec:f-p-f} we show how
this generalized data governs integration formulae localized at the
fixed points.

\subsection{A $\sbt$-transverse path and its data}

We begin with a straightforward generalization of the notion of a
transverse path, and its associated data.  Recall that $X$ is a
symplectic manifold, with an action of the torus $T$, and with
associated moment map $\mu:X\to\td$.  Let $\sbt\subset T$ be a
subtorus. In section~\ref{sec:con-w-c-c} we saw how $\Lie(T/\sbt)^*$
can be considered to be a subspace of $\td$ via a natural embedding
(it is a subspace of dimension $\dim T-\dim\sbt$). Recall also that
$X^\sbt$, the set of points fixed by $\sbt$, is a closed symplectic
submanifold of $X$. 

\begin{fact}
If $X^\sbt_i$ is any connected component of $X^\sbt$, then we have:
\begin{enumerate}
\item The restriction of $\mu$ to $X^\sbt_i$ gives a moment map for the
$T$-action on $X^\sbt_i$;
\item The image $\mu(X^\sbt_i)$ lies in an affine translate $\sbs\subset\td$ of
$\Lie(T/\sbt)$;
\item The $T$-action on  $X^\sbt_i$ descends to a $T/\sbt$-action; and
\item Composing the restriction of $\mu$ with an
identification of $\sbs$ with $\Lie(T/\sbt)^*$ gives a moment map for
the $T/\sbt$-action on $X^\sbt_i$.
\end{enumerate}
\end{fact}

Hence we define, in analogy with section \ref{sec:con-w-c-c}

\begin{definition}\label{def:t-reg}
Given $q\in\td$, set $\sbs:=q+\Lie(T/\sbt)^*$.  We say $q$ is
\sdef{$\sbt$-regular} if $\mu$ maps some component of $X^\sbt$ to
$\sbs$, and for each such component, the point $q$ is regular value
for the restriction of $\mu$, thought of as a map
to $\sbs$.
\end{definition}

For example, using the notions of `wall' and `interior' from definition
\ref{def:int}, if $H\iso S^1$ is a subgroup of $T$, and if $q$ lies in
a wall corresponding to $H$, then $q$ is $H$-regular iff $q$ lies in
the interior of this wall.

\begin{definition}\label{def:t-tran}
Let $\sbs$ be an affine translate of $\Lie(T/\sbt)^*$, and suppose
$q_0, q_1\in\sbs$ are $\sbt$-regular values. Then a path $Z\subset\td$
from $q_0$ to $q_1$ is \sdef{$\sbt$-transverse} if it is contained in the
subspace $\sbs$, and for each component of $X^\sbt$ which $\mu$ maps
to $\sbs$, the path $Z$ is transverse to the restriction of $\mu$,
thought of as a map to
$\sbs$.
\end{definition}

\begin{definition}\label{def:t-data}
Suppose $Z\subset\sbs$ is a $\sbt$-transverse path, with endpoints the
$\sbt$-regular values $q_0$ and $q_1$. We define the
\sdef{wall-crossing data} for $Z$ to be the set
\begin{equation*}
\data(Z) := \{ (H,q) \mid \text{$H$ is a subtorus of $T$ with
$\sbt\subset H$, and $q\in Z\cap\mu(X^H)$} \}
\end{equation*}
Applying proposition \ref{pr:muinvZ}, it follows that $H/\sbt\iso
S^1$, and we orient $H/\sbt$ as in definition \ref{def:w-c-d}, that
is, we orient $Z$ so that the positive direction goes from $q_0$ to
$q_1$, and we orient $H/\sbt$ compatibly.
\end{definition}

\subsection{The module of relations}

We now define a module which records the data from all possible
$\sbt$-transverse paths simultaneously.

\begin{definition}
An \sdef{oriented $\sbt$-flag of subtori} in $T$ is a collection of
subtori
\begin{equation*}
\flagt = (1 = H_0 \subset H_1 \subset H_2 \subset \h.. \subset H_k =
\sbt \subset T),
\end{equation*}
such that $H_i$ is an $i$-torus, and each $H_i/H_{i-1}\iso S^1$ is
given an orientation.
\end{definition}

\begin{definition}
We define the $Z$-module $\CA$ by
\begin{equation*}
\CA:=\bigoplus_{\sbt\subset T}\CA_\sbt,
\end{equation*}
as $\sbt$ runs through all subtori of $T$, where
\begin{equation*}
\CA_\sbt:=\bigoplus \Z(\flagt,q)
\end{equation*}
is the set of formal linear combinations of pairs $(\flagt,q)$, where
$q$ is $\sbt$-regular and $\flagt$ is an oriented $\sbt$-flag of
subtori.
\end{definition}

Note that $\CA_\sbt$ will be nontrivial for only finitely many $\sbt$,
namely those for which there exists a $\sbt$-regular value. These
correspond to the $\sbt$ such that there is some point $x\in X$ whose
stabilizer subgroup has identity component $\sbt$ (the fact that there
are only finitely many such $\sbt$ is a standard fact in the theory of
group actions on manifolds \cite{bre:int-com,kaw:the-tra}). We also
note that $\CA_T$ corresponds to the $T$-fixed points of $X$: if
$(\flagt,q)\in\CA_T$ then $q\in\td$ is one of the finite set of points
in the set $\mu(X^T)\subset\td$. 

\begin{definition}
We now define the submodule of relations $\CR\subset\CA$.  There are
two kinds of generators of $\CR$. The first kind comes from a pair
consisting of a $\sbt$-transverse path $Z$ and an oriented $\sbt$-flag
of subtori $\flagt$, for any choice of subtorus $\sbt$. The associated
generator of $\CR$ is the sum
\begin{equation*}
- (\flagt,q_0) + (\flagt,q_1) + \sum_{(H,r)\in\data(Z)} (\flagt\union H,r),
\end{equation*}
where $q_0$ and $q_1$ are the endpoints of $Z$, and $\flagt\union H$
denotes the oriented $H$-flag defined by concatenating $\flagt$ and
$H$, with $H/\sbt$ oriented as in $\data(Z)$. The second kind of
generator of $\CR$ corresponds to points which are outside the image
of $\mu$: for any subtorus $\sbt\subset T$, suppose $q$ is a
$\sbt$-regular value and let $\flagt$ be an oriented $\sbt$-flag. If
$q\notin\mu(X^\sbt)$ then
\begin{equation*}
(\flagt,q)
\end{equation*}
is a generator of $\CR$.
Finally, given $(\flagt,q)\in\CA$, we write $[\flagt,q]$ for its
equivalence class in the quotient module $\CA/\CR$.
\end{definition}

Since $X$ is compact, for any regular value $p_0\in\td$, there is a
path $Z$ starting at $p_0$ and ending outside the image of the moment
map. The corresponding fact is true for each $X^\sbt\subset X$. Hence

\begin{lemma}\label{lem:descent}
For any $(\flagt,q)\in\CA$ we have
\begin{equation*}
[\flagt,q] = \sum_{i\in I} [\flagt_i,v_i]
\end{equation*}
in $\CA/\CR$, 
where $I$ is a finite indexing set, and each $(\flagt_i,v_i)\in\CA_T$.
\end{lemma}

\section{Cobordisms between symplectic quotients and bundles over the 
fixed points}\label{sec:f-p-c}

\begin{figure}[!hbtp]\label{fig:f-p-c}
\psfrag{j}{\f$\mu\inv(p_0)$}
\psfrag{k}{\f$\mu\inv(p_1)$}
\psfrag{a}{\f$\syq{X}{p_0}{T}$}
\psfrag{b}{\f$\syq{X}{p_1}{T}$}
\psfrag{p}{\f$p_0$}
\psfrag{s}{\f$p_1$}
\psfrag{q}{\f$q_0$}
\psfrag{r}{\f$q_1$}
\psfrag{m}{\f$\mu$}
\psfrag{x}{\f$X$}
\psfrag{z}{\f$Z$}
\psfrag{w}{\f$W$}
\psfrag{d}{\f$/T$}
\psfrag{c}{\f$W/T$}
\psfrag{0}{\f$Z_0$}
\psfrag{1}{\f$Z_1$}
\psfrag{2}{\f$q_2$}
\psfrag{3}{\f$q_3$}
\psfrag{4}{\f$r_0$}
\psfrag{5}{\f$r_1$}
\psfrag{6}{\f$S_0$}
\psfrag{7}{\f$S_1$}
\psfrag{8}{\f$\Lie(T/\sbt_0)^*$}
\psfrag{9}{\f$\Lie(T/\sbt_1)^*$}
\psfrag{t}{$\td$}
\begin{center}
\includegraphics[width=5.5in]{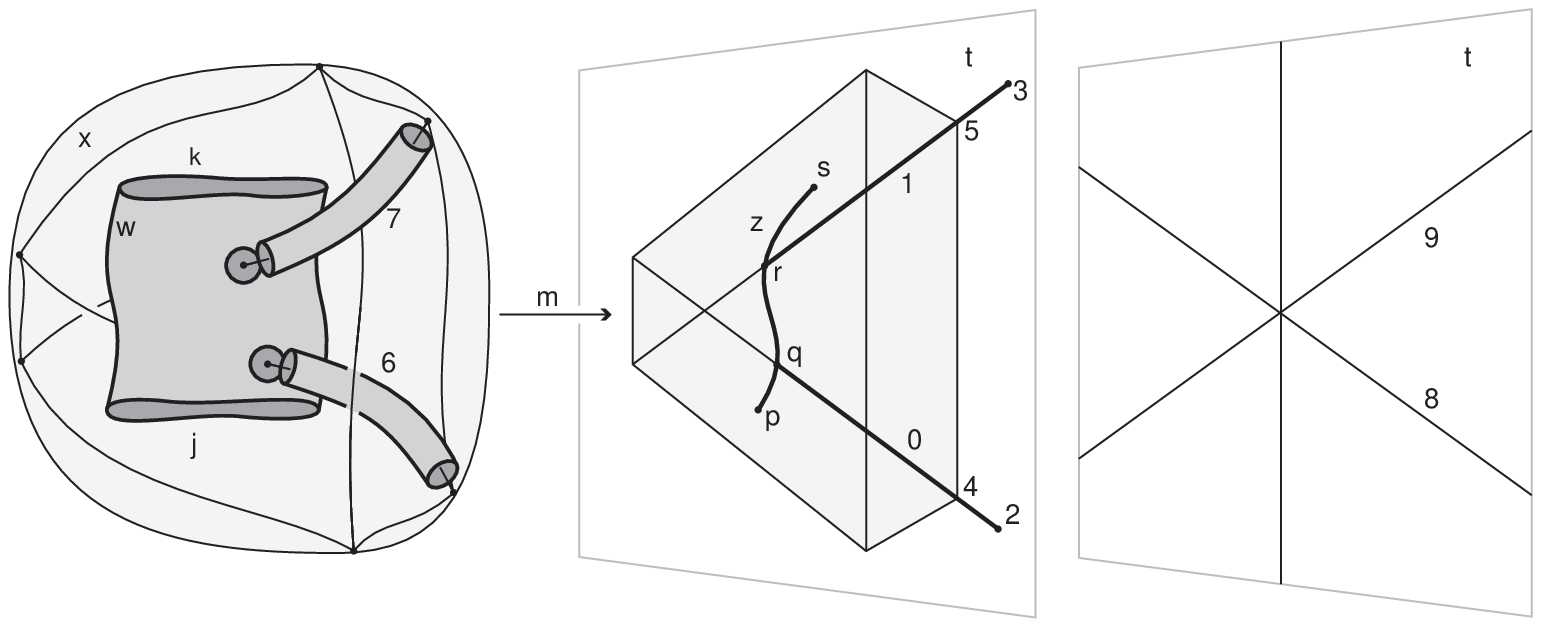}
\end{center}
\caption{The definitions of this section: $Z_0$ is
a $\sbt_0$-transverse path, with endpoints the $\sbt_0$-regular values
$q_0,q_2$. Since $\sbt_0$ is a $1$-torus, there is only one
$\sbt_0$-flag, namely $\flagt_0:=(1\subset\sbt_0)$. The wall-crossing
data of $Z_0$ is the pair $(T,r_0)$. Now associated to $Z_0$ is a
submanifold-with-boundary $W_0\subset X^{\sbt_0}$, and the space
labelled $S_0$ is $S_{\flagt_0}(\eval{\nu X^{\sbt_0}}_{W_0})$ (as
described in the proof of theorem~C). An analogous description holds
for $Z_1$.}
\end{figure}

In this section we show how the relations defined in the previous
section correspond to cobordisms. We begin by defining, for each
generator  $(\flagt,q)$ of $\CA$, a space $P_{(\flagt,q)}$. We will
then show how `relations', i.e. finite sums in the submodule $\CR$,
correspond to cobordisms between these spaces. The constructions in
this section are illustrated in figure~\ref{fig:f-p-c}.

\subsection{The spaces involved}

For every pair $(\flagt,q)$, where $\flagt$ is a $\sbt$-flag and
$q\in\td$ is a $\sbt$-regular value, we will define an associated
space $P_{(\flagt,q)}$. We first describe $P_{(\flagt,q)}$ in two
special cases, and then give the general definition.  In the case that
$\sbt=\{1\}$ is the trivial group, then the only $\sbt$-flag is the
trivial flag, which we denote by $1\subset T$, and a $\sbt$-regular value is
just a regular value of the moment map $\mu:X\to\td$. In this case
\begin{equation*}
P_{(1\subset T,q)} = \syq{X}{q}{T}.
\end{equation*}
If $Z\subset\td$ is a transverse path, and $(H,q)\in\data(Z)$ is one
of its wall-crossing pairs, then it follows that $q$ is an $H$-regular
value, and $H\iso S^1$ defines the oriented $H$-flag $1\subset
H\subset T$, and we have
\begin{equation*}
P_{(1\subset H \subset T,q)} = P_{(H,q)},
\end{equation*}
where the space on the right is the wall-crossing space defined in
equation~\eqref{eq:PHq}.

\begin{definition}
Suppose the torus $\sbt$ acts on the complex vector space $V$, with
$0\in V$ the only point fixed by $\sbt$. Then associated to every flag
of subtori of $\sbt$ is a submanifold of $V$ on which the
$\sbt$-action is locally free (this submanifold may be empty).  To
define the submanifold, we first define a canonical decomposition of
$V$.  Let $\flagt=(1=H_0\subset H_1\subset \h.. \subset H_k = \sbt)$
be a $\sbt$-flag, that is, a full flag of subtori of $\sbt$. There is
an associated flag of subspaces of $V$, stable under the
$\sbt$-action:
\begin{equation*}
V=V^{H_0} \supset V^{H_1} \supset \h.. \supset V^{H_k} = \{0\}
\end{equation*}
where $V^{H_i}$ is the subspace fixed by $H_i$. We define $V_i\subset
V$ to be the orthogonal complement to $V^{H_i}$ in $V^{H_{i-1}}$,
relative to a $\sbt$-invariant metric, for $1\le i\le k$. Then
$V_i\iso V^{H_{i-1}}/V^{H_i}$, and these subspaces define a
decomposition of $V$ into subrepresentations
\begin{equation*}
V = V_1 \oplus V_2 \oplus \h.. \oplus V_k.
\end{equation*}
We set 
\begin{equation*}
S_\flagt(V):= S(V_1) \ti S(V_2) \ti \h.. \ti S(V_k) \quad \subset V
\end{equation*}
where $S(V_i)$ is the unit sphere, relative to an invariant
metric. Note that $S_\flagt(V)$ will be nonempty precisely when each $V_i$
is nontrivial, that is, when each inclusion is strict in the flag of
subspaces $V^{H_0} \supset V^{H_1} \supset \h.. \supset V^{H_k}$. 

Finally, we define
\begin{equation*}
P_\flagt(V):=S_\flagt(V)/\sbt.
\end{equation*}
This is a locally free quotient, and hence has the structure of an
orbifold.
An orientation of $V$ induces an orientation on
$P_\flagt(V)$ as follows. We fix an orientation of each $V_i$ so that
the product orientation equals the given orientation of $V$. We then
orient each $S(V_i)/T_i$ by applying the formula of
definition~\ref{def:ort-V}, and give $P_\flagt(V)$ the induced product
orientation (see the end of this section, where the structure of
$P_\flagt(V)$ is described in more detail).
\end{definition}

\begin{remarks}
\begin{enumerate}
\item To see that the $\sbt$-action is locally free on $S_\flagt(V)$
we choose a  decomposition of $\sbt$ which is compatible with
$\flagt$, that is
\begin{equation*}
\sbt = T_1 \ti T_2 \ti \h.. \ti T_k,
\end{equation*}
where each $T_i\iso H_i/H_{i-1}\iso S^1$. Then the above
decomposition of $V$ has the property that the $T_i$-action on $V_i$
leaves only $0\in V_i$ fixed, so that the $T_i$-action on $S(V_i)$ is
locally free.
\item The quotient $P_\flagt(V)$ can be described as a $k$-fold
`tower' of weighted projective bundles, where $k=\dim\sbt$. We make
some remarks about this at the end of this section.
\end{enumerate}
\end{remarks}

\begin{definition}
We now observe that we can apply the above construction both fibrewise and
equivariantly. Suppose $T\supset \sbt$ acts on a manifold $Y$, and the
action lifts to a complex vector bundle $V\to Y$. Moreover, suppose
that the stabilizer subgroup of each point $y\in Y$ is $\sbt$. Then
each fibre $V_y$ is a $\sbt$-representation and, if $0\in V_y$ is the
only point fixed by $\sbt$, we define the submanifold
$S_\flagt(V_y)\subset V_y$ by applying the above
construction. Applying this to each fibre simultaneously, relative to
a $T$-invariant metric, gives a submanifold 
\begin{equation*}
S_\flagt(V)\subset V
\end{equation*}
which is stable under the action of $T$.

We now apply this fibrewise construction
to the symplectic manifold $X$, with $T$-moment map
$\mu$. Given a pair $(\flagt,q)$, where $\flagt$ is a $\sbt$-flag and
$q\in\td$ is a $\sbt$-regular value, we let $S_{(\flagt,q)}$ be the
result of applying the above construction with
$Y:=X^\sbt\cap\mu\inv(q)$ and $V:=\eval{\nu X^\sbt}_{Y}$, with a
$T$-invariant almost complex structure, compatible with the symplectic
form, giving $V$ the structure of a complex vector bundle.
That is
\begin{equation*}
S_{(\flagt,q)}:= 
S_{\flagt}\(\eval{\nu X^\sbt}_{X^\sbt\cap\mu\inv(q)}\to 
   X^\sbt\cap\mu\inv(q)\).
\end{equation*}
Using an equivariant exponential map to identify a neighbourhood of
the zero-section of $\nu X^\sbt$ with a neighbourhood of $X^\sbt$ in
$X$ we can consider $S_{(\flagt,q)}$ to be a submanifold of $X$. It
follows from the above construction and the fact that $q$ is
$\sbt$-regular that the $T$-action is locally free on
$S_{(\flagt,q)}$.
We then define 
\begin{equation*} 
P_{(\flagt,q)} := S_{(\flagt,q)}/T
\end{equation*}
which we see is the total space of a bundle over the symplectic
quotient $\syq{X^\sbt}{q}{T}$ with fibre $P_\flagt(\nu_x X^\sbt)$. We
note that in the case that the symplectic quotient
$\syq{X^\sbt}{q}{T}$ is smooth, this is an honest fibre bundle, but in
general, the symplectic quotient $\syq{X^\sbt}{q}{T}$ may have
orbifold singularities, in which case the above construction defines
$P_{(\flagt,q)}\to\syq{X^\sbt}{q}{T}$ as an orbibundle.
\end{definition}

\subsection{The cobordism theorem}

\begin{theoremC}
Suppose
\begin{equation*}
\sum_i c_i[\flagt_i,q_i] \  = \ 0\ \in\ \CA/\CR, \qquad c_i\in\Z.
\end{equation*}
Then there exists an oriented manifold $W$, with a
locally free action of $T$, and a $T$-equivariant map
\begin{equation*}
W\to X
\end{equation*}
such that
\begin{equation*}
\bd(W/T)\  \iso\  \bigsqcup_i c_i P_{(\flagt_i,q_i)}.
\end{equation*}

In particular, for any regular value $p\in\td$ of the moment map, the
symplectic quotient $\syq{X}{p}{T}$ is cobordant in the above sense to a
union of spaces $P_{(\flagt_i,v_i)}$, for $(\flagt_i,v_i)\in\CA_T$,
and such spaces can be described as towers of weighted projective
bundles over components of the fixed points $X^T$.
\end{theoremC}

\begin{proof}
Since we can glue together oriented cobordisms along their boundaries, it is
enough to show the above result in the case that $\sum_i
c_i(\flagt_i,q_i)$ is one of the relations which generate $\CR$. 

Each such relation comes from a  $\sbt$-transverse path $Z$, and a
choice of $\sbt$-flag $\flagt$, and so we fix such a $Z$ and
$\flagt$. Then we wish to find a manifold $W$ with a locally free
$T$-action, together with an equivariant map $W\to X$, such that 
\begin{equation*}
\bd(W/T) \iso 
- P_{(\flagt,q_0)} + P_{(\flagt,q_1)} + \sum_{(H,r)\in\data(Z)} 
 P_{(\flagt\union H,r)}.
\end{equation*}

In fact we can construct a submanifold $W\subset X$ with this
property. The first step is to apply theorem~A to $Z$. Explicitly, $Z$
lies in a subspace $\sbs\subset\td$, which we can identify with
$\Lie(T/\sbt)^*$. We then apply theorem~A, where the symplectic
manifold consists of those components of $X^\sbt$ which $\mu$ maps to
$\sbs$, the torus is $T/\sbt$, and the moment map is given by $\mu$
with the identification of $\sbs$ with $\Lie(T/\sbt)^*$. This gives a
submanifold-with-boundary $W'\subset X^\sbt$, with a locally free
action of $T/\sbt$, and with boundary
\begin{equation*}
-X^\sbt\cap\mu\inv(q_0) \sqcup X^\sbt\cap\mu\inv(q_1) \sqcup
\bigsqcup_{(H,r)\in\data(Z)} \eval{S(\nu X^H:X^\sbt)}_{X^H\cap\mu\inv(r)}
\end{equation*}
where $\nu X^H:X^\sbt$ denotes the normal bundle to $X^H$ in
$X^\sbt$, and $q_0,q_1$ are the endpoints of $Z$.

But, since $W'\subset X^\sbt$ is a submanifold-with-boundary, with a
locally free action of $T/\sbt$, it follows that
\begin{equation*}
W := S_\flagt \( \eval{\nu X^\sbt}_{W'} \to W' \)
\end{equation*}
defines a submanifold of $X$ with a locally free action of $T$, and
$\bd W =S_\flagt \( \eval{\nu X^\sbt}_{\bd W'} \to \bd W' \)$. 

Finally, using the fact that
\begin{equation*}
S_\flagt \( \eval{\nu X^\sbt}_{S(\nu X^H:X^\sbt)} \) = 
S_{\flagt\cup H} (\nu X^H),
\end{equation*}
we see that $W/T$ has the desired boundary, thus proving the result.
\end{proof}

\subsection{The structure of the spaces $P_{(\flagt,q)}$}

Let $(\flagt,q)\in\CA_\sbt$, that is, $q$ is a $\sbt$-regular value
and $\flagt$ is a $\sbt$-flag.

\begin{proposition}
The space $P_{(\flagt,q)}$ is the total space of a tower
\begin{equation*}
P_{(\flagt,q)} = P_1 \xra{\pi_1} P_2 \xra{\pi_2} \h.. \xra{\pi_{k-1}} 
P_k \xra{\pi_k} \syq{X^\sbt}{q}{T}
\end{equation*}
where $k=\dim\sbt$, and each $\pi_i$ is an orbibundle projection with
fibre a weighted projective space.
\end{proposition}

We can identify the spaces $P_i$ explicitly (see below). The explicit
formulae for cohomology pairings in the next section follow from these
identifications (although they can also be deduced by inductively
applying theorem~B).

\begin{proof} For simplicity of notation we treat explicitly the case
in which $\sbt=T$, so that $P_{(\flagt,q)}$ is a bundle over certain
components of the fixed point set, and we assume such components
consist of a single point. Adapting these arguments to deal with the
general case is straightforward.

Letting $x$ be the point in question, we set $V=T_xX$, so that $V$ is
a complex representation of $T$.

We choose 
a decomposition of $\sbt=T$ which is compatible
with $\flagt$, that is
\begin{equation*}
\sbt = T_1 \ti T_2 \ti \h.. \ti T_k,
\end{equation*}
where each $T_i\iso H_i/H_{i-1}\iso S^1$.

Then, tracing through the definitions, we see that
\begin{enumerate}
\item For $1\le i,j\le k$, each $T_i$ acts on each $V_j$;
\item If $j > i$ then $T_i$ acts trivially on $V_j$;
\item The $T_i$ action on $V_i$ leaves only $0$ fixed.
\end{enumerate}

We now note the following general fact.

Fact: Suppose $Y_1\ti Y_2$ is acted on by $T_1\ti T_2$, such that the $T_1$-action
is free on $Y_1$ and trivial on $Y_2$, and the $T_2$-action is free on
$Y_2$. Then the projection $Y_1\ti Y_2$ descends to a projection
\begin{equation*}
(Y_1\ti Y_2)/(T_1\ti T_2) \to Y_2/T_2
\end{equation*}
with fibre $Y_1/T_1$.

Hence, defining
\begin{equation*}
\begin{split}
S_i &:= S(V_i) \ti S(V_{i+1}) \ti \h.. \ti S(V_k),\quad\text{and}\\
P_i &:= S_i/(T_i\ti T_{i+1} \ti \h.. \ti T_k),
\end{split}
\end{equation*}
we see that the natural projection $S_i\to S_{i+1}$
descends to a projection $\pi_i:P_i\to P_{i+1}$, with  fibre
$S(V_i)/T_i$. As in Proposition~\ref{pr:PHq}, we can thus express
$\pi_i:P_i\to P_{i+1}$ as the weighted projectivization of the  complex vector
bundle induced by $V_i\ti S_{i+1} \to S_{i+1}$.
\end{proof}

\begin{definition}
We can use the above description to orient $P_{(\flagt,q)}$. Recall
that $\flagt$ is an \textit{oriented} flag: this is equivalent to the
statement that each $T_i\iso S^1$ is oriented. Since each $V_i$ is a
complex subrepresentation of $V$, each $V_i$ has an orientation. We
thus use the formula of definition~\ref{def:ort-V} to orient each
$S(V_i)/T_i$, and we give $P_{(\flagt,q)}$ the induced product
orientation.
\end{definition}

\section{Localizating integration formulae to the fixed
points}\label{sec:f-p-f}

In this section we show how the relations defined in section
\ref{sec:rels} correspond to integration formulae.
We begin by defining a map which generalizes the localization map
$\loc_H$ defined in section~\ref{sec:w-c-f}. We then state theorem~D
in terms of this map. We then give an explicit formula for this
localization map in terms of characteristic classes.

\begin{definition}
Let $\sbt$ be a subtorus of $T$, and let $\flagt$ be an oriented $\sbt$-flag.
Then we define the map
\begin{equation*}
\loc_\flagt : \HH_T^*(X) \to \HH_{T/\sbt}^*(X^\sbt)
\end{equation*}
as follows. Firstly, in the case that $\sbt=\{1\}$ is the trivial
subtorus, so that $\flagt=(1\subset T)$ is the trivial flag, then we
define $\loc_\flagt$ to be the identity map. Otherwise we set
\begin{equation*}
\loc_\flagt:= \loc_{H_k/H_{k-1}} \o \h.. \loc_{H_2/H_1} \o \loc_{H_1}.
\end{equation*}
Here $H_i$ is the subtorus in the flag $\flagt$:
\begin{equation*}
\flagt = (1 = H_0 \subset H_1 \subset H_2 \subset \h.. \subset H_k =
\sbt \subset T),
\end{equation*}
and 
\begin{equation*}
\loc_{H_i/H_{i-1}}:\HH^*_{T/H_{i-1}}(X^{H_{i-1}}) \to 
\HH^*_{T/H_{i}}(X^{H_{i}})
\end{equation*}
 is the localization map of 
definition~\ref{def:loc}, with data consisting of the triple\\
$(X^{H_{i-1}},T/H_{i-1},H_i/H_{i-1})$. Recall that $H_i/H_{i-1}\iso
S^1$ is assumed to be oriented.
\end{definition}

After stating theorem~D we will give an explicit formula for
$\loc_\flagt$ using a decomposition of $T$ and characteristic classes.

Note that $\loc_\flagt$ can equivalently be defined via integration
over the fibre of the bundle
\begin{equation*}
P_\flagt(\nu X^\sbt) \to X^\sbt
\end{equation*}
in an analogous way to the definition of $\loc_H$
(definition~\ref{def:loc}).

\begin{theoremD}
Suppose
\begin{equation*}
\sum_i c_i[\flagt_i,q_i] \  =\  0\ \in\ \CA/\CR, \quad c_i\in\Z.
\end{equation*}
Then for any $a\in\HH^*_T(X)$,
\begin{equation*}
\sum_i c_i 
  \ig{ \kappa(\loc_{\flagt_i}(a))}{\syq{X^{\sbt_i}}{q_i}{T}}\ =\ 0.
\end{equation*}
where, for each $i$, the flag $\flagt_i$ is a $\sbt_i$-flag, and where
$\kappa$ is the relevant natural map from the equivariant cohomology
of a manifold to the ordinary cohomology of its symplectic quotient,
as described in the notation section of the Introduction.

Moreover, for each flag $\flagt_i$, the class $\loc_{\flagt_i}(a)$ only
depends on the restriction of $a$ to the submanifold $X^{\sbt_i}$.
\end{theoremD}

The proof consists of straightforward unwinding of the definitions, and can
be seen to either follow from theorem~C, or from theorem~B, using
inductive arguments analogous to those in the proof of theorem~C.  We
give a concrete application of this theorem in section~\ref{sec:cp2n},
in which we calculate some cohomology pairings on the symplectic
reduction of products of $\cp^2$.

\subsection{A formula for $\loc_\flagt$ in terms of characteristic
classes}

Suppose $\flagt$ is an (oriented) $T$-flag of subtori (that is, we
suppose $\sbt=T$). We consider the map 
\begin{equation*}
\loc_\flagt:\HH_T^*(X)\to \HH^*(X^T).
\end{equation*}
We first observe that, for any component $F\subset X^T$ of the fixed
point set and any class $a\in \HH_T^*(X)$, the restriction of
$\loc_\flagt(a)$ to $F$ only depends on the restriction of $a$ to $F$
(this follows from the definition of $\loc_H$).

Since $T$ acts trivially on $F$, we have
$\HH_T^*(F)\iso\HH^*(F)\ot\HH_T^*(\pt)$. We
choose a  decomposition
\begin{equation*}
T = T_1 \ti T_2 \ti \h.. \ti T_d
\end{equation*}
compatible with the flag $\flagt$, that is, where each 
$T_i\iso H_i/H_{i-1}\iso S^1$. This gives a set of generators
$\{u_1,u_2,\h..,u_d\}$ of $\HH_T^*(\pt)$ so that
\begin{equation*}
\HH_T^*(F)\iso\HH^*(F)\ot\Q[u_1,u_2,\h..,u_d].
\end{equation*}
Explicitly, $u_i$ is the equivariant first Chern class of the
representation of $T$ on $\C$ where $T_i$ acts with weight $1$ (recall
$T_i$ is oriented), and the other $T_j$ act trivially.

We now define the map
\begin{equation*}
\begin{split}
\ell_i:\Q[u_i] &\to \HH^*(F)\ot\Q[u_{i+1},\h..,u_d], \quad\text{by}\\
u_i^{j+k_i} &\mapsto s^{T_{i+1}\ti\h..\ti T_d}_j(V_i,T_i)\\
\end{split}
\end{equation*}
where $k_i+1=\rank V_i$ and $s^{T_{i+1}\ti\h..\ti T_d}_j(V_i,T_i)$ is
the equivariant weighted Segre class (equivariant with respect to
$T_{i+1}\ti\h..\ti T_d$) of the bundle $V_i\to F$. 

Then $\ell_i$ extends to a map
\begin{equation*}
\tilde\ell_i:\HH^*(F)\ot\Q[u_{i},\h..,u_d] \to 
\HH^*(F)\ot\Q[u_{i+1},\h..,u_d]
\end{equation*}
by tensoring with the identity map on the complement of
$\Q[u_i]$. Thus $\tilde\ell_i$ is a homomorphism of
$\HH^*(F)\ot\Q[u_{i+1},\h..,u_d]$-modules.

Now for any $a\in\HH^*_T(X)$, the restriction $\eval{a}_{F}$ can be
decomposed 

We then have
\begin{proposition}
\begin{equation*}
\loc_\flagt(a)=o_T(X)\c.\tilde\ell_d \o \tilde\ell_{d-1} \o \h.. \o
\tilde\ell_1(\eval{a}_{F}).
\end{equation*}
where $o_T(X)$ is the order of the maximal subgroup of $T$ which fixes
every point in $X$.
\end{proposition}

We will use this formula in the explicit calculations of
section~\ref{sec:cp2n}.

\begin{proof}
This follows by repeatedly applying proposition~\ref{pr:expl}, using
explicit identifications coming from the choice of decomposition of
$T$. For example, we have
\begin{equation*}
H_i = T_1 \ti T_2 \ti \h.. \ti T_i
\end{equation*}
and so on.
\end{proof}

\section{A more refined look at the module of relations}

In section~\ref{sec:rels}, we gave a number of definitions,
culminating in the definitions of the modules $\CA$ and $\CR$. The aim
of those definitions was to keep track of the relations arising from
paths as simply as possible.
In this section we give `improved' versions of these
definitions. The result of these improved definitions will be that
$\CA$ and $\CA/\CR$ will be much smaller, and should have properties
which more accurately reflect the manifold $X$. The cost of this
improvement is that the definitions are somewhat more subtle.

This section contains no new results: its only aim is to give
alternative definitions which may be useful in some applications.
Theorems C and D are still true with the improved definitions given in
this section.

\begin{definition}
We say an action of a Lie group $G$ on a manifold $Y$ is \sdef{locally
effective} if there is some point in $Y$ whose stabilizer subgroup is
finite.
\end{definition}

\begin{definition}
Let $\sbt\subset T$ be a subtorus. We denote by $X^{[\sbt]}\subset
X^\sbt$ the connected components of $X^\sbt$ on which the
$T/\sbt$-action is locally effective.
\end{definition}

Note that there are only finitely many subtori $\sbt$ for which $X^{[\sbt]}$
is nonempty.

Given this definition, we redefine the notions of a
\sdef{$\sbt$-regular} point $q$, a \sdef{$\sbt$-transverse} path $Z$,
and the \sdef{wall-crossing data} of a $\sbt$-transverse path by
substituting $X^{[\sbt]}$ for $X^\sbt$ in definitions \ref{def:t-reg},
\ref{def:t-tran} and \ref{def:t-data}.

\begin{definition}
Let $V$ be a $\sbt$-representation, such that $0\in V$ is the only
fixed point. A $\sbt$-flag $\flagt=(1=H_0\subset H_1\subset
H_2\subset\h..\subset H_k=\sbt)$ is called \sdef{$V$-admissible} if
the associated flag in $V$
\begin{equation*}
V=V^{H_0} \supset V^{H_1} \supset \h.. \supset V^{H_k} = \{0\}
\end{equation*}
has each inclusion a strict inclusion.
\end{definition}

\begin{definition}
Let $\sbt$ be a subtorus of $T$, and let $\flagt$ be a $\sbt$-flag and
$q$ a $\sbt$-regular value (using the version of $\sbt$-regular
defined in this section). We say the pair $(\flagt,q)$ is
\sdef{admissible} if there is some point $x\in
X^{[\sbt]}\cap\mu\inv(q)$ such that $\flagt$ is $\nu_x
X^{[\sbt]}$-admissible.
\end{definition}

It is easy to see that the admissible pairs are precisely those pairs
$(\flagt,q)$ for which the space $P_{(\flagt,q)}$ is nonempty.

\begin{definition}
We now redefine $\CA$ to have generators the set of admissible pairs
$(\flagt,q)$. We will redefine the submodule of relations $\CR$ to
come from the data for $\sbt$-transverse paths, as $\sbt$ runs through
all subtori, in the same way as before. However there is one
difference: some of the pairs which arise from the data of a path may
not be admissible, and we simply discard these pairs and construct
relations from the pairs that remain. (The point is that these pairs
correspond to empty spaces, so there is no harm in discarding them).
Explicitly, if $\sbt\subset T$ is a subtorus, $Z$ is a
$\sbt$-transverse path and $\flagt$ is an oriented $\sbt$-flag, then
we take the sum
\begin{equation*}
- (\flagt,q_0) + (\flagt,q_1) + \sum_{(H,r)\in\data(Z)} (\flagt\union H,r),
\end{equation*}
and throw out any terms in this sum which are not admissible pairs. We
then define the resulting sum to be a generator of $\CR$.  (Here, as before,
$q_0$ and $q_1$ are the endpoints of $Z$, and $\flagt\union H$ denotes
the oriented $H$-flag defined by concatenating $\flagt$ and $H$, with
$H/\sbt$ oriented as in $\data(Z)$.)
\end{definition}

The statement that every element of $\CA$ 
can be localized to the fixed points becomes

\begin{proposition}\label{pr:descent}
$\CA/\CR$ is generated by $\CA_T/(\CR\cap\CA_T)$.
\end{proposition}

I conjecture

\begin{conjecture}
For any $0\le i \le \dim T$, let $\CA_i=\bigoplus_{\dim\sbt=i}
\CA_\sbt$. Then $\CA/\CR$ is generated by $\CA_i/(\CR\cap\CA_i)$.
\end{conjecture}

\begin{question}
Using  the `improved' definitions of this section, is the following
`converse' to theorem~D true: Given $a\in\HH^*_T(X^T)$, suppose that,
for every relation $r\in\CR\cap\CA_T$, the sum of integrals of classes
induced by $a$ vanishes (as in theorem~D). Then does $a$ extend to a class $\tilde
a\in\HH^*_T(X)$?
\end{question}

\section{Calculations I: cohomology pairings on 
symplectic quotients of $(S^2)^n$}\label{sec:s2n}

\newcommand{\nth}{\underline{n}}
\newcommand{\sth}{\underline{s}}

Consider the unit sphere $S^2\subset \R^3$. The Euclidean volume on
$\R^3$ restricts to give a symplectic form on $S^2$ (with respect to
which its volume is $4\pi$).  $SO(3)$ acts naturally on $S^2$, and
this action is Hamiltonian (it is possible to identify $\R^3$ with
$\Lie(SO(3))^*$ such that the inclusion of $S^2$ is a moment map).  We
choose a maximal torus $S^1\subset SO(3)$ to be the subgroup which fixes
the north and south poles, and normalize so that the positive
direction in $S^1$ rotates the sphere counterclockwise, as seen from
the north pole.

Let $X=(S^2)^n$, the $n$-fold product, with the diagonal action of
$SO(3)$, and hence $S^1$.  The symplectic form on $X$ is given by the
direct sum of the symplectic forms on the factors. We will fix $n$ to
be odd, and calculate cohomology pairings on $X\d/S^1(0)$. We will
also invoke a formula proved in~\cite{skm:g-t}
to use these pairings to
determine cohomology pairings on $(S^2)^n\d/SO(3)$.

\subsection{The moment map}

The action of $S^1$ on $S^2$ is Hamiltonian, with moment map given by
the height function
\begin{equation*}
\begin{split} 
\mu:S^2 &\to \R \\
x &\mapsto \opht(x). \\
\end{split} 
\end{equation*}
We have $\mu(S^2)=[-1,1]$, $\mu(\text{north pole})=1$,
$\mu(\text{south pole})=-1$. Choosing a compatible almost complex
structure (for example the standard one), the weight of the action on
the tangent space at
the north pole is $1$, and the weight at the south pole is $-1$.

A point in $(S^2)^n$ is given by an $n$-tuple $(x_1,\h..,x_n)$ where
each $x_i\in S^2$.  The action of $S^1$ on $(S^2)^n$ has moment map
given by summing the heights on each of the factors:
\begin{equation*}
\mu(x_1,\h..,x_n) = \opht(x_1) +\h.. + \opht(x_n).
\end{equation*}
Thus $\mu((S^2)^n)=[-n,n]$. 

A point in $X$ is fixed if and only if
each $x_i$ is either the north pole or the south pole. 

\begin{definition}
Let $I$ be any subset of the set $\{1,\h..,n\}$. Then we define the
point $f_I\in(S^2)^n$ by setting $x_i$ to be the south pole if $i\in
I$ and the north pole otherwise.
\end{definition}

This defines a one-to-one correspondence
between the fixed points and the subsets of $\{1,\h..,n\}$. In
particular, the fixed points are isolated, and we have
\begin{equation*}
\mu(f_I) = n-2|I|.
\end{equation*}

Hence $0$ is a regular value of $\mu$ when $n$ is odd. 

\subsection{The integration formula relating the symplectic quotients
by a nonabelian group and by its maximal torus}

Suppose $X$ is a symplectic manifold with a Hamiltonian action of the
nonabelian Lie group $G$, having moment map $\mu_G:X\to\Lie(G)^*$. The
inclusion $T\hra G$ induces a projection $\Lie(G)^*\onto\td$, and
composing of $\mu_G$ with this projection gives a moment map
$\mu_T:X\to\td$ for the action of $T$ on $X$. In the companion
paper~\cite{skm:g-t} the following formula is proven,  relating
integrals on $\syq{X}{0}{G}=\frac{\mu_G\inv(0)}{G}$ to integrals on
$\syq{X}{0}{T}=\frac{\mu_T\inv(0)}{T}$.

\begin{proposition}\label{pr:g-t}
For any $a\in\HH_G^*(X)$, 
\begin{equation*}
\ig{\kappa(a)}{\syq{X}{0}{G}} = \frac{1}{|W|}
\ig{\kappa(a)\smile\prod_{\al\in\roots} c_1(L_\al)}{\syq{X}{0}{T}}.
\end{equation*}
Here $|W|$ denotes the order of the Weyl group, and 
$\roots\subset\td$ denotes the set of roots of $G$ (both
positive and negative). Given a root $\al$, then the complex 
line bundle $L_\al\to\syq{X}{0}{T}$ is the line bundle associated to
the fibering $\mu_T\inv(0)\to\syq{X}{0}{T}$ and the $1$-dimensional
$T$-representation of weight $\al$. 
\end{proposition}

\subsection{The volume of the symplectic quotient}

\begin{definition}
The \sdef{symplectic volume} of a compact symplectic manifold (or orbifold)
$(M^{2n},\om)$ is the integral $\int_M{\om^n/n!}$. From now on we will
refer to the symplectic volume as simply the \sdef{volume}.
\end{definition}

We will now go through the calculations necessary to prove

\begin{proposition}\label{pr:volS2}
For $n$ odd,
\begin{equation*}
\vol((S^2)^n\d/S^1(0)) = \frac{(2\pi)^{n-1}}{(n-1)!}
  \sum_{k=0}^{\frac{n-1}{2}} 
  (-1)^k \choose{n}{k}(n-2k)^{n-1},
\end{equation*}
and
\begin{equation*}
\vol((S^2)^n\d/SO(3)) = -\frac{(2\pi)^{n-3}}{(n-3)!}
  \half\sum_{k=0}^{\frac{n-1}{2}} 
  (-1)^k \choose{n}{k}(n-2k)^{n-3}.
\end{equation*}
\end{proposition}

$X$ is endowed with a line bundle $\CL\to X$ (known in the literature
as the prequantum line
bundle), with a connection whose curvature is $-i\om$. Hence
$c_1(\CL)=\left[\frac{\om}{2\pi}\right]$. The action on $X$ lifts to
an action on $\CL$. Hence the volume of $X\d/S^1(0)$ is given by
\begin{equation*}
\vol(X\d/S^1(0))=\frac{(2\pi)^{n-1}}{(n-1)!}
 \<\kappa(c^{S^1}_1(\CL)^{n-1}),[X\d/S^1(0)]\>
\end{equation*}
where $c^{S^1}_1$ denotes the equivariant first Chern class (and, as
usual, $\kappa:\HH^*_T(X)\to\HH^*(X\d/S^1)$ is the map described in the
introduction.)

In order to evaluate classes on $X\d/SO(3)$ we use the integration
formula, proposition~\ref{pr:g-t} of the companion paper~\cite{skm:g-t}.
We first need a definition.
Let $\al$ be a weight of $S^1$. Then we denote by
$\C_{(\al)}$ the representation induced by $\al$, and set
$\underline{\C}_{(\al)}:=X\ti\C_{(\al)}$, thought of as an
equivariant line bundle over $X$. $\underline{\C}_{(\al)}$ induces a
line bundle on the quotient $X\d/S^1(0)$, which we denote by  
$L_{(\al)}$.
Applying the integration formula, proposition~\ref{pr:g-t},
we have
\begin{multline*}
\vol(X\d/SO3)(0))= \\
\frac{(2\pi)^{n-3}}{(n-3)!}\half \ig{\kappa(c^{S^1}_1(\CL)^{n-3}\smile 
  c^{S^1}_1(\ul{\C}_{(1)})
  \smile c^{S^1}_1(\ul{\C}_{(-1)}))}{X\d/S^1(0)}.
\end{multline*}

\subsection{The calculation}

We now go through the steps necessary to evaluate cohomology classes
on $X\d/T(0)$. Steps 1-3 are independent of the particular class we
wish to evaluate, and steps 4 and 5 depend on the class.

\paragraph{Step 1: Fix $Z$, and enumerate the components $X_{q,i}$.}  
We fix our submanifold $Z$ to be the interval $[0,n+1]$, with 
$p_0=0$ and $p_1=n+1$, which is outside the image of $\mu$. Then
$Z\cap\{\text{walls}\}=\{n-2k \mid k=0\h..\left[\frac{n}{2}\right]\}$. For
$q=n-2k$, 
\begin{equation*}
\{ X_{q,i} \} = \{ f_I \mid |I|=k \}.
\end{equation*}
(There are $\choose{n}{k}$ such points.)

\paragraph{Step 2: Identify $\nu X_{q,i}$.}
Our submanifolds are the points $f_I$. The normal bundle to $f_I$ is
the direct sum of copies of $T_{{\nth}} S^2$ and 
$T_{{\sth}} S^2$.
Hence, with $k=|I|$,
\begin{equation*}
\nu f_I \iso \C^{n-k}_{(1)} \oplus \C^k_{(-1)}
\end{equation*}
Here $\C^m_{(w)}$ denotes $\C^m$ with the $S^1$-action with
weight $w$. Note that the weights are determined by the isomorphism $S^1\xra{\iso} S^1$ induced by
the orientation of $Z$; in our case this is the identity
isomorphism. (If we had instead chosen $Z=[-n-1,0]$, we would have the
orientation-reversing isomorphism.)

\paragraph{Step 3: Calculate $s^w(\nu X_{q,i})$.}
The weighted Segre class lies in $\HH^*_{T/H}(X_{q,i})$. In our case
this is just $\HH^*(f_I)$. We have
\begin{equation*}
s^w_j(\nu f_I) =\begin{cases}
  (\prod\{\text{wts}\})\inv=(-1)^{k} & j=0, \\
  0 & j>0. \\
  \end{cases}
\end{equation*}
Since all the weights are $\pm 1$, we have
$\hcf\{|\text{wts}|\}=1$.

\paragraph{Step 4: Calculate $a$ in terms of local bases.}

\begin{equation*}
\CL|_{f_I} \iso \C_{(\mu(f_I))} = \C_{(n-2k)}.
\end{equation*}
Hence
\begin{equation*}
c^{S^1}_1(\CL)|_{f_I} = (n-2k)u,
\end{equation*}
where $u$ is the positive generator of $\HH^*_{S^1}(\pt)$. And
\begin{equation*}
c^{S^1}_1(\ul{\C}_{(w)})|_{f_I} = wu,\qquad w\in\Z.
\end{equation*}

The two classes we wish to evaluate are
\begin{equation*}
a_1 := c^{S^1}_1(\CL)^{n-1}
\end{equation*}
which gives us the degree of $X\d/S^1(0)$, and 
\begin{equation*}
a_2 := \half c^{S^1}_1(\CL)^{n-3}\smile 
  c^{S^1}_1(\ul{\C}_{(1)})
  \smile c^{S^1}_1(\ul{\C}_{(-1)})
\end{equation*}
which gives the degree of $X\d/SO(3)$. 
We thus have
\begin{equation*}
a_1|_{f_I} = (n-2k)^{n-1} u^{n-1}
\end{equation*}
and
\begin{equation*}
a_2|_{f_I} = - \half (n-2k)^{n-3} u^{n-1}.
\end{equation*}

\paragraph{Step 5: Apply the formula.}

Using the Segre classes calculated above, we have
\begin{equation*}
\loc_{S^1}(a_1)|_{f_I} = (-1)^k (n-2k)^{n-1}
\end{equation*}
and
\begin{equation*}
\loc_{S^1}(a_2)|_{f_I} = -\half (-1)^k (n-2k)^{n-3}
\end{equation*}
Hence, summing over the $\choose{n}{k}$ points $f_I$ with $|I|=k$, and
letting $k$ run from $0$ to $\left[\frac{n}{2}\right]$, we get
proposition \ref{pr:volS2}.

\subsection{Cohomology classes on $(S^2)^n$.}

We wish to describe cohomology classes on $X$, and in particular
understand their restrictions to the fixed points. There are some
standard results which will help us greatly. We first recall these
general results.

Let $G$ be a compact Lie group, with $T\subset G$ the maximal torus,
and $W$ the Weyl group. For any $G$-space $Y$, there is a natural 
action of $W$ on $\HH^*_T(Y)$, and the natural map
$\HH^*_G(X)\to\HH^*_T(X)$ defines
an isomorphism 
\begin{equation*}
\HH^*_G(X) \xra{\iso} \HH^*_T(X)^W
\end{equation*}
(see e.g. \cite[Equation 2.11]{ati-bot:mom-map}).

Suppose $X_1$ and $X_2$ are symplectic manifolds, with Hamiltonian
actions of the group  $G$. Then the homotopy quotients $(X_1)_G$ and
$(X_2)_G$ are cohomologially trivial as bundles over $BG$
\cite[Proposition 5.8]{kir:coh-quo}. This means
that the Serre spectral sequence of the fibering $(X_i)_G\to BG$
degenerates at the $E_2$ term. We give the product $X_1\ti X_2$ the
diagonal $G$-action. Then it follows that
\begin{equation}\label{eq:tens_Ham}
\HH^*_G(X_1\ti X_2) \iso \HH^*_G(X_1) \ot_{\HH^*(BG)} \HH^*_G(X_2).
\end{equation}
In order to see this, we first note that $X_1\ti X_2$ is both a
Hamiltonian $G$-manifold (with moment map given by the sums of the
respective moment maps), and a Hamiltonian $G\ti G$ manifold. 
Now consider the diagonal map
\begin{equation*}
j:BG\hra BG\ti BG.
\end{equation*}
This induces the ring homomorphism
\begin{equation*}
\begin{split}
\HH^*((X_1)_G) \ot_{\HH^*(BG)} \HH^*((X_2)_G) 
  &\to \HH^*((X_1)_G \ti_{BG} (X_2)_G) \\
 a, b &\mapsto j^*(a\ot b).
\end{split}
\end{equation*}
But by degeneracy of the relevant spectral sequences, this must be an
isomorphism of groups, and hence an isomorphism of rings. 
Thus we have equation \eqref{eq:tens_Ham} above.

This means we can represent an equivariant cohomology class on $X_1\ti
X_2$ as a sum of tensor products of equivariant classes on $X_1$ and
$X_2$. Now for our calculations we will only need to know the
restriction of a class to the fixed points. However restriction
commutes with the above tensor
product (and the fixed point
set of $X_1\ti X_2$ is the product of the fixed points of $X_1$ with
the fixed points of $X_2$).

We now specialize to the case at hand. Let
$\nth$ denote the
north pole and ${\sth}$ the south pole of $S^2$.
The restriction map
\begin{equation*}
\HH^*_{S^1}(S^2) \to \HH^*_{S^1}(\{{\nth},{\sth}\}) 
\end{equation*}
is injective. We set 
\begin{equation*}
\HH^*_{S^1}(\{{\nth},{\sth}\}) = \Q[u_{\nth}]\oplus\Q[u_{\sth}]
\end{equation*}
so that for example $u_{\nth}$ is the degree-$2$ generator of the equivariant
cohomology of $\{{\nth}\}$. 

Then the image of the restriction consists of those pairs of polynomials
whose degree-zero terms agree.  (One can see this e.g. by thinking
about the topology of the homotopy quotients.)

The $SO(3)$-equivariant cohomology is the subring invariant under the
Weyl group, in this case $\Z/2\Z$. The nontrivial element of $W$
permutes the north and south poles and simultaneously acts via the
involution on $S^1$. This identifies $u_{\nth}$ with $-u_{\sth}$. 
Hence, fixing a normalization,
\begin{equation*}
\HH^*_{SO(3)}(S^2) \iso \Q[v]
\end{equation*}
with the inclusion given by
\begin{equation*}
\begin{split}
\Q[v] &\hra \Q[u_{\nth}]\oplus\Q[u_{\sth}] \\
a(v) &\mapsto a(u_{\nth})\oplus a(-u_{\sth}). \\
\end{split}
\end{equation*}

Applying equation~\ref{eq:tens_Ham}, we represent a class on $X$ as a sum
of tensor products of classes on $S^2$. 
Hence, consider the class
\begin{equation*}
a^{(1)} \ot a^{(2)} \ot \h.. \ot a^{(n)}
\end{equation*}
where $a^{(i)}$ is an equivariant cohomology class on the $i$-th copy
of $S^2$. We represent $a^{(i)}$ as the pair of polynomials
$(a^{(i)}_{{\nth}},a^{(i)}_{{\sth}})$. A fixed point $f_I$ is
an element of the set $\{{\nth},{\sth}\}^n$, and 
the restriction of $a^{(1)} \ot a^{(2)}
\ot \h.. \ot a^{(n)}$ to $f_I$ is simply the product of the
appropriate polynomials. 
(For example if $I=\emptyset$, then
$f_I=({\nth},\h..,{\nth})$, and the restriction
to $f_I$ is $\prod_{i=1}^n a^{(i)}_{{\nth}}$.)

We shall concentrate on evaluating classes on $X\d/SO(3)$. Any such
class is a linear combination of classes of the form
\begin{equation*}
a=v_1^{l_1} \ot v_2^{l_2} \ot \h.. \ot v_n^{l_n}.
\end{equation*}
$a$ is of top degree when $l_1 + \h.. + l_n = n-3$.
We now calculate $a|_{f_I}$. Define, for $i\in\{1\h..n\}$,
\begin{equation*}
\sg(i)=\begin{cases} 1, &i\notin I;\\
      -1, &i\in I.\\
     \end{cases}
\end{equation*}
Then
\begin{equation*}
\begin{split}
a|_{f_I} &= (\sg(1)u)^{l_1} \c. \h.. \c. (\sg(n)u)^{l_n} \\
         &= u^{\sum l_i} \prod_{i=1}^n \sg(i)^{l_i}. \\
\end{split}
\end{equation*}
Applying the integration formula, and assuming $\sum l_i = n-3$,
\begin{equation*}
\ig{a}{X\d/SO(3)} = -\half\sum_{k=0}^{\frac{n-1}{2}} 
  (-1)^k \sum_{|I|=k} \(\prod_{i=1}^n \sg(i)^{l_i}\).
\end{equation*}
It is clear that this expression only depends on the parity of the
$l_i$, and is invariant under permuting the $S^2$ factors. This will
allow us to deduce quite a lot, but for the moment we will press on
and derive an explicit formula. Define
$J\subset\{1\h..n\}$ by 
\begin{equation*}
i\in J\iff l_i \text{ odd}
\end{equation*}
and set $m=|J|$. Then
\begin{equation*}
\begin{split}
\ig{a}{X\d/SO(3)} &= -\half\sum_{k=0}^{\frac{n-1}{2}} 
  (-1)^k \sum_{|I|=k} (-1)^{|J\cap I|} \\
  &= -\half\sum_{|I|\le\frac{n-1}{2}} 
  (-1)^{|I|} (-1)^{|J\cap I|}. \\
\end{split}
\end{equation*}
Now since $\sum l_i=n-3$, at least one $l_i=0$. By
invariance we may as well assume $J=\{1\h..m\}$, and hence $l_n=0$. 
We can split the above sum
into those $I$ which contain $n$ and those which don't. The resulting
cancellations leave us with
\begin{equation*}
\begin{split}
\ig{a}{X\d/SO(3)} &= -\half
 \sum_{\substack{K\subset\{1\h..n-1\} \\ |K| = \frac{n-1}{2}}}
  (-1)^{|K|} (-1)^{|K\cap\{1\h..m\}|} \\
 &= -\half (-1)^{\frac{n-1}{2}}
 \sum_{\substack{K\subset\{1\h..n-1\} \\ |K| = \frac{n-1}{2}}}
  (-1)^{|K\cap\{1\h..m\}|}. \\
\end{split}
\end{equation*}
From this description one can easily derive explicit computational
formul{\ae}. Alternatively, using the easily-described product structure
in $\HH^*_{SO(3)}(X)$ and Poincar\'e duality in the symplectic
quotient, we can see some classes whose image must vanish on the
symplectic quotient.

\begin{proposition}
Using the identification described above
\begin{equation*}
\HH^*_{SO(3)}((S^2)^n)\iso\Q[v_1,v_2,\h..,v_n]
\end{equation*}
and the natural ring homomorphism
$\kappa:\HH^*_{SO(3)}((S^2)^n)\to\HH^*((S^2)^n//SO(3))$, we have
\begin{multline*}
\ig{\kappa(v_1^{l_1}v_2^{l_2}\h.. v_n^{l_n})}{(S^2)^n\d/SO(3)} \\
\begin{split}
  &= -\half (-1)^{\frac{n-1}{2}}
 \sum_{\substack{K\subset\{1\h..n-1\} \\ |K| = \frac{n-1}{2}}}
  (-1)^{|K\cap\{1\h..m\}|} \\
  &= \half (-1)^{\frac{n-1}{2}} \( \choose{n-1}{\frac{n-1}{2}} - 
2\sum_{j=0}^{\frac{m}{2}}
   \choose{m}{2j} \choose{n-1-m}{\frac{n-1}{2} - 2j}\) \\
\end{split}
\end{multline*}
where $\sum_i l_i = n-3$ and $m$ is equal to the number of odd $l_i$. 
\end{proposition}

It follows, for example, that the ideal $\ker(\kappa)$ contains the elements
$v_i^2-v_j^2$.

\section{Calculations II: volume of the symplectic 
quotient of $(\cp^2)^n$}\label{sec:cp2n}

\subsection{Generalities on $\cp^{k-1}$}

Consider the defining representation of $\Uk$ on $\C^k$.

The maximal
torus $T^k\subset \Uk$ consists of the diagonal matrices 
\begin{equation*}
\{\begin{pmatrix}
e^{it_1} \\
 &e^{it_2} \\
 && \ddots \\
 &&& e^{it_k}
\end{pmatrix}\mid t_i\in\R\}.
\end{equation*}
The moment map for the action of the maximal torus on $\C^k$ is
\begin{equation*}
\mu(z_1,\h.. z_k) = -\half(|z_1|^2,\h..,|z_k|^2).
\end{equation*}

The centre
\begin{equation*}
Z(\Uk) = 
\{\begin{pmatrix}
e^{it} \\
& \ddots \\
&& e^{it}
\end{pmatrix}\mid t\in\R\}
\end{equation*}
acts, with moment map 
\begin{equation*}
\mu_{Z}(z) = -\half\sum |z_i|^2.
\end{equation*}

We can take the symplectic quotient of $\C^k$ by $Z(\Uk)$ at any
negative value, to get $\cp^{k-1}$ (with a scaled symplectic form).
Henceforth, we let $\cp^{k-1}$ denote the symplectic
manifold $\C^k\d/Z(\Uk)(-k)$. This is endowed with prequantum line
bundle $\CL\to\cp^{k-1}$ of degree $k$. $P\Uk$ acts on $\cp^{k-1}$,
and the action lifts to $\CL$.

Let $T$ denote the $(k-1)$-torus. We identify $T$ with the
maximal torus of $P\Uk$ via the inclusion into $\Uk$
\begin{equation*}
T := 
\{\begin{pmatrix}
e^{it_1} \\
 & \ddots \\
 &&e^{it_{k-1}} \\
\end{pmatrix}\}
\hra
\{\begin{pmatrix}
e^{it_1} \\
 & \ddots \\
 &&e^{it_{k-1}} \\
 &&& 1 \\
\end{pmatrix}\}
\end{equation*}
The image of this inclusion is a \emph{slice}:
every element of $T^k$ decomposes in a unique way as a product of
elements of $Z(\Uk)$ and $T$, thus identifying $T$ with
the maximal torus of $P\Uk$. 

Let $\td \iso\R^{k-1}$ have standard basis
$\{e_1,\h..,e_{k-1}\}$. Then $e_j$ corresponds to the representation
\begin{equation*}
e_j : \begin{pmatrix}
e^{it_1} \\
 & \ddots \\
 &&e^{it_{k-1}} \\
\end{pmatrix}
\mapsto e^{it_j}
\end{equation*}

The $T$-action on $\cp^{k-1}$ has fixed points $\{F_j \mid
j=1\h..k\}$, where $F_j$ denotes the point $[0:\h..:1:\h..:0]$
(only the $j$-th coordinate nonzero). We henceforth let $\mu$ denote
the moment map for the action of $T$ on $\cp^{k-1}$. The image of
$\mu$ is the convex hull of the points $\mu(F_j)$. And 
\begin{equation*}
\mu(F_j) = \begin{cases} \(\sum_{i=1}^{k-1} e_i\) - ke_j, &j<k \\
                            \sum_{i=1}^{k-1} e_i\, &j=k \\
           \end{cases}
\end{equation*}

The walls for $\mu$ have corresponding subgroups $H_j\iso S^1$, for
$j=1\h..k$. $H_j$ stabilizes the wall equal to the convex hull of the
points $\{\mu(F_i) \mid i\ne j\}$.
We have
\begin{equation*}
H_j=\{\begin{pmatrix}
1 \\
 & \ddots \\
 &&e^{it} \\
 &&& \ddots \\
 &&&& 1 \\
\end{pmatrix}, t\in\R\}
\qquad\text{for } j=1\h..k-1
\end{equation*}
and
\begin{equation*}
H_k=\{\begin{pmatrix}
e^{-it} \\
 & \ddots \\
 &&e^{-it} \\
\end{pmatrix}, t\in\R\}
\end{equation*}
We let $H_j$ denote the above subgroup, with the isomorphism
$H_j\xra{\iso}S^1$ implied by the above coordinates on $H_j$. We write
$\bar{H}_j$ for the subgroup with the opposite isomorphism with $S^1$.

A set of positive roots for $P\Uk$ is given by $\{e_i-e_j \mid i<j\le
k-1\} \union \{e_i \mid i=1\h..k-1\}$. 

The action of $T$ on the normal bundle to the fixed point $F_j$ is
given by
\begin{equation*}
\nu F_j \iso 
  \begin{cases} 
    \bigoplus_{i\ne j} \C_{(e_i-e_j)} \oplus \C_{(-e_j)},&j\le k-1\\
    \bigoplus_{i\ne k} \C_{(e_i)},&j=k \\
  \end{cases}
\end{equation*}

We now consider the diagonal action of $P\Uk$ on $(\cp^{k-1})^n$, and
hence of $T\subset P\Uk$. The fixed points under the $T$-action
are simply elements of the $n$-fold product of the fixed points in
$\cp^{k-1}$. Thus they correspond to partitions
\begin{equation*}
\{1\h..n\} = I_1 \sqcup I_2 \sqcup \h.. \sqcup I_k
\end{equation*}
in the obvious way. We denote such a fixed point by
\begin{equation*}
F_{I_1,\h..,I_k} \in (\cp^{k-1})^n
\end{equation*}

\subsection{Calculations on $(\cp^2)^n$}

We now specialize to $\cp^2$, setting $X=(\cp^2)^n$. We will calculate
invariants of $X\d/T(0)$ and $X\d/P\U3$, for $n$ not a multiple of $3$.

\begin{figure}[!htbp]\label{fig:cp24}
\begin{center}
\includegraphics{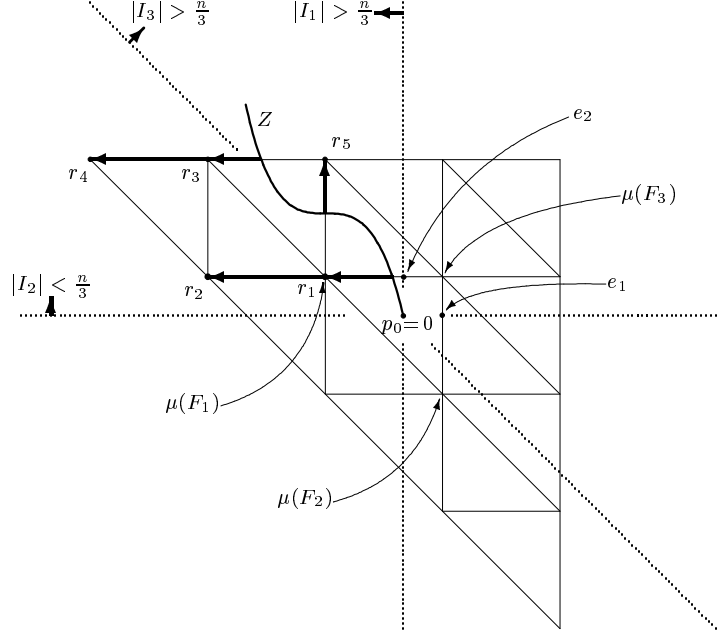}
\end{center}
\caption{The moment map for $(\cp^2)^4$, showing the transverse
paths used in the calculation, with their wall-crossings.}
\end{figure}

The fixed points correspond to
partitions
\begin{equation*}
\{1\h..n\} = I_1 \sqcup I_2 \sqcup I_3.
\end{equation*}
For $n=1$ we have
\begin{equation*}
\begin{split}
\mu(F_1) &= e_2 - 2e_1 \\
\mu(F_2) &= e_1 - 2e_2 \\
\mu(F_3) &= e_1 + e_2 \\
\end{split}
\end{equation*}
Hence, in general
\begin{equation*}
\mu(F_{I_1,I_2,I_3})= (-2|I_1|+|I_2|+|I_3|)e_1 + (|I_1|-2|I_2|+|I_3|)e_2.
\end{equation*}
It follows that $0$ is a regular value as long as $n$ is not a multiple of
$3$.

We start by taking the path $Z$, as depicted in figure
\ref{fig:cp24}. In the case of $(\cp^2)^4$, $Z$ has $3$
wall-crossings. The horizontal walls (in the figure) correspond to the subgroup
\begin{equation*}
H_2=\{\begin{pmatrix}
 1 \\
 &e^{it} \\
\end{pmatrix}\mid t\in\R\}
\end{equation*}
$Z$ crosses these walls in in the same direction as $e_2$, so that the
isomorphism  $H_2\xra{\iso} S^1$ is the standard one for $H_2$, as
described above. On the other hand, the vertical walls correspond to
the subgroup $H_1$, and the direction of the crossing by $Z$
corresponds to the oriented subgroup $\bar{H}_1$.

Let $\flag_1$ denote the flag
\begin{equation*}
\flag_1 = (H_2, H_2 \ti \bar{H}_1)
\end{equation*}
and let $\flag_2$ denote the flag
\begin{equation*}
\flag_2 = (\bar{H}_1, \bar{H}_1 \ti H_2 )
\end{equation*}
We then have, in the case $n=4$,
\begin{equation*}
p_0 \sim (r_1,\flag_1) + (r_2,\flag_1) + (r_3,\flag_1) + (r_4,\flag_1)
  + (r_5,\flag_2)
\end{equation*}

In general, let $R_1$ denote the set of vertices corresponding
to fixed points $F_{I_1,I_2,I_3}$ with $|I_1|>\frac{n}{3}$ and
$|I_3|>\frac{n}{3}$, and $R_2$ the vertices corresponding
to fixed points $F_{I_1,I_2,I_3}$ with $|I_2|<\frac{n}{3}$ and
$|I_3|<\frac{n}{3}$. We then have
\begin{equation}\label{eq:cp2tofp}
p_0 \sim \sum_{r\in R_1} (r,\flag_1) + \sum_{r\in R_2} (r,\flag_2).
\end{equation}

Fixing attention on the point $F_{I_1,I_2,I_3}$, we now calculate the
maps
\begin{equation*}
\loc_{\flag_i} : \HH^*_T(\pt)=\Q[u_1,u_2] \to \HH^*(\pt)=\Q
\end{equation*}
where $u_1$ and $u_2$ are the generators corresponding to $H_1$ and
$H_2$. We have
\begin{equation*}
V := \nu F_{I_1,I_2,I_3} \iso \C^{|I_1|}_{(e_2-e_1)} \oplus
\C^{|I_1|}_{(-e_1)} \oplus \C^{|I_2|}_{(e_1-e_2)} \oplus 
\C^{|I_2|}_{(-e_2)} \oplus \C^{|I_3|}_{(e_1)} \oplus 
\C^{|I_3|}_{(e_2)}.
\end{equation*}
The subbundle stabilized by $H_2$ is
\begin{equation*}
V^{H_2} = \C^{|I_1|}_{(-e_1)} \oplus \C^{|I_3|}_{(e_1)}.
\end{equation*}
Hence, to calculate $\loc_{\flag_1}$ we need the equivariant weighted
Segre classes of  $V/V^{H_2}$.
\begin{equation*}
V/V^{H_2} \iso \C^{|I_1|}_{(e_2-e_1)} \oplus
\C^{|I_2|}_{(e_1-e_2)} \oplus \C^{|I_2|}_{(-e_2)} \oplus
\C^{|I_3|}_{(e_2)}.
\end{equation*}
Now the weighted Chern class
\begin{equation*}
c_1^{H_1}(\C_{(ke_1 + le_2)}) = ku_1
\end{equation*}
and hence
\begin{equation*}
s^w_{H_1}(\C_{(ke_1 + le_2)}) = (l+ku_1)^{-1}
\end{equation*}
Therefore
\begin{equation*}
\begin{split}
s^w_{H_1}(V/V^{H_2})
  &= (1-u_1)^{-|I_1|} (u_1-1)^{-|I_2|} (-1)^{-|I_2|} (1)^{-|I_3|} \\
  &= (1-u_1)^{-(|I_1|+|I_2|)} \\
\end{split}
\end{equation*}
and $\rank(V/V^{H_2})=|I_1|+2|I_2|+|I_3|$. Hence, setting
$k=\rank(V/V^{H_2})$, and $l=|I_1|+|I_2|$,
\begin{equation*}
\loc(H_2,H_1,V/V^{H_2}) : u_2^j \mapsto 
  \begin{cases}
    0, &j<k-1 \\
    \choose{j+1-k+l-1}{l-1}
      u_1^{j+1-k}, &j\ge k -1\\
  \end{cases}
\end{equation*}
By the functorial properties of integration over the fibre,
this map commutes with multiplication by $u_1^{j_1}$. To get
$\loc_{\flag_1}$ we must compose with $\loc(\bar{H}_1,1,V^{H_2})$,
which is equal to $-\loc(H_1,1,V^{H_2})$. We have
\begin{equation*}
\loc(H_1,1,V^{H_2}): u_1^j \mapsto 
  \begin{cases}
    (-1)^{|I_1|}, &j=|I_1|+|I_3|-1 \\
    0 &\text{otherwise}\\
  \end{cases}
\end{equation*}
Hence
\begin{equation*}
\loc_{\flag_1}(u_1^{j_1}u_2^{j_2}) = 
  \begin{cases}
    (-1)^{|I_1|+1}\choose{j_2-|I_2|-|I_3|}{|I_1|+|I_2|-1},
  &\substack{j_2\ge|I_1|+2|I_2|+|I_3|-1\\\text{and } j_1+j_2=2m-2} \\
    0 &\text{otherwise.}\\
  \end{cases}
\end{equation*}
And similarly
\begin{equation*}
\loc_{\flag_2}(u_1^{j_1}u_2^{j_2}) = 
  \begin{cases}
    (-1)^{|I_2|+1}\choose{j_1-|I_1|-|I_3|}{|I_1|+|I_2|-1},
  &\substack{j_1\ge2|I_1|+|I_2|+|I_3|-1\\\text{and } j_1+j_2=2m-2} \\
    0 &\text{otherwise.}\\
  \end{cases}
\end{equation*}

\subsection{The Volume}

We can now easily write down formul{\ae} for the evaluation of cohomology
classes on $(\cp^2)^n\d/T$. And, by applying the integration formula from
the companion paper~\cite{skm:g-t}
relating evaluation of classes on
$G$-symplectic-quotients to evaluation on $T$-symplectic-quotients, we
can write down formul{\ae} for the evaluation of classes on
$(\cp^2)^n\d/P\U3$. As an example, we will give a formula for the
volume of $(\cp^2)^n\d/P\U3$.

As usual, the `prequantum line bundle' $\CL\to(\cp^2)^n$, which has first
Chern class equal to $\left[\frac{\om}{2\pi}\right]$, descends to a
line bundle over the symplectic quotient, which we also denote by
$\CL$. The dimension of $(\cp^2)^n\d/P\U3$ is $4n-16$, and hence the
volume is equal to the evaluation of the class $\frac{(2\pi
c_1(\CL))^{2n-8}}{(2n-8)!}$
against the fundamental class.

We define
\begin{equation*}
a := \frac{1}{6(2n-8)!} (2\pi c_1^{S^1}(\CL))^{2n-8} \smile 
  \prod_{\al\in\roots} c_1^{S^1}(\C_{(\al)}),
\end{equation*}
where $\roots=\{e_1-e_2, e_1, e_2, e_2-e_1, -e_1, -e_2\}$ is the set
of roots.
Then we have
\begin{equation*}
\vol((\cp^2)^n\d/P\U3) = \ig{a}{(\cp^2)^n\d/T(0)}
\end{equation*}
and, applying Theorem~D and equation
\eqref{eq:cp2tofp},
\begin{equation*}
\vol((\cp^2)^n\d/P\U3)
  =\sum_{\substack{(I_1,I_2,I_3)\\|I_1|>\frac{n}{3},|I_3|>\frac{n}{3}}}
    \loc_{\flag_1}(\eval{a}_{F_{I_1,I_2,I_3}})
  + \sum_{\substack{(I_1,I_2,I_3)\\|I_2|<\frac{n}{3},|I_3|<\frac{n}{3}}}
    \loc_{\flag_2}(\eval{a}_{F_{I_1,I_2,I_3}})
\end{equation*}
where the triples $(I_1,I_2,I_3)$ run through partitions of
$\{1\h..n\}$. Given such a partition, set
$i_1=|I_1|,i_2=|I_2|,i_3=|I_3|$. Then 
\begin{equation*}
\begin{split}
\eval{\CL}_{F_{I_1,I_2,I_3}} &\iso\C_{(\mu({F_{I_1,I_2,I_3}}))} \\
  &=\C_{((-2i_1+i_2+i_3)e_1+(i_1-2i_2+i_3)e_2)}.
\end{split}
\end{equation*}
Hence
\begin{equation*}
\eval{c_1^{S^1}(\CL)}_{F_{I_1,I_2,I_3}} = 
  (-2i_1+i_2+i_3)u_1+(i_1-2i_2+i_3)u_2,
\end{equation*}
so that
\begin{equation*}
\eval{a}_{F_{I_1,I_2,I_3}}= 
((-2i_1+i_2+i_3)u_1+(i_1-2i_2+i_3)u_2)^{2n-8}\c.
  (2u_1^3u_2^3-u_1^4u_2^2-u_1^2u_2^4).
\end{equation*}
Applying the formul{\ae} for $\loc_{\flag_1}$ and $\loc_{\flag_2}$, and
using the identity
$2\choose{m}{k}-\choose{m+1}{k}-\choose{m-1}{k}=-\choose{m-1}{k-2}$, we
easily derive the unilluminating but nontheless computable formula
\begin{multline*}
\frac{(2n-8)!}{(2\pi)^{2n-8}}\vol((\cp^2)^n\d/P\U3) = \\
\sum_{\substack{i_1>\frac{n}{3},i_3>\frac{n}{3}\\i_1+i_3\le n}}
  \frac{n!(-1)^{i_1+1}}{i_1!i_3!(n-i_1-i_3)!}\c. \\
  \( \choose{2n-8}{i_1+i_3-4}(n-3i_1)^{i_1+i_3-4}
     (3i_1+3i_3-n)^{2n-4-i_1-i_3}(2+i_3-n) - \right.\\
     \choose{2n-8}{i_1+i_3-3}(n-3i_1)^{i_1+i_3-3}
     (3i_1+3i_3-n)^{2n-5-i_1-i_3} - \\
     \left. \sum_{j=0}^{i_1+i_3-5} \choose{n+i_1-6-j}{n-i_3-3}
       \choose{2n-8}{j}(n-3i_1)^j(3i_1+3i_3-n)^{2n-8-j}     \) \\
+ \sum_{\substack{i_1<\frac{n}{3},i_2<\frac{n}{3}\\i_1+i_2\le n}}
  \frac{n!(-1)^{i_1+1}}{i_1!i_2!(n-i_1-i_2)!}\c. \\
  \( \choose{2n-8}{i_1+i_2-4}(n-3i_1)^{i_1+i_2-4}
     (3i_1+3i_2-n)^{2n-4-i_1-i_2}(2+i_2-n) - \right.\\
     \choose{2n-8}{i_1+i_2-3}(n-3i_1)^{i_1+i_2-3}
     (3i_1+3i_2-n)^{2n-5-i_1-i_2} - \\
     \left.\sum_{j=0}^{i_1+i_2-5} \choose{n+i_1-6-j}{n-i_2-3}
       \choose{2n-8}{j}(n-3i_1)^j(3i_1+3i_2-n)^{2n-8-j}     \). \\
\end{multline*}

\appendix

\section{Orbifolds, orbifold-fibre-bundles, and integration over the
fibre}\label{app:orb}

\spre{The purpose of this appendix is to collect together a number of facts
about orbifolds which we use in the paper.  These are all
straightforward generalizations of standard results.}

An orbifold is a generalization of a manifold, and can roughly be
thought of as follows: whereas an $n$-dimensional manifold is locally
modelled on $\R^n$, an $n$-dimensional orbifold is locally modelled on
the quotient of $\R^n$ by a finite group.  Orbifolds were first
defined and studied by Satake in his announcement~\cite{sat:gen-not}
and his paper~\cite{sat:gau-bon} (Satake used the term `$V$-manifold';
the term `orbifold' is due to Thurston).  Our interest in orbifolds
comes from the fact that the wall-crossing-cobordism and its boundary
are in general orbifolds (even if we are interested in a symplectic
quotient which is smooth, we may encounter orbifold singularities
after crossing a wall).

In this appendix we collect together facts involving orbifolds which
we need in the rest of the paper. These facts all involve integration on
orbifolds, in one form or another, and can be seen as straightforward
generalizations of standard facts involving manifolds. These
generalizations exist because an orbifold is a `rational (co)homology
manifold', which basically means that, if we take rational
coefficients, it possesses the same homological and cohomological
properties as a manifold.

We begin by giving Satake's definition of an orbifold, as well as his
generalizations to oriented and symplectic orbifolds. We then state
the various facts involving orbifolds, and indicate how these facts
follow from results in the literature.

\subsection{The definition of an orbifold}

We now give Satake's definitions. We do this to set up notation which
we refer to in the rest of the appendix, but also to make explicit
some of the subtleties in the definition. These subtleties are
necessary for orbifolds to have the good properties that we need (such
as a rational fundamental class).

\begin{definition}[Satake~\cite{sat:gen-not,sat:gau-bon}]
Let $M$ be a Hausdorff topological space. A ($C^\infty$)
\sdef{orbifold structure} on $M$ consists of 
a covering $\CU$ of $M$ by open sets, and for each open set
$U\in\CU$, an associated triple $(\tU, G_U, \vp_U)$, where
\begin{description}
\item[$\tU$] is a connected open subset of $\R^n$;
\item[$G_U$] is a finite group of linear transformations mapping $\tU$
to itself, such that the set of points fixed by $G_U$ has codimension
$\ge 2$; and
\item[$\vp_U$] is a continuous map $\tU\to U$ such that, for every
$x\in\tU$ and $g\in G_U$, $\vp_U(gx)=\vp_U(x)$. We assume that the
induced map $G_U\!\!\setminus\!\! \tU \to U$ is a homeomorphism. 
\end{description}
Moreover, if $U, V\in\CU$ are open sets such that $U\subset
V$, then we are given an injective group homomorphism
$\bt_{UV}:G_U\hra G_V$, and an inclusion $i_{UV}:\tU\hra\tilde V$
which is a diffeomorphism onto its image, and which is equivariant with respect to
the action of $G_U$ (and its image in $G_V$), and such that
$\vp_U=\vp_V\o i_{UV}$.  Finally, we assume that the open sets in
$\CU$ form a basis for the topology of $M$. (It is fairly standard to
refer to $\tU$ as a local cover, $G_U$ as a local group, and $\vp_U$
as a local covering map.)

An \sdef{orbifold}, then, is a space $M$ together with an equivalence
class of orbifold structures on $M$ (see Satake~\cite{sat:gen-not} for
details of the straightforward notion of when two such sets of data
define the same orbifold structure).

By enhancing the definition of an orbifold structure, we can define an
\sdef{oriented orbifold}: we ask that each $\tU$ be given an
orientation which is preserved by the action of the group $G_U$, and
that such orientations be compatible with the inclusions
$i_{UV}:\tU\hra \tilde V$.

Similarly, we define a \sdef{symplectic orbifold} by asking that each
$\tU$ be given a symplectic form, with the same invariance and
compatibility conditions.
\end{definition}

\begin{definition}
A point $x$ of an orbifold $M$ is a \sdef{smooth point} if there
exists some open set $U\in\CU$ containing $x$, and such that the
associated group $G_U$ is the trivial group.  The set of points which
are not smooth points are called \sdef{singular points}.
\end{definition}

\begin{remark}\label{rem:smooth-connected}
The set of smooth points of an orbifold $M$ is connected (within each
component of $M$).  More precisely, given any open set $U\in\CU$ with
associated triple $(\tU, G_U, \vp_U)$, then the set of singular
points in $U$ is the image, under $\vp_U$, of a finite union of
submanifolds of $\tU$ having codimension~$\ge 2$. Each of these
submanifolds is the submanifold of points fixed by some nontrivial
element $g\in G_U$. (A straightforward argument by contradiction shows
that the codimension restriction on the fixed points of each local
group $G_U$ implies the same restriction for each nontrivial subgroup of
$G_U$, and hence for each nontrivial $g\in G_U$).
\end{remark}

\subsection{The fundamental class of an oriented orbifold}

\begin{fact}\label{fact:int}
Let $M$ be an $n$-dimensional compact oriented orbifold (without
boundary). Then the orientation defines a rational fundamental class
$[M]\in\HH_n(M)$ (recall that we are taking homology and cohomology
with rational coefficients throughout this paper).  Moreover, $M$
satisfies rational Poincar\'e duality, which can be expressed as the
fact that the pairing
$\HH^i(M)\ti\HH^{n-i}(M) \to \Q$ given by
$(a,b) \mapsto \ig{a\smile b}{M}$
is a dual pairing on the rational cohomology of $M$.
\end{fact}

The relationship between the orientation and the fundamental class is
as follows. At any smooth point $x\in M$, we use the orientation to
define a generator $1_x \in \HH_n(M,M\setminus\{x\})$ via the
identification with $\HH_n(\R^n,\R^n\setminus\{0\})\iso\Q$ given by
excision (using an oriented chart). Then the fundamental class is
the unique class $[M]\in\HH_n(M)$ whose image under the natural map 
$\HH_n(M) \to \HH_n(M,M\setminus\{x\})$ has image $1_x$, for each
smooth point $x$. (Since the set of smooth points is connected, we
actually only need to use one smooth point for each
component of $M$ to get the right normalization.)

\begin{proof}[Sketch of proof]
There are two different approaches to the proof. Satake's approach
\cite{sat:gen-not,sat:gau-bon} is to define an orbifold version of the
de Rham complex%
\footnote{\label{foot:odf}
A differential form on an orbifold $M$ is a collection of
differential forms on the sets $\tU$, invariant under the local
groups $G_U$, and compatible with the inclusion maps in the
obvious way; integration is defined using a partition of unity and
adding up integrals on sets $\tU$ multiplied by the factors
$1/|G_U|$.}  and to prove de Rham's theorem: that the orbifold de Rham
cohomology is canonically isomorphic to the singular cohomology of $M$
(with real coefficients). The fundamental class is then defined in
terms of integration.

The other approach is to use the notion of a `rational homology
manifold', as described by Borel in
\cite[chapters~I--II]{bor:sem-tra}%
\footnote{A \sdef{rational homology $n$-manifold} is a space whose
local homology, with rational coefficients, agrees with that of an
$n$-manifold (where the local homology at $x\in M$ is
$\HH_*(M,M\setminus\{x\})$. It's an easy calculation to show that an
orbifold is a rational homology manifold.  The construction of the
rational fundamental class of a rational homology manifold mimcs the
usual construction: one shows that an orientation gives a constant
section of the local homology sheaf, and then applies a Mayer-Vietoris
patching argument,(as in \cite[section 5]{bot-tu:dif-for} or
\cite[section 6.3]{spa:alg-top}).}.
An orbifold is a rational homology manifold,
and Borel shows how various properties of the homology of manifolds go
over to rational homology manifolds, including the existence of a
(rational) fundamental class and (rational) Poincar\'e duality.
\end{proof}

\subsection{Oriented orbifolds with boundary and Stokes's theorem}

Satake defines an orbifold-with-boundary in \cite[section
3.4]{sat:gau-bon}. His definition is equivalent to modifying the
definition of orbifold by allowing the open covers $\tU$ to be
open subsets of $\R^n$ or of the halfspace $\R^{n-1}\ti[0,\infty)$
(but keeping the same conditions with respect to $G_U$ and $\vp_U$). We
then have

\begin{fact}\label{fact:orb-bd}
Let $M$ be an $n$-dimensional compact oriented
orbifold-with-boundary. Then the boundary $\bd M$ is an
$(n-1)$-dimensional orbifold, with a natural orientation induced from
the orientation of $M$, and $\bd M$ is null-homologous in $M$ (that
is, the image of the fundamental class $[\bd M]$ is zero in
$\HH_{n-1}(M)$).
\end{fact}

\begin{proof}[Sketch of proof]
In the language of differential forms, this is just Stokes's theorem,
and the standard local argument applies (e.g.~\cite[theorem
3.5]{bot-tu:dif-for}). Alternatively, using the rational (co)homology
manifold approach, this fact follows from
Poincar\'e-Lefschetz duality~\cite[chapter II]{bor:sem-tra}.
\end{proof}

\subsection{Orbibundles and integration over the fibre}

An `orbibundle' is the natural orbifold version of a fibre
bundle. Satake defined orbibundles (he called them $V$-bundles). 

\begin{definition}[Satake~\cite{sat:gau-bon}]
Let $M$ be an orbifold, with orbifold structure defined by the open
cover $\CU$.  An \sdef{orbibundle} over $M$ is defined by giving, for
each open set $U\in\CU$ (with associated triple $(\tU,G_U,\vp_U)$) a
$G_U$-equivariant fibre bundle $\tilde E \to \tU$. (Each inclusion
$i_{UV}$ must lift to a $G_U$-equivariant bundle map, which is an
isomorphism on the fibres.)  Given an orbibundle over $M$, there is an
associated topological space (which we will refer to as the
\sdef{total space}) $E$ with a map $E\xra{\pi} M$ defined so that
$\pi\inv(U)=\tilde E/G_U$. An orbibundle is \sdef{oriented} if the
fibres of each bundle $\tilde E\to\tU$ are oriented (these
orientations must be preserved by the local groups $G_U$ and 
compatible with inclusion maps).
\end{definition}

\begin{remarks}
\begin{enumerate}
\item Although the fibre of an orbibundle may be any space, in
our applications the fibre will always be an orbifold.
\item The total space $E$ of an orbibundle $E\xra{\pi} M$ is not in
general a fibre bundle: if $x$ is a smooth point of $M$ then
$\pi\inv(x)$ will be a copy of the fibre $F$, but if $x$ is an
orbifold point of $M$, then $\pi\inv(x)$ may be the quotient of $F$ by
a finite group.
\end{enumerate}
\end{remarks}

We will describe the properties of a map on cohomology known as
`integration over the fibre', but in order to do this, we must define
the notion of a suborbifold.

\begin{definition}
Given an orbifold $M$, then a 
\sdef{suborbifold} $M'$ of $M$ is defined by giving a submanifold of
each $\tU$, stable under $G_U$ and compatible with the inclusion maps,
and such that the restriction of the orbifold structure on $M$ defines
an orbifold structure on $M'$ (in particular, in each submanifold the
set of points fixed by $G_U$ should have codimension~$\ge~2$).
\end{definition}

It is important to note that, with this definition, a suborbifold $M'$
of $M$ consists mainly of smooth points of $M$ (more precisely, those
points of $M'$ which are smooth in $M$ make up a dense open subset of
$M'$). This is consistent with Satake's definition of an orbifold,
which forces most points of an orbifold to be smooth
points\footnote{It would be possible to give an alternative definition
of an orbifold which removed these restrictions. Specifically, given a
local triple $(\tU, G_U, \vp_U)$, we could remove the restriction that
the set of points fixed by the $G_U$-action on $\tU$ have codimension
$\ge 2$, and alter the rest of the definition in a compatible manner.
This alternative definition would be more natural in some respects,
but it would also be more involved, since we would then need to take
into account various numerical factors.}.

\begin{fact}\label{fact:int-fib}
Let $E\xra{\pi} M$ be an orbibundle, with fibre the compact
oriented orbifold $F$. Then there is a map
\begin{equation*}
\pi_* : \HH^*(E) \to \HH^{*-\dim F}(M)
\end{equation*}
known as \sdef{integration over the fibre}%
\footnote{often referred to as the \sdef{Gysin map} (it generalizes
the Gysin map defined for a sphere bundle) or, in a more general setting, the
\sdef{pushforward}.}
having the following properties:
\begin{enumerate}
\item Integration over the fibre is a module homomorphism of
$\HH^*(M)$-modules (the module structure is given by pullback via
$\pi$ followed by cup product). This is equivalent to the `push-pull
formula'
\begin{equation*} \pi_*(\pi^*(a)\smile  b) =
  a\smile \pi_*(b), \qquad\forall a\in\HH^*(E),b\in\HH^*(M).
\end{equation*}
\item Let $i:M'\hra M$ be the inclusion of a suborbifold
of $M$ , and let $E'\xra{\pi'} M'$ denote the
orbibundle over $M'$ defined by the restriction of $E$. Then
the following square commutes:
\begin{equation*}
\xymatrix{\HH^*(E')  \ar[d]^{\pi'_*} & %
  \HH^*(E) \ar[d]^{\pi_*} \ar[l]_{\tilde{i}{}^*} \\
\HH^{*-\dim F}(M')  & \HH^{*-\dim F}(M). \ar[l]_{i^*}}
\end{equation*}
(where $\tilde i:E'\to E$ is the lift of $i$).
\item If $E, M$ and $F$ are compact oriented orbifolds, and the
orientation of $E$ equals the product of the orientations of $M$ and
$F$, then for any class $a\in\HH^*(E)$ we have
\begin{equation*}
\int_E a = \int_M \pi_*(a).
\end{equation*}
\end{enumerate}
\end{fact}

\begin{proof}[Sketch of proof]
We again indicate two different proofs. Using differential forms, the
usual formula for integration over the fibre is well-defined on the
local bundles $\tilde E\to\tU$ (this was defined by 
Lichnerowicz~\cite{lic:the-hom}, and is also explained by Bott and
Tu~\cite[p. 61]{bot-tu:dif-for}; of course we are using fact
\ref{fact:int}, allowing us to integrate over the orbifold fibres). It
is easy to check that this gives $G_U$-invariant differential forms on
the sets $\tU$, and hence orbifold differential forms on $M$ (see
fotnote \ref{foot:odf}). The advantage of this approach is that the
three properties we have listed above follow immediately from the
definition.

Alternatively, in the manifold case, integration over the fibre can be
defined using the Leray-Serre spectral sequence of the fibration
(described for sphere bundles quite explicitly in Bott and Tu
\cite[pp. 177--179]{bot-tu:dif-for}). For an orbibundle $E\xra{\pi}M$
we use the Leray spectral sequence (with rational coefficients) of the
map $\pi$ \cite[pp. 179--182]{bot-tu:dif-for}, trivializing the the
top cohomology sheaf of the fibres by the rational fundamental classes
on the local covers $\tilde E\to\tU$. Finally, the algebraic and
naturality properties of the Leray-Serre spectral sequence which imply
properties 1-3 above also carry over to the Leray spectral sequence
(see e.g. McCleary \cite{mcc:use-gui}).
\end{proof}

\begin{remark}\label{rem:int-fib-cw}
We also need a related result concerning integration over the fibre:
this time for an (honest) fibre bundle $E\xra{\pi} B$, with fibre an
oriented orbifold $F$, but where the base space $B$ may be any
CW-complex. Using the same arguments as above, it is easy to show that
integration over the fibre $\pi_*$ is well-defined for
such bundles, and satisifes properties 1 and 2.
\end{remark}

\subsection{How orbifold-fibre-bundles can arise as locally free
quotients of manifolds}

\begin{fact}
Suppose the compact connected Lie group $G$ acts on a compact oriented
manifold $N$ with a locally free action (that is, the stabilizer
subgroup of each point is finite). Then the quotient space $N/G$ can
be given an oriented orbifold structure (the orientation 
is fixed by orienting $G$).
\end{fact}

This orbifold structure on $N/G$ is constructed by taking local slices
for the action (for the existence and properties of local slices, see
for example Bredon~\cite[Chapter IV]{bre:int-com},
Kawakubo~\cite[section 4.4]{kaw:the-tra}, or the chapter by
Palais~\cite[chapter VIII]{bor:sem-tra}). Specifically, given a point
$x\in N$, then there exists a \sdef{linear slice} for the $G$-action
at $x$: a submanifold $S\subset N$ which is transverse to the
$G$-orbits, is mapped to itself by the stabilizer subgroup $G_x$, and
is equivariantly identified with an open subset of $\R^n$ with respect
to a linear action of $G_x$ on $\R^n$. Letting $F$ denote the subgroup
of $G_x$ which fixes every point in $S$, then the triple
$(S,G_x/F,\vp)$ defines the orbifold structure at $[x]\in N/G$ (where
$\vp$ maps $S$ to $S\c.G\subset N/G$).

The following existence facts follow easily from the definition of
orbifold together with simple arguments involving local slices.

\begin{facts}
\begin{enumerate}
\item If $N'$ is an oriented submanifold of $N$, stable under $G$ and
transverse to the submanifolds $N^H$, for each finite subgroup
$H\subset G$, then $N'/G$ is a suborbifold of $N$.
\item If $N$ is an oriented manifold-with-boundary on which the
compact connected Lie group $G$ acts, with a locally free action, then $N/G$
is an oriented orbifold-with-boundary.
\item Suppose $E$ and $N$ are oriented manifolds, and $E\to N$ is a fibre
bundle. Then  if $G$ and $H$ are compact connected Lie groups, and
$G\ti H$ acts on $E$, covering an action of $H$ on $N$, and these
actions are locally free, then $E/(G\ti H)\to N/H$ is an orbibundle.
\end{enumerate}
\end{facts}

\section{Cohomology and integration formulae for weighted projective
bundles}\label{app:wtd-proj-b}

The purpose of this appendix is to give generalizations of two
classical formulae concerning projective bundles. Let $Y$ be a
CW-complex, let $V\to Y$ be a complex vector bundle, and let $\P(V)\to
Y$ be its projectivization \cite[p. 269]{bot-tu:dif-for}.  The first
classical formula describes the cohomology of $\P(V)$, and the second
(and perhaps less well-known) calculates integrals over the fibres of the
bundle $\P(V)\to Y$.

The generalizations we give apply to bundles constructed as
follows. Let $V\to Y$ be a complex vector bundle, and suppose $S^1$
acts on $V$, such that the action is linear on the fibres of $V$ (that
is, the action covers the trivial action on $Y$), and such that the set
of fixed points equals the zero section. We consider the bundle
$S(V)/S^1\xra{\pi} Y$, where $S(V)$ denotes the unit sphere bundle in
$V$, relative to some invariant metric.

These bundles can be considered as generalizations of projective
bundles in the following sense. If $S^1$ acts with `weight one' on the
fibres (i.e.\ the standard multiplication action of $S^1\subset
\C^\ti$), then each $S^1$-orbit lies in precisely one line in $V$, and
identifying $S^1$-orbits with lines induces a isomorphism
$S(V)/S^1\iso\P(V)$. The general case that we consider allows any
combination of positive and negative weights. This general case
includes `weighted projectivizations' 
which correspond to $S^1$ actions having only
positive weights (Kawasaki calculates the cohomology of weighted
projective spaces in \cite{kaw:coh-twi}; 
for some definitions and results in algebraic
geometry on weighted projective spaces, see \cite{dol:wei-pro}).

We begin by reviewing the cohomology and integration formulae in the
case of projective bundles. We then state and prove the general
cohomology formula, followed by the general integration formula.
Finally, using the homotopy quotient construction, we will observe
that all the definitions, formulae, and proofs naturally extend to the
case in which an auxilliary group $G$ acts on $V$ and $Y$, commuting
with the $S^1$-action and with the projection.

\subsection{Projective bundles}

The projectivization $\P(V)$ possesses a distinguished cohomology
class $h\in\HH^2(\P(V))$, which is usually defined as follows. Let
$S\to\P(V)$ denote the \sdef{tautological line bundle} (where the
fibre of $S$ over a point is just the corresponding line in $V$), and
define $h$ to be the first Chern class of the dual line bundle,
$h=c_1(S^*)$.  

Then the cohomology of $\P(V)$ is given by the formula\footnote{Bott
and Tu \cite[pp. 269-271]{bot-tu:dif-for} describe the
projectivization and the tautological line bundle, and following
Grothedieck, they define the Chern classes in terms of the cohomology
formula.}
\begin{equation*}
\HH^*(\P(V)) \iso \frac{\HH^*(Y)[h]}%
{\<c_0(V)h^r + c_1(V)h^{r-1} + \h.. + c_r(V)\>}.
\end{equation*}
where $c_i(V)\in\HH^{2i}(Y)$ is the $i$-th Chern class, and
$r=\rank(V)$. In this formula the product $ah^i$ (where
$a\in\HH^*(Y)$) is identified with the class
$(\pi^*a)h^i\in\HH^*(\P(V))$.

The vector bundle $V\to Y$ has associated  Segre classes
$s_i(V)\in\HH^{2i}(Y)$. The total Chern class and the total Segre
class are multiplicative inverses to each other (in the cohomology
ring of $Y$), that is
\begin{equation*}
c(V)s(V) = 1, 
\end{equation*}
and this can be used to define the Segre classes in terms of the Chern
classes. (As an example, consider the tautological line bundle
$S\to\cp^n$. Then $c(S)=1-h$, where $h$ is the generator of
$\HH^*(\cp^n)$, and $s(S)=(1-h)\inv=1+h+h^2+\h..+h^n$.)

The integration formula expresses integrals over the fibres in terms
of Segre classes:
\begin{equation*}
\pi_*(h^i) = \left\{\begin{array}{ll}
0 & i<\rank(V)-1,\\
s_{i-\rank(V)+1}(V),\qquad & i\ge \rank(V)-1,
\end{array} \right.
\end{equation*}
where $\pi_*$ denotes integration over the fibre (see fact 
\ref{fact:int-fib}). (This formula is sufficient to
calculate the integral over the fibres of any class on $\P(V)$, since
every class can be expressed in the form $(\pi^*a)h^i$, and we have
$\pi_*((\pi^*a)h^i)=a\pi_*(h^i)$.)\footnote{The integration formula
might appear to be overkill: since it follows from the cohomology
formula that every class on $\P(V)$ can be expressed as $(\pi^*a)h^i$
for $0\le i \le \rank(V)-1$, in fact we only need to observe that
$\pi_*(h^i)=0$ for $0\le i \le \rank(V)-1$, and
$\pi_*(h^{\rank(V)-1})=1$. However in applications we are often given
a class on $\P(V)$ expressed as $(\pi^*a)h^i$ where $i$ is not
necessarily in this range.  Using the cohomology formula, we could
rewrite such a class in terms of the cohomology of $Y$ and the classes
$\{1,h,h..,h^{\rank(V)-1}\}$, in which case the integral over the
fibres would be the coefficient of $h^{\rank(V)-1}$. The integration
formula is simply the answer one gets by following this process.}

\subsection{Weighted Chern classes and the cohomology formula}

We now return to the general case: $V\to Y$ is a complex vector
bundle, with an action of $S^1$ on $V$, covering the trivial action on
$Y$, and such that the set of fixed points equals the zero section.

We first define the weighted Chern class of the pair $(V,S^1)$
(although we will sometimes abuse notation and simply refer to this as
the weighted Chern class of $V$). We will
then state and prove a formula for the cohomology of the total space
of the bundle $S(V)/S^1\xra{\pi} Y$ (where $S(V)$ denotes the unit
sphere bundle in $V$ relative to some invariant metric).

\begin{definition} \label{def:wtd-c-c}
The quick definition of the weighted Chern class is
this: the weighted Chern class $c^w$ is multiplicative under direct 
sum of bundles, and commutes with pullbacks (so that the splitting
principle applies), and for a line bundle
$L$ acted on with weight $i$, is given by $c^w(L)=i+c_1(L)$ (where
$c_1(L)$ is the regular first Chern class).

Explicitly, under the $S^1$
action, $V$ splits into `isotypic' subbundles
\begin{equation*}
V\iso \bigoplus_{i\in\Z} V_i,
\end{equation*}
where $S^1$ acts with weight $i$ on $V_i$ (that is, $\lb\in S^1$ acts
on $V_i$ by multiplying the fibre coordinates by $\lb^{i}$).  Then the
\sdef{weighted Chern class} of $(V,S^1)$, which we denote
$c^w(V)\in\HH^*(Y)$, is the product
\begin{equation*}
c^w(V) := \prod_i c^w(V_i),
\end{equation*}
where, setting  $r$  equal to the rank of $V_i$,
\begin{equation*}
c^w(V_i) = i^{r} + i^{r-1}c_1(V_i) + i^{r-2}c_2(V_i) + \h.. + c_{r}(V_i)
\end{equation*}
(here $c_j(V_i)$ is the regular $j$-th Chern class of $V_i$).
It follows from the properties of the regular Chern class
that the weighed Chern class is natural with respect
to pullbacks, and multiplicative with respect to direct sum (it is
easiest to think of the $S^1$-action as simply decomposing $V$ into a
direct sum of subbundles, each of which is labelled with an integer,
and to note that this decomposition commutes with pullback and direct
sum in an obvious way).
\end{definition}

\begin{proposition}\label{pr:c-wtd-p-b}
Let $V\to Y$ be a complex vector bundle with an action of $S^1$ as
above. Define $h\in\HH^2(S(V)/S^1)$ to be the first Chern class of the
principal orbifold bundle $S(V)\to S(V)/S^1$ (see remark
\ref{rem:define-h} below). Then there is a ring isomorphism
\begin{equation*}
\HH^*(S(V)/S^1) \iso \frac{\HH^*(Y)[h]}%
{\<c_0^w(V)h^r + c_1^w(V)h^{r-1} + \h.. + c_r^w(V)\>},
\end{equation*}
induced by identifying a product $ah^i$, where $a\in\HH^*(Y)$, with
the class $(\pi^*a)h^i\in\HH^*(S(V)/S^1)$. 
\end{proposition}

\begin{remark}\label{rem:define-h}
Suppose $S^1$ acts with weight one on the fibres, so that we have a
natural isomorphism $S(V)/S^1\iso\P(V)$. Then the two definitions of
the class $h$ agree: the classical definition, as the first Chern
class of the dual of the tautological line bundle over $\P(V)$, and
the definition in the above proposition.  (The above definition of $h$
is equivalent to defining $h$ as the first Chern class of the
associated orbifold line bundle $S(V)\ti_{S^1}\C_{(1)}\to S(V)/S^1$,
where $\C_{(1)}$ denotes $\C$ with the weight one action of $S^1$. In
the classical case, it is easy to show that this associated line
bundle is isomorphic to the dual of the tautological line bundle.)
\end{remark}

\begin{proof}[Proof of Proposition \ref{pr:c-wtd-p-b}]
This proof comprises two steps. We first identify the weighted Chern
classes of $(V,S^1)$ as certain coefficients of an equivariant Euler class. 
We then show how this equivariant Euler class appears in a standard
long exact sequence, and how the properties of this long exact
sequence give us the proposition.

\paragraph{Step 1: Relating $c^w(V)$ to an equivariant Euler class.}
The $S^1$-equivariant bundle $V\to Y$ has an $S^1$-equivariant Euler
class
\begin{equation*}
e_{S^1}(V) \in \HH^{*}_{S^1}(Y) \iso \HH^{*}(Y)\ot
\HH^*(BS^1),
\end{equation*}
which we claim is given by 
\begin{equation}\label{eq:euler-chern}
e_{S^1}(V) = c_0^w(V) u^r + c_1^w(V) u^{r-1} + \h.. +
c_r^w(V),
\end{equation}
where $u\in\HH^2(BS^1)$ denotes the positive integral
generator. (We briefly recall the definition of the equivariant Euler
class. The equivariant cohomology of $Y$ is defined to be the regular
cohomology of the homotopy quotient $Y_{S^1}=(Y\ti ES^1)/S^1$. An
equivariant vector bundle $V\to Y$ pulls back to an equivariant vector
bundle over $Y\ti ES^1$, and by the quotient construction 
induces a regular vector bundle
over $Y_{S^1}$; the equivariant Euler class of $V$ is defined to be
the regular Euler class of this induced bundle.)

To show the above relationship between $e_{S^1}(V)$ and $c^w(V)$, 
we first show that it holds for line bundles. We then appeal to the
splitting principle to extend this to vector bundles.

Suppose $L\to Y$ is a
complex line bundle, possessing an action of $S^1$ covering a trivial
action on $Y$. Let $i\in\Z$ equal the weight of the action of
$S^1$ on the fibres of $L$. Let $L_{(0)}\to Y$ denote the same
line bundle, but with a trivial action of $S^1$, and let
$\ul\C_{(i)}\to Y$ denote the
trivial line bundle with a weight-$i$ action of $S^1$. Then 
\begin{equation*}
L \iso L_{(0)} \ot \ul\C_{(i)}
\end{equation*}
(as $S^1$-equivariant line bundles). Hence, since Euler classes add
when we tensor line bundles,
\begin{equation*}
\begin{array}{rl}
e_{S^1}(L) &= e_{S^1}(L_{(0)}) + e_{S^1}(\ul\C_{(i)}) \\
&= c_1(L) + i u\\
&= c_1^w(L) + c_0^w(L)u.\\
\end{array}
\end{equation*}
This proves our claim (equation \eqref{eq:euler-chern}) for line
bundles, and the general case follows from the splitting principle,
together with the observation that both sides of equation
\eqref{eq:euler-chern} are multiplicative with respect to direct sum
of the vector bundles we are considering.

\paragraph{Step 2: The map $\HH^{*}_{S^1}(Y) \to \HH^{*}(S(V)/S^1)$.}

Let $p$ and $\pi$ denote the maps
\begin{equation*}
\xymatrix{S(V) \ar[rr]^{/S^1} \ar[rd]_p && 
  S(V)/S^1 \ar[ld]^{\pi} \\
  & Y}
\end{equation*}
and let $/S^1$ denote the natural identification in equivariant
cohomology
$\HH^*_{S^1}(S(V)) \xra[\iso]{/S^1} \HH^*(S(V)/S^1)$.

Then, by naturality of this isomorphism, together with the definition
of $h$, we have
\begin{equation*}
(p^*(au^i))/S^1 = (\pi^*a)h^i
\end{equation*}
for any $a\in\HH^*(Y)$. 

But the natural map $(p^*\c.)/S^1$
fits into a short exact sequence of $\HH^*(Y)$-modules
\begin{equation}\label{ses}
\xymatrix{
0 \ar[r] &\<e_{S^1}(V)\> \ar@{^{(}->}[r] &\HH^{*}_{S^1}(Y) 
\ar@{>>}[rr]^{(p^*\c.)/S^1} &&\HH^{*}(S(V)/S^1) \ar[r] &0,\\
}
\end{equation}
where $\<e_{S^1}(V)\>\subset \HH^{*}_{S^1}(Y)$ denotes the
ideal generated by $e_{S^1}(V)$. 

These properties follow from the existence of  the long exact sequence in
equivariant cohomology for the
pair $(V,S(V))$, together with the following identifications:
\begin{equation*}
\xymatrix{
\h.. \ar[r] &\HH^{*}_{S^1}(V,S(V)) \ar[r] &\HH^{*}_{S^1}(V)
\ar[d]_\iso \ar[r] &\HH^{*}_{S^1}(S(V)) \ar[d]_\iso^{/S^1} \ar[r] &\h..\\
&\HH^{*}_{S^1}(Y) \ar[u]^\iso_{\smile \Phi} \ar[r]^{\smile e_{S^1}(V)} 
&\HH^{*}_{S^1}(Y) \ar[ru]^{p^*}
\ar[r] &\HH^{*}(S(V)/S^1) \\
}
\end{equation*}
Here the leftmost identification (denoted $\smile \Phi$) is the Thom
isomorphism in equivariant cohomology, with $\Phi$ the Thom class (see
\cite[section 2]{ati-bot:mom-map} for more on this identification); the next
identification is induced by restriction to the zero-section of $V$, 
and is an isomorphism because of the homotopy equivalence between $V$
and $Y$. The restriction of the Thom class $\Phi$ to the zero section
equals the equivariant Euler class $e_{S^1}(V)$, and hence the
composition of the Thom isomorphism with the restriction is given by
multiplication by the equivariant Euler class on $\HH^*_{S^1}(Y)$. 
The remaining maps are easily identified as labelled. 
Finally, using our explicit identification of the Euler class
$e_{S^1}(V)$, we see that multiplication by this Euler class is 
injective, and thus the sequence is short exact.

Hence we have 
\begin{equation*}
\HH^*(S(V)/S^1) \iso \frac{\HH^{*}_{S^1}(Y)}{\<e_{S^1}(V)\>},
\end{equation*}
and, substituting our formula for $e_{S^1}(V)$, we have proven the proposition.
\end{proof}

\subsection{Weighted Segre classes and the integration formula}

We now prove a formula which calculates integrals over the fibres of
the bundle $S(V)/S^1\xra{\pi} Y$. This formula involves the `weighted
Segre classes' of the pair $(V,S^1)$, which we define. We must also define an
orientation of the fibres of $\pi$, so that integration over the fibre
is well-defined.  In the case that $S(V)/S^1$ can be naturally
identified with a weighted projective bundle (i.e.\ if the weights of
the $S^1$-action are all positive)
this orientation agrees with the standard orientation induced by the
complex structure on the fibres.

\begin{definition}\label{def:wtd-s-c}
Let $V\to Y$ be a complex vector bundle with an action of $S^1$ as
above. The condition that the set of points fixed by the action equals
the zero section is equivalent to the condition that no subbundle of
$V$ be acted on with weight zero. It follows that the total weighted Chern
class of $(V,S^1)$ is invertible in the rational cohomology ring of $Y$
(since the degree-zero component is nonzero), and  we define the
\sdef{weighted Segre class} to be its multiplicative inverse:
\begin{equation*}
s^w(V) c^w(V) = 1.
\end{equation*}
\end{definition}

\begin{definition}\label{def:ort-app}
Given any point $y\in Y$, let  $V_y$ denote the fibre of $V$ over the point
$y$. Then, for any $v\in S(V_y)$, we have the
isomorphism
\begin{equation*}
T_{S^1\c.v}(S(V_y)/S^1) \oplus \R^+\c.v \oplus \la{s}^1 \iso V_y,
\end{equation*}
where $\R^+\c.v \subset V_y$ denotes the ray from the origin through
$v$, and $\la{s}^1$ is the Lie algebra of $S^1$, identified with $\R$
in the standard way.  We define the orientation of $S(V_y)/S^1$ to be
that orientation which is compatible with the above isomorphism
together with the given orientations of $\R^+$, $\la{s}^1$, and
$V_y$ (where $V_y$ has the standard orientation defined by its complex
structure, as in equation \eqref{eq:complex-ort}).
\end{definition}

\begin{proposition} \label{pr:segre-int}
Let $Y$ be connected and $V\to Y$ be a complex vector bundle with an
action of $S^1$ as above. Consider the bundle $S(V)/S^1\xra{\pi}Y$,
and define $h\in\HH^2(S(V)/S^1)$ to be the first Chern class of the
principal orbifold bundle $S(V)\to S(V)/S^1$ as in proposition
\ref{pr:c-wtd-p-b} above.  Then, for any $a\in\HH^*(Y)$,
\begin{equation*}
\pi_*\((\pi^*a)\smile h^i\) = \left\{ 
\begin{array}{ll}
0 & i<\rank(V)-1,\\
k a\smile s^w_{i-\rank(V)+1}(V),\qquad & i\ge \rank(V)-1.
\end{array}
\right.
\end{equation*}
Here $\pi_*$ denotes integration over the fibre with respect to the
orientation defined  above, and $k$ is the
greatest common divisor of the absolute values of the weights
appearing in the $S^1$ action on the fibres of $V$. 
\end{proposition}

\begin{proof}

This proof consists of two steps. In step 1 we relate the rational
fundamental class of the fibres with the fundamental class of complex
projective space. Then, in step 2, we use the formula
from proposition \ref{pr:c-wtd-p-b} above.

\paragraph{Step 1: The rational fundamental class of the fibres of 
$S(V)/S^1\to Y$.}

Given $y\in Y$, let $V_y$ denote the fibre of $V$ over the point
$y$. Then $S^1$ acts on $V_y$, and we can make an $S^1$-equivariant
identification 
\begin{equation*}
V_y \iso \C^r_{(i_1,i_2,\h..,i_r)},
\end{equation*}
where $\C^r_{(i_1,i_2,\h..,i_r)}$ denotes $\C^r$ with the
weight-$(i_1,i_2,\h..,i_r)$ action of $S^1$ (that is, $\lb\in
S^1\subset\C^\ti$ acts by
$\lb\c.(z_1,\h..,z_r) = (\lb^{i_1}z_1,\h..,\lb^{i_r}z_r)$.)
Moreover, we can arrange that $i_1,\h..,i_n < 0$ and $i_{n+1},\h..,i_r
> 0$. Then the map
\begin{equation*}
\begin{array}{rl}
\tilde{\vp} : \C^r_{(1,1,\h..,1)} &\to V_y = \C^r_{(i_1,i_2,\h..,i_r)} \\
   (z_1,\h..,z_r) &\mapsto
   ({\bar{z}_1}^{|i_1|},\h..,{\bar{z}_n}^{|i_n|},{z_{n+1}}^{i_{n+1}},{z_r}^{i_r})\\
\end{array}
\end{equation*}
is smooth and intertwines the $S^1$-actions.

There is an obvious $S^1$-invariant metric on $V_y
=\C^r_{(i_1,i_2,\h..,i_r)}$ such that $\tilde{\vp}$ maps the standard
unit sphere in $\C^r_{(1,1,\h..,1)}$ to the unit sphere in $V_y$.

Hence $\tilde{\vp}$ descends to a map
\begin{equation*}
\vp  : S(\C^r_{(1,1,\h..,1)})/S^1 = \cp^{r-1} \to 
S(V_y)/S^1.
\end{equation*}

We can now relate the rational fundamental class of $S(V_y)/S^1$ to
the fundamental class of $\cp^{r-1}$ by calculating the oriented
degree of $\vp$ (that is, the topological degree of $\vp$, multiplied by
$\pm 1$ according to whether $\vp$ preserves or reverses orientation).

We easily see that the oriented degree of $\tilde{\vp}$ equals 
$\prod_{j=1}^r i_j = c_0^w(V)$ (and this also equals the oriented
degree of the restriction of $\tilde{\vp}$ to the unit sphere).
To calculate the degree of $\vp$, we must divide this number
by the degree with which a generic $S^1$-orbit in
$S(\C^r_{(1,1,\h..,1)})$ covers its image. It is easy to see that this
degree equals the greatest common divisor of the absolute values of the
$i_j$. Hence, setting
\begin{equation*}
k:= \gcd (|i_1|,|i_2|,\h..,|i_r|)
\end{equation*}
then the oriented degree of $\vp$ is given by
\begin{equation}
\deg(\vp) = k\inv c_0^w(V).
\end{equation}

Now consider the maps 
\begin{equation*}
\xymatrix{
{\cp^{r-1}} \ar[r]^\vp 
&S(V_y)/S^1 \ar@{^{(}->}[r]^\psi \ar[d]_{\pi'} &S(V)/S^1 \ar[d]^\pi \\
&y \ar@{^{(}->}[r] &Y \\
}
\end{equation*}
We have defined the class $h\in\HH^*(S(V)/S^1)$ to be the first Chern
class of the orbifold $S^1$-bundle $S(V)\to S(V)/S^1$. (Or
equivalently, $h$ is the first Chern class of the associated orbifold
line bundle $S(V)\ti_{S^1} \C_{(1)} \to S(V)/S^1$, where $\C_{(1)}$
denotes $\C$ with the weight-one action of $S^1$.) By naturality of
this definition, we
see that $h$ pulls back to the integral generator of the cohomology of
$\cp^{r-1}$, so that
\begin{equation*}
\ig{(\vp^*\psi^*h)^{r-1}}{\cp^{r-1}} =1.
\end{equation*}
Using the degree of $\vp$, we
thus have
\begin{equation*}
\pi'_*((\psi^*h)^{r-1}) = k c_0^w(V)\inv,
\end{equation*}
and hence, since integration over the fibre commutes with restriction,
and the result is a
degree-zero cohomology class, we have
\begin{equation}\label{int-h-r-1}
\pi_*(h^{r-1}) = k c_0^w(V)\inv \qquad \in\HH^0_G(Y).
\end{equation}
(Of course $\pi_*(h^i)=0$ if $i<r-1$, for degree reasons.)

\paragraph{Step 2: Using the relation in cohomology to extend this
formula to all powers of $h$.}

We now calculate
\begin{equation}\label{product}
\pi_*\((\pi^*c^w(V))\smile(1 + h + h^2 + h^3 + \h..)\).
\end{equation}
Since $\pi_*$ lowers degree by $2r-2$, we only need to consider terms
in the product $(\pi^*c^w(V))\smile(1 + h + h^2 + h^3 + \h..)$
of degree $2r-2$ and greater. The degree $2r-2$ term
is
\begin{equation*}
\pi^*c^w_{r-1}(V) + \pi^*c^w_{r-2}(V)h + \h.. + \pi^*c^w_0(V)h^{r-1},
\end{equation*}
and applying $\pi_*$ to this term gives the coefficient of $h^{r-1}$,
multiplied by $\pi_*(h^{r-1})$ (which we have calculated in equation
\eqref{int-h-r-1} above). (We are using the fact 
that $\pi_*$ is a  homomorphism of $\HH^*(Y)$-modules.)
Hence, the integral over the fibre of the
degree $2r-2$ term of the product \eqref{product} equals $k\in\HH^0(Y)$. 

The degree $2r$ term of the product \eqref{product} is
\begin{equation*}
\pi^*c^w_{r}(V) + \pi^*c^w_{r-1}(V)h + \h.. + \pi^*c^w_0(V)h^r.
\end{equation*}
Comparing with our explicit identification of $e_{S^1}(V)$ in equation
\eqref{eq:euler-chern} above, we see that this term is exactly the
class $p^*\(e_{S^1}(V)\)/S^1$, and, using the short exact sequence
\eqref{ses}, this term vanishes.  Similarly, the degree $2(r+j)$ term
equals $p^*\(u^j e_{S^1}(V)\)/S^1$ and hence also vanishes.
Thus we have
\begin{equation*}
\pi_*\((\pi^*c^w(V))\smile(1 + h + h^2 + h^3 + \h..)\) = k
\end{equation*}
and hence (using the module-homomorphism property of $\pi_*$)
\begin{equation*}
\pi_*(1 + h + h^2 + h^3 + \h..) = k c^w(V)\inv = k s^w(V).
\end{equation*}
The proposition now follows by identifying terms by degree.
\end{proof}

\subsection{Equivariant weighted Segre classes, and the equivariant
integration formula}

Suppose $V\to Y$ is a vector bundle, with an action of $S^1$ as above,
and suppose moreover that an auxilliary group $G$ acts on $V$ and $Y$,
commuting with the projection and with the action of $S^1$. Then we
can generalize the definition of weighted Chern classes and weighted
Segre classes to the $G$-equivariant case as follows.

Recall that the homotopy quotient construction replaces a $G$-space
$Y$ with the space $Y_G:=EG \ti_G Y$, and the equivariant cohomology
of $Y$ is defined to be the ordinary cohomology of $Y_G$.  Given a
$G$-equivariant vector bundle $V\to Y$ the same construction gives a
vector bundle $V_G \to Y_G$ (this is explained by Atiyah and Bott in
\cite[section 2: equation (2.1) and remark (1)]{ati-bot:mom-map}), and
the $G$-equivariant characteristic classes of $V$ are then taken to be
the ordinary characteristic classes of $V_G \to Y_G$, which thus take
values in the $G$-equivariant cohomology of $Y$.

\begin{definition}\label{def:ewtd}
In our case, in the presence of a commuting $S^1$-action, the bundle
$V_G\to Y_G$ has an induced $S^1$-action, and we define the
\sdef{$G$-equivariant weighted Chern classes} and the
\sdef{$G$-equivariant weighted Segre classes} of the pair $(V,S^1)$ to
be the weighted Chern classes and weighted Segre classes of the pair
$(V_G\to Y_G,S^1)$ (definitions \ref{def:wtd-c-c} and \ref{def:wtd-s-c}).
\end{definition}

Applying
the cohomology formula (Proposition \ref{pr:c-wtd-p-b}) to
the bundle $V_G\to Y_G$ and making the obvious identifications, we
thus get the equivariant version of the cohomology formula:
\begin{equation}\label{eq-coh}
\HH_G^*(S(V)/S^1) \iso \frac{\HH_G^*(Y)[h]}%
{\<c_0^w(V)h^r + c_1^w(V)h^{r-1} + \h.. + c_r^w(V)\>},
\end{equation}
where $c_i^W(V)$ now denotes the $i$-th $G$-equivariant weighted Chern
class of $(V,S^1)$.
Similarly, the integration formula (Proposition \ref{pr:segre-int}) 
applied to $V_G\to Y_G$ gives the (formally identical) 
equivariant formula:
\begin{equation}\label{eq-int}
\pi_*\((\pi^*a)\smile h^i\) = \left\{ 
\begin{array}{ll}
0 & i<\rank(V)-1,\\
k a\smile s^w_{i-\rank(V)+1}(V),\qquad & i\ge \rank(V)-1,
\end{array}
\right.
\end{equation}
where $s_i^W(V)$ now denotes the $i$-th $G$-equivariant weighted Segre
class of $(V,S^1)$.

\section{Proof of the orientation lemma}\label{app:ort-bd}

We now give the proof of lemma \ref{lem:ort-bd}, which describes the
orientations of the boundary components of the
wall-crossing-cobordism.

\begin{lemma}[Lemma \ref{lem:ort-bd}]
Let the wall-crossing-cobordism $W/T$ be oriented as in definition
\ref{def:ort-w-c-c}. Then the induced boundary orientation of 
$\syq{X}{p_0}{T}$ is $ -(\om_{p_0}^k)$, and of
$\syq{X}{p_1}{T}$ is $\om_{p_1}^k$
(where $\om_{p_i}$ denote the respective induced symplectic forms), 
and the induced boundary orientation of each $P_{(H,q)}$ is equal to
the product orientation defined in \ref{def:ort-PHq} above.
\end{lemma}

\begin{proof}
Before beginning the proof proper, we fix three conventions which will
hold throughout the proof.

1. Most of the steps in this proof consist of exhibiting
isomorphisms of the form
\begin{equation*}
V_1 \oplus V_2 \iso V_3,
\end{equation*}
where the $V_i$ are vector spaces. For each such isomorphism, we will
be using orientations of two of the vector spaces to induce an
orientation on the third, in the obvious manner (explicitly: so that
concatenating oriented bases for $V_1$ and $V_2$ gives an oriented
basis for $V_3$).

2. When we decompose tangent spaces, we will assume without
explicit mention that these decompositions are orthogonal
decompositions relative to some choice of invariant metric; and it
will always be the case that the induced orientations are independent
of the choices made.

3. Finally, given a symplectic form (that is, a nondegenerate
$2$-form) on any vector space, then the \sdef{symplectic orientation}
of that vector space will mean the orientation defined by the top
power of the symplectic form. 

We break the proof up into two steps. In step 1 we assume that $T$ is
$1$-dimensional, and in step 2 we reduce the general case to to the
case treated in step 1.

\paragraph{Step 1: Assuming $T$ is $1$-dimensional.}

If $T$ is $1$-dimensional, then $Z$ is just a closed subinterval of
$\t$, bounded by $p_0$ and $p_1$. The orientation of $Z$ (definition
\ref{def:w-c-d}) and of each wall-crossing subgroup (which is just $T$
itself) orients $\t$ so that $p_0<p_1$. 
We orient $\t^*$ compatibly (that is, so that the duality pairing
between a positive vector in $\t$ and a positive vector in $\t^*$ is
positive). 

Then, restating definition \ref{def:ort-w-c-c} in this case,  
we have oriented the wall-crossing-cobordism $W/T$ via the
isomorphism 
\begin{equation}\label{eq:bd1}
T_{[x]}(W/T) \iso T_{[x]}\syq{X}{p}{T} \oplus \t^*
\end{equation}
for each $x\in W$ (with $p=\mu(x)$), 
relative to the orientation of $\td$ described
above and the symplectic orientation of $\syq{X}{p}{T}$.

The key calculation in step 1 is to 
compare this orientation of $W/T$ with the symplectic
orientation of $X$. Now at each $x\in W$, we can decompose $T_xW=T_xX$
into the orbit direction and its orthgonal complement. Using the
natural identifications, then we claim that
\begin{equation}\label{eq:bd2}
T_xX = T_xW\iso T_{[x]}(W/T) \oplus \ol{\t},
\end{equation}
is orientation-preserving, where $\ol{\t}$ denotes $\t$ with the
opposite orientation. Using equation \eqref{eq:bd1} above, this is
equivalent to showing that the isomorphism
\begin{equation}\label{eq:bd3}
T_xX \iso  T_{[x]}\syq{X}{p}{T} \oplus \t^* \oplus \ol{\t}
\end{equation}
is orientation-preserving. Here we think of the spaces on the right as
subspaces of $T_xX$. Explicitly, we let 
\begin{equation}\label{eq:bd4}
\begin{array}{rl}
h:T_{[x]}\syq{X}{p}{T} &\hra T_{[x]}X, \\
i:\td &\hra T_xX, \quad\text{and} \\
j:\t &\hra T_xX \\
\end{array}
\end{equation}
denote these
identifications (so that $i\inv=\eval{d\mu}_{\im(i)}$ and $j$ is given
by the infinitesimal action of $T$ at $x$, and $h$ is the
identification of a complement to $j(\t)$ in $\mu\inv(p)$ with a slice
at $x$). Let $\psi\in\td$ and
$\xi\in\t$ be positive with respect to the orientations we have chosen
for $\td$ and $\t$ (so that $\xi$ is negative with respect to the
orientation of $\ol{\t}$).

Then our sign convention for the moment map condition (equation
\eqref{eq:moment-map}) implies that 
\begin{equation}\label{eq:bd5}
\begin{array}{rl}
\om(i(\psi),-j(\xi))
  &= \om(j(\xi),i(\psi)) \\
  &= \<d\mu(i(\psi)),\xi\> \\
  &= \<\psi,\xi\> \\
  & >0.
\end{array}
\end{equation}
This means that $(i(\psi),-j(\xi))$ is a positively oriented basis of
$i(\td)\oplus j(\ol{\t})\subset T_xX$  with respect to the restriction
of the symplectic form on $X$. 

Now, recall that the symplectic form on
$\syq{X}{p}{T}$ is induced by restricting the symplectic form on $X$ to
$\mu\inv(p)$, where it is degenerate in the orbit directions, and
hence descends to $\syq{X}{p}{T}$. In terms of the maps in
\eqref{eq:bd4} above, this means that the symplectic form on
$T_{[x]}\syq{X}{p}{T}$ agrees with the pullback, via $h$, of the
symplectic form on $T_xX$. 

Thus we have shown that the identification in equation \eqref{eq:bd3} is
orientation-preserving.  (Being completely explicit: if $(v_1,
\h.. ,v_k)$ is an oriented basis of $T_{[x]}\syq{X}{p}{T}$, then\\
$(v_1,\h..,\v_k,\psi,-\xi)$ is an oriented basis of the right-hand-side
of equation \eqref{eq:bd3}, and\\ $(h(v_1),\h..,h(v_k),i(\psi),j(-\xi))$
is an oriented basis of the left-hand-side.)

Having derived this alternative description of the orientation of
$W/T$, we can now calculate the induced orientations of the various
boundary components. Now $W/T$ is an odd-dimensional manifold, and
hence the induced boundary orientation is defined to be that
orientation of $\bd(W/T)$ which is compatible with the isomorphism
\begin{equation*}
T_{[x]}(W/T) \iso T_{[x]} \bd(W/T) \oplus \R\c.\outn
\end{equation*}
where $\outn$ is an outward-pointing normal vector.  (Here we are
using the convention that makes Stokes's theorem sign-free, as
explained in \cite[page 31]{bot-tu:dif-for}; for the boundary of an 
even-dimensional manifold we would need to use the inward-pointing normal
vector in the above equation.)
Combining this with equation \eqref{eq:bd2}, this means that we can calculate the
orientation of $\bd(W/T)$ via the isomorphism
\begin{equation*}
T_xX \iso T_{[x]} \bd(W/T) \oplus \R\c.\outn \oplus \ol{\t},
\end{equation*}
and the orientations of $\R$, $\ol{\t}$, and the symplectic
orientation of $T_xX$.

For the boundary component $\syq{X}{p_0}{T}$, the isomorphism
$\R\c.\outn\iso \ol{\td}$ is orientation preserving. Thus we use the
isomorphism 
\begin{equation*}
T_xX \iso T_{[x]}\syq{X}{p_0}{T} \oplus \ol{\td} \oplus \ol{\t} 
\end{equation*}
together with the orientations of $\ol{\td}, \ol{\t}$, and $T_xX$ to
orient this boundary component. Using the reasoning in equations
\eqref{eq:bd5}, this gives
\begin{equation*}
-(\om_{p_0})^k.
\end{equation*}
A similar argument (except with $\td$ replacing $\ol{\td}$) shows that
the induced boundary orientation of $\syq{X}{p_1}{T}$ is equal to its
symplectic orientation:
\begin{equation*}
\om_{p_1}^k.
\end{equation*}

Finally we come to $P_{(H,q)}$, which in our case is $P_{(T,q)}$,
since $H=T$. To
conform with the notation of definition \ref{def:ort-PHq}, let $x\in
X^T=X_{(T,q)}$ and $v\in S(\nu_x X^T)$, where we identify the point
$(x,v)$ with a point in $X$ via an equivariant exponential map. Then
$P_{(T,q)}$ is oriented by the isomorphism
\begin{equation*}
T_{(x,v)}X \iso T_{[x,v]}P_{(T,q)} \oplus \R\c.(-v) \oplus \ol{\t}
\end{equation*}
(with the symplectic orientation of $T_{(x,v)}X$). Since $X^T$ is a
symplectic submanifold, we can decompose the symplectic form according
to the isomorphism
\begin{equation*}
T_{(x,v)}X \iso T_x X^T \oplus \nu_x X^T,
\end{equation*}
so that this can be viewed as an orientation-preserving isomorphism
with respect to the induced symplectic forms on all three
spaces. Then, since $\R\c.(-v)\oplus\ol{\t}$ gives the same
orientation as $\R\c.v\oplus\t$, we have recovered the orientation
of definitions~\ref{def:ort-V} and \ref{def:ort-PHq}.

\paragraph{Step 2: Reducing the general case to the case of step 1.}

In step 1 we assumed that the torus $T$ was $1$-dimensional. It is
easy to reduce the general case to the case of step 1, as follows.

We first observe that the orientation of $W/T$ is locally defined (in
terms of a codimension-$1$ foliation by symplectic orbifolds).  In
order to reduce the general orientation calculation, we need only
consider the wall-crossing-cobordism in a neighbourhood of a boundary
component. We will describe the construction for the boundary
components $P_{(H,q)}$; the case of the components $\syq{X}{p_i}{T}$
is analogous.

We fix attention on a single wall-crossing, with 
wall-crossing data $(H,q)$, and associated boundary $P_{(H,q)}$.
Choose $T'\subset T$ so that $T=T'\ti H$. Then we have the following
inclusions and associated dual projections
\begin{equation*}
\begin{array}{rlrl}
 & \t'\hra\t \qquad\qquad\qquad & & \la{h}\hra\t \\
\phi: &\td\onto\t'{}^* & \psi: & \td\onto\la{h}\\
\end{array}
\end{equation*}
We define $q':=\phi(q)$, and $\mu':=\phi\o\mu$, so that $\mu'$ is a
moment map for the action of $T'$. 

Now suppose that, in some neighbourhood of $q$, $Z$ is parallel to
$\t'{}^\perp$ (this can easily be arranged by deforming $Z$). Then
\begin{equation*}
\mu\inv(Z) = \mu'{}\inv(q')
\end{equation*}
in a neighbourhood of $\mu\inv(q)$.
Now the $T$-action on $X$ descends to an action of $H$ on
$\syq{X}{q'}{T'}$, with moment map given by the restriction of
$\psi\o\mu$. 

It is now easy to see that, in a neighbourhood of $P_{H,q}$, 
the wall-crossing-cobordism $W(X,T,\mu,Z)$
constructed from the data $X,T,\mu,Z$ coincides with the
wall-crossing-cobordism $W(\syq{X}{q'}{T'},H,\psi\o\mu,\psi(Z))$. These
are foliated by the same symplectic suborbifolds, since
\begin{equation*}
\syq{X}{p}{T} \iso \syq{(\syq{X}{q'}{T'})}{\psi(p)}{H}
\end{equation*}
is an isomorphism of symplectic stratified spaces.

Since $Z$ is transverse to $\mu$ at $q$, it follows that
$\syq{X}{q'}{T'}$ is a symplectic orbifold in a neighbourhood of
$\mu\inv(q)$, with a Hamiltonian action of the $1$-dimensional torus
$H$, and we have thus reduced our calculation to the case of step 1.
\end{proof}

\bibliographystyle{alpha}

\end{document}